\documentclass[10pt]{amsart}
\usepackage[top=1in,bottom=1.25in,left=1.25in,right=1.25in]{geometry}
\usepackage{xr}
\usepackage{hyperref}
\usepackage{graphicx}
\usepackage[foot]{amsaddr}
\usepackage{amsmath}
\usepackage{amssymb}
\usepackage{amsthm}
\usepackage{xcolor}
\usepackage{float}
\usepackage{indentfirst}
\usepackage{placeins}
\usepackage{booktabs}
\usepackage{mathtools}
\usepackage{cleveref}
\usepackage{subcaption}
\usepackage{tikz}

\crefname{figure}{Figure}{Figures}
\crefname{table}{Table}{Figures}
\newcommand{\C}{\mathbb{C}}
\newcommand{\R}{\mathbb{R}}
\newcommand{\A}{\textup{~\AA}}
\newcommand{\cB}{\mathcal{B}}
\newcommand{\bS}{\mathbb{S}}
\newcommand{\cO}{\mathcal{O}}

\newcommand{\cN}{\mathcal{N}}
\newcommand{\cS}{\mathcal{S}}
\newcommand{\cJ}{\mathcal{J}}
\newcommand{\iid}{\overset{i.i.d.}{\sim}}

\newcommand{\cM}{\mathcal{M}}
\newcommand{\Mest}{\overline M}
\newcommand{\cI}{\mathcal{I}}
\newcommand{\Int}{\int_{\mathcal{SO}(3)}}
\newtheorem*{remark}{Remark}

\title[Subspace method of moments]{Subspace method of moments for  ab initio 3-D single-particle Cryo-EM reconstruction}

\author[J. Hoskins]{Jeremy Hoskins}
\address{Department of Statistics and CCAM, 
    University of Chicago, Chicago, IL 60637 USA.}
\email{jeremyhoskins@uchicago.edu}
\author[Y. Khoo]{Yuehaw Khoo}
\address{Department of Statistics and CCAM, 
    University of Chicago, Chicago, IL 60637 USA.}
\email{ykhoo@uchicago.edu}
\author[O. Mickelin]{Oscar Mickelin}
\address{Program in Applied and Computational Mathematics, Princeton University, Princeton,  NJ 08544 USA.}
\email{hm6655@princeton.edu}
\author[A. Singer]{Amit Singer}
\address{Department of Mathematics and Program in Applied and Computational Mathematics, Princeton University, Princeton, NJ 08544 USA.}
\email{amits@math.princeton.edu}
\author[Y. Wang]{Yuguan Wang}
\address{Department of Statistics, University of Chicago,  Chicago, IL 60637 USA.}
\email{yuguanw@uchicago.edu}

\begin{document}

\begin{abstract}
Cryo-electron microscopy (cryo-EM) is a widely used technique for recovering the 3-D structure of biological molecules from a large number of experimentally generated noisy 2-D tomographic projection images of the 3-D structure, taken from unknown viewing angles. Through computationally intensive algorithms, these observed images are processed to reconstruct the 3-D structures. Many popular computational methods rely on estimating the unknown angles as part of the reconstruction process, which becomes particularly challenging at low signal-to-noise ratios. The method of moments (MoM) offers an alternative approach that circumvents the estimation of viewing orientations of individual projection images by instead estimating the underlying distribution of the viewing angles, and is robust to noise given sufficiently many images. However, the method of moments typically entails computing higher-order moments of the projection images, incurring significant  computational and memory costs. To mitigate this, we propose a new approach called the subspace method of moments (SubspaceMoM), which compresses the first three moments using data-driven low-rank tensor techniques as well as expansion into a suitable function basis. The compressed moments can be efficiently computed from the set of projection images using numerical quadrature and can be employed to jointly reconstruct the 3-D structure and the distribution of viewing orientations. We illustrate the practical applicability of SubspaceMoM through numerical experiments using up to the third-order moment  on  synthetic datasets with a  simplified cryo-EM image formation model, which significantly improves the reconstruction resolution  compared to previous MoM approaches.
\end{abstract}

\maketitle

\section{Introduction}
\label{sec:introduction}
\begin{figure}[!ht]
 \centering
 \begin{tikzpicture}[scale=0.99]
  \node[inner sep=0] at (-17,-5.5) {\includegraphics[width=0.49\textwidth]{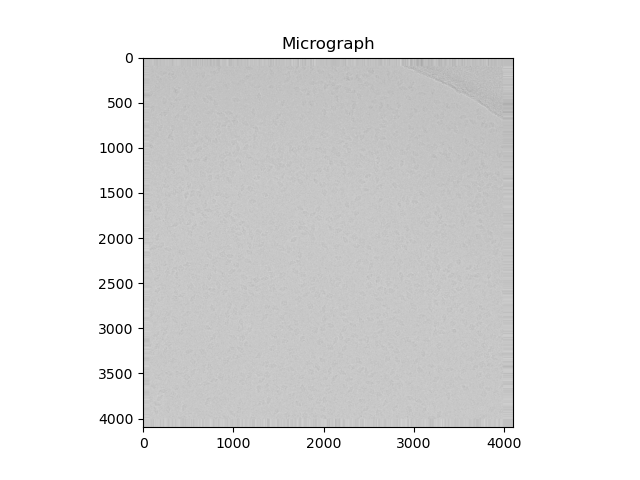}}; 
  \node[inner sep=0] at (-9,-5.5) {\includegraphics[width=0.49\textwidth]{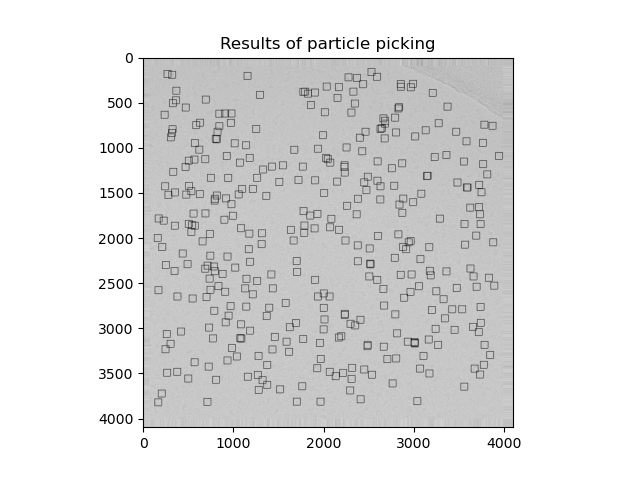}}; 
  \node[inner sep=0] at (-17,-8.5) {\footnotesize (A)};
  \node[inner sep=0] at (-9,-8.5) {\footnotesize (B)};
 \end{tikzpicture}
\caption{(A) A cryo-EM micrograph from the Beta-galactosidase Falcon-II  dataset~\cite{SCHERES2015114}. (B) The result of the particle picking procedure using the APPLE picker~\cite{HEIMOWITZ2018215} implemented in the ASPIRE software \cite{aspire}. }
 \label{fig:particle_picking}
\end{figure}

Cryo-EM \cite{Shen2018} is a powerful imaging technique that enables the study of biological molecules, including proteins and nucleic acids, at near-atomic resolution. This technique involves freezing many copies of a given molecule in a thin layer of ice. The sample is then imaged using an electron microscope. The resulting projection image, termed a \emph{micrograph}, consists of many two-dimensional projections of the three-dimensional molecule.  Due to thermal motion, the copies of the molecules in the sample rotate before freezing. The projections in the micrograph therefore consist of tomographic projections of randomly rotated copies of the molecule. A \textit{particle picking} procedure  \cite{Bepler2019, HEIMOWITZ2018215} is then employed to extract the individual projections from the micrograph.  An example of a micrograph and the result of particle picking are shown in~\cref{fig:particle_picking}.  Once extracted, the goal of cryo-EM is to reconstruct the three-dimensional structure of the molecule from the projection images. A significant impediment to the practical application of cryo-EM is the large amount of noise present in the images. This noise, arising due to radiation damage of the sample, can easily obscure critical features, making accurate reconstruction difficult.

One popular method for cryo-EM reconstruction is based on maximum likelihood estimation \cite{SIGWORTH2010263,10.1214/23-AOS2292}, typically implemented via the Expectation–Maximization (EM) algorithm and referred to as 3-D iterative refinement~\cite{SCHERES2012519}.  This method refines a 3-D model by iteratively aligning it to the observed 2-D projection images until it closely matches the experimental data. The EM approach has been extremely successful in determining high resolution structures of large macromolecules. However, EM fails to produce meaningful reconstruction for molecules with molecular weights below $\sim$40 kDa. This limitation is primarily due to the low signal-to-noise ratio (SNR), which renders orientation assignment unreliable.   Furthermore, cryo-EM reconstruction is inherently a high-dimensional inverse problem, and typical datasets consist of hundreds of thousands of images, each of which must be utilized in every iteration of the refinement process.  As a result, the computational cost is substantial. In fact, the number of projection images required for accurate reconstruction is inversely proportional to the SNR.  The likelihood optimization problem is also highly nonconvex, further complicating the reconstruction. To address these challenges, techniques such as stochastic optimization~\cite{punjani2017cryosparc} have been proposed. Additionally, \emph{ab initio} modeling \cite{Reboul2018SingleparticleCA,doi:10.1137/090767777} is often employed to produce an initial structure sufficiently close to the ground truth to enable successful iterative refinement.

The method of moments (MoM) approach for single particle reconstruction, proposed over forty years ago by Zvi Kam \cite{KAM198015}, is a promising approach for \textit{ab initio} modeling. One particularly attractive feature of it is that it requires only low-order moments of the 2-D projection images of the volume. Given a large number of images, together with appropriate debiasing tools, the moments can be estimated quite accurately, making MoM less sensitive to noise compared to other \textit{ab initio} approaches. Initially, Kam  assumed a uniform distribution of the molecule orientations within the samples. However, he found that the second moment alone is insufficient and suggested using certain slices of the third-order and fourth-order moments to improve performance. In general, the third-order moment is required in the case of uniform distribution of orientations  \cite{BANDEIRA2023236}. However, \textit{ab initio} modeling with only the second-order moment is still possible under uniformity when there are one or two noiseless projections provided \cite{8363873,doi:10.1137/22M1498218}, which can be found by denoising the noisy images \cite{SHI2022107018,9506435,BHAMRE201672,marshall_mickelin_shi_singer_2023}. More recent works \cite{pnas.2216507120,bendory2023sample} also showed that reconstruction is possible using the second-order moment when the molecule has a sufficiently sparse representation. Unfortunately, the uniform assumption is frequently invalid in practice, and the empirical orientations of the molecules can be quite non-uniform. Unlike the widespread belief that uniform distribution of orientation is optimal for reconstruction~\cite{BALDWIN2020160}, non-uniformity is actually helpful, and it has been shown that the first two moments become sufficient for recovering a suitable approximation of the molecular structure if the distribution is known and very non-uniform~\cite{Sharon_2020}. It was also found in \cite{Sharon_2020} that the distribution can be estimated jointly with the 3-D structure from the first two moments, but the resolution of the estimated structure is limited.  For the iterative refinement approach, one can also sidestep the orientation estimation by matching the data in a distributional sense~\cite{9483649}. For a simpler but closely related problem called multi-reference alignment (MRA) \cite{doi:10.1137/18M1214317,doi:10.1137/20M1354994,Fan2020LikelihoodLA} (that often serves as a mathematical abstraction of cryo-EM without the projections), the third-order moment is required in the case of the uniform distribution, but the first two moments suffice for generic non-uniform distributions~\cite{Abbe2017MultireferenceAI}. Additionally, higher order moments were found to be important for MRA when using the generalized method of moments~\cite{Abas2021TheGM}, 3-D reconstruction from images generated by X-ray free electron lasers (XFELs)~\cite{Zhao:zf5022} and reconstructing high-dimensional noisy curves~\cite{lo2024methodmomentsestimationnoisy}. Moreover, it is shown in \cite{BENDORY2022391} that a signal in $\C^N$ can be uniquely determined by a few measurements of its \emph{bispectrum} and \emph{trispectrum}, which are some special entries of the third-order and fourth-order moments. Thus, one can hope that incorporating the third-order moment would significantly improve the quality of cryo-EM reconstruction for the case of an unknown non-uniform distribution of viewing orientations, which is the typical case for experimental cryo-EM data. In general, the third-order moment is not practically useful because it is too expensive to compute and too large to store. The focus of this paper is to develop an efficient approach to compress and represent the third-order moment, without sacrificing the reconstruction quality.

We propose a flexible computational framework, called the \emph{subspace method of moments (SubspaceMoM)}, which enables joint reconstruction of the 3-D molecular structure (often referred to as ``volume'' by practitioners) and the unknown non-uniform distribution from a compressed representation of the first three moments, which are the \textit{subspace moments}. Following \cite{Sharon_2020}, we assume the volume and the non-uniform distribution can be expanded using some basis functions. We project the  sample moments computed from the noisy images into low-dimensional subspace moments, which  are expressed as sums of simple rank-one components that can be evaluated efficiently. The coefficients of the volume and the non-uniform density are then reconstructed from the subspace moments through numerical optimization. We provide a particular choice of the basis functions, which serves as a baseline for SubspaceMoM. We also provide a way of finding subspaces by computing low-rank (tensor) decompositions of the second and  third moments using a randomized method.  We demonstrate the proposed method on simulated datasets that include the effect of the contrast transfer function (CTF) of the electron microscope~\cite{MINDELL2003334,ROHOU2015216}.  The numerical experiments show that we can successfully obtain \textit{ab initio} models from the first three subspace moments  when the distribution is non-uniform. We also show that for high-noise levels, the \emph{ab initio} models obtained by SubspaceMoM have better resolution than those obtained by stochastic optimization and can reliably initialize a high-resolution 3-D refinement algorithm.  Extending the numerical benchmarks to experimentally obtained projection images would require taking into account several factors that could potentially lead to inaccurately estimated moments. This includes the fact that the projection images can be imperfectly centered, that the particle picking procedure can potentially  select non-particles, that molecules are often flexible and exhibit conformational heterogeneity, and that detector response~\cite{VULOVIC201319} can affect the images.  We defer these experimental considerations to future work.   The Python implementation of our method is  publicly available at \href{https://github.com/wangyuguan/subspace_MoM}{https://github.com/wangyuguan/subspace\_MoM}.

\section{Mathematical Preliminaries}
\label{sec:Preliminaries}
In this section, we review some preliminaries necessary for describing our method. In Section \ref{sec:tensornotations}, we first introduce convenient notations for representing and manipulating tensors. Following this, in Section \ref{sec:randomrangefinding}, we introduce relevant techniques from randomized numerical linear algebra for matrix and tensor decompositions.  Then in Section \ref{sec:intro_cryo_em}, we introduce the cryo-EM model and the reconstruction problem considered in this paper. Finally, in Section \ref{sec:MoMforCryoEM}, we conclude with a brief introduction to the method of moments in cryo-EM.

 \subsection{Tensor notations}
\label{sec:tensornotations}
Let $M \in \C^{d_1\times d_2 \times \ldots \times d_n}$ be a tensor of order $n \ge 3$ where $d_1,\ldots,d_n$ are positive integers denoting the dimensions of the tensor. Let $U \in \C^{d_k\times r}$ be a matrix where $1\le k \le n$.  The \textit{mode-$k$ contraction} between the tensor $M$ and the matrix $U$ is defined as  
\begin{align}
\label{eqn:modeprod}
    (M \times_k U)_{i_1,\ldots,i_{k-1},j,i_{k+1},\ldots,i_n} = \sum_{i_k=1}^{d_k} M_{i_1,\ldots,i_{k-1},i_k,i_{k+1},\ldots,i_n} U_{i_k,j}
\end{align}
for $1\le i_1\le d_1,\ldots,1\le i_n \le d_n$ and $1\le j \le r$. The result of the contraction $$M \times_k U \in \C^{d_1 \times \ldots \times d_{k-1} \times r \times d_{k+1} \times \ldots \times d_n}$$ is also a tensor of order $n$. 

The \textit{mode-$k$ unfolding matrix} of  $M$ is defined to be a matrix $M_{[k]} \in \C^{d_k\times \prod_{i\not =k}^n d_i}$ obtained by ``flattening'' all but the $k$-th dimension. Specifically, the $(j,l)$-th entry of the matrix $M_{[k]}$ corresponds to the tensor entry $ M_{i_1,\ldots,i_k=j,\ldots,i_n}$, where the column index  $l$ is computed via the mapping:
\begin{align}
\label{eqn:unfolding}
\ell = 1 + \sum_{\substack{k' = 1, \ k' \neq k}}^{n} (i_{k'} - 1) \prod_{\substack{l' = 1, \ l' \neq k}}^{k'-1} d_{l'}.
\end{align}
 
We now introduce a standard low-rank representation format for tensors, which is called the  \textit{Tucker decomposition}~\cite{Tuck1966c}. The Tucker decomposition   of a tensor $M$ approximates  $M$  in terms of matrices $U^{(k)}\in \C^{d_k \times r_k}$ where $r_k\le d_k,  k=1,\ldots,n$  and a core tensor $M_c \in \C^{r_1\times r_2 \times \ldots \times r_n}$ by
\begin{align}
\label{eqn:tucker}
    M \approx M_c \times_1 U^{(1)} \times_2 U^{(2)} \times_3 \ldots \times_n U^{(n)}.
\end{align}
In this paper, we further require the matrices $U^{(k)},k=1,\ldots,n$ to have orthonormal columns. 
For each $k$, the matrix $U^{(k)}$ can be obtained by finding the range of the mode-$k$ unfolding matrix $M_{[k]}$. Specifically, the columns of $U^{(k)}$ are chosen to form an orthonormal basis that approximately spans  $\text{Range}(M_{[k]})$, i.e. 
\begin{align}
\label{eqn:rangefindingtucker}
     U^{(k)} {U^{(k)}}^* M_{[k]} \approx M_{[k]},\quad k=1,\ldots,n.
\end{align}
Once the factor matrices $U^{(k)}$ are determined,  the core tensor $M_c$ is computed via
\begin{align}
\label{eqn:coretensor}
    M_c = M \times_1 (U^{(1)})^* \times_2 (U^{(2)})^*  \times_3 \ldots \times_n (U^{(n)})^*.
\end{align}
In practice, the range approximation in~\eqref{eqn:rangefindingtucker} often carried out using  a truncated singular value decomposition (SVD) of $M_{[k]}$, where $U^{(k)}$ is formed from the leading left singular vectors.  This procedure is called the higher-order singular value decomposition (HOSVD) \cite{HOSVD}. Other methods for computing the Tucker decomposition include higher-order orthogonal iteration (HOOI) \cite{HOOI} and gradient-based optimization \cite{doi:10.1137/070688316}.  When the tensor $M \in \C^{d\times d \times \ldots \times d}$ is symmetric, all the unfolding matrices coincide.  In that case, we only need to carry out HOSVD on one unfolding matrix and obtain $U \in \C^{d\times r}$. Then we set $U^{(1)}=U^{(2)}=\ldots =U^{(n)}=U$ in \eqref{eqn:coretensor} with the resulting core tensor $M_c\in \C^{r\times r \times \ldots \times r}$ also being symmetric. For further details on tensor decomposition algorithms, the reader is referred to the comprehensive overview provided in \cite{tamara}. 

Lastly, we define the Frobenius norm of a tensor $M\in \C^{d_1\times d_2 \times \ldots \times d_n}$ to be 
\begin{align*}
    \Vert M \Vert_F = \sqrt{\sum_{i_1,\ldots,i_n=1}^{d_1,\ldots,d_n}|M_{i_1,\ldots,i_n}|^2}. 
\end{align*}

\subsection{Randomized range-finding algorithms}
\label{sec:randomrangefinding}
For large tensors, computing the truncated SVD of unfolding matrices can be prohibitively expensive. The randomized range-finding algorithm \cite{HMT} efficiently finds an appropriate approximation $U^{(k)}$ via randomized sketching.

Consider a scenario where we have a large matrix $M \in \C^{d\times m}$, and our goal is to find a matrix $U\in \C^{d\times r}$ with $r \ll m,d$ such that
$UU^* M \approx M.$
Following \cite{HMT}, the randomized range-finding algorithm is carried out in three steps:
\begin{enumerate}
\item Draw a random sketch matrix $S \in \C^{m \times r}$, for instance a Gaussian random matrix.
\item Form the $d\times r$ matrix $Y=MS$.
\item Compute $U$ via the QR decomposition $Y=UR$.
\end{enumerate}
The computational complexity of the randomized range-finding algorithm is only $\cO(mr+dmr+dr^2)$ compared to finding the range via a truncated SVD, which has complexity $\cO(md^2)$.  

\subsection{CUR decomposition}
\label{sec:CUR}
Another type of low-rank decomposition used in this work is the CUR decomposition \cite{doi:10.1073/pnas.0803205106}. Given a matrix $M \in \C^{d\times m}$ and a target rank (greater than or equal to) $r$,  one constructs the CUR approximation by  selecting index sets $\cS \subset [d]$ and $\cJ \in [m]$, with $|\cS|,|\cJ|\ge r$, sampled without replacement. The CUR approximation is then given by
$$
M \approx CU^\dag R
$$
where $C=M(:,\cJ), U=M(\cS,\cJ)$ and $R=M(\cS,:)$. Here,  $U^\dag$ denotes the Moore–Penrose pseudoinverse  of $U.$  In the special case where $M$ is Hermitian,  a Hermitian CUR approximation can be obtained by setting $\cS=\cJ$ and $R=C^*$.

This approach extends naturally to tensors $M \in \C^{d_1 \times d_2\times \cdots \times d_n}$, known as  \emph{mode-wise CUR decomposition}~\cite{10.5555/3546258.3546443}. For each mode $k = 1, \dots, n$, we sample index set $\cS_k \subset [d_k]$ and $\cJ_k \subset [\prod_{j \neq k} d_j]$, with $|\cS_k|, |\cJ_k| \ge r$ and  extract the following sub-matrices and sub-tensor 
$$
C^{(k)}= M_{[k]}(:,\cJ_k), \quad U^{(k)}=C^{(k)}(\cS_k,:) \quad \text{ and } \quad R=M(\cS_1,\cS_2,\cdots,\cS_n).
$$
 The approximation in Tucker format  is then given by 
\begin{align*}
    M\approx R \times_1 (C^{(1)} (U^{(1)})^\dag)^T \times_2 (C^{(2)} (U^{(2)})^\dag)^T \times \cdots \times (C^{(n)} (U^{(n)})^\dag)^T.   
\end{align*}
If the tensor $M$ is symmetric,  a symmetric approximation can be more easily constructed by selecting  $\cS_1=\cS_2=\cdots =\cS_n$ and $\cJ_1=\cJ_2=\cdots =\cJ_n$, so that $C^{(1)}=C^{(2)}=\cdots=C^{(n)}, U^{(1)}=U^{(2)}=\ldots=U^{(n)}.$ 

A widely used strategy for selecting rows or columns  is \textit{norm-weighted sampling}, in which indices are sampled with probabilities proportional to the squared Frobenius norm of the corresponding rows or columns~\cite{Deshpande2006AdaptiveSA,10.1145/1039488.1039494}.  This approach extends naturally to tensors as discussed in \cite{doi:10.1137/060665336}. For instance, when sampling rows from a matrix $A\in \C^{m\times n}$, each row index $i\in [m]$ is chosen with probability
\begin{align}
\label{eqn:cur_sampling}
    p(i) = \frac{\|A(i,:)\|_F^2}{\|A\|_F^2}
\end{align}
which yields stable low-rank approximations with provable error guarantees.

\subsection{Introduction to the cryo-EM reconstruction problem}
\label{sec:intro_cryo_em}
We now give a brief introduction to the cryo-EM reconstruction problem. For a more comprehensive introduction to this topic, we refer the reader to \cite{Singer2018MATHEMATICSFC,9016106}. 

Let $V_\star:\mathbb{R}^3 \to \mathbb{R}$ denote the electrostatic  potential created by a molecule, which we will refer to as the \textit{volume} throughout this work. We consider a simplified cryo-EM imaging model, assuming perfect centering of the volume in the projected images. Under this assumption, the observed images, denoted by $I_j, j=1,\ldots,N$, are generated according to
\begin{align}
\label{eqn:cryoemimages}
    I_j(x,y) = \mathcal{H}_j * \cI[V_\star,R_j](x,y) + \epsilon_j(x,y), \quad (x,y)\in \mathcal{X}, \quad j=1,\ldots,N,
\end{align}
where $\mathcal{X}$  denotes an $m\times m$ equispaced Cartesian 2-D grid, and
\begin{align*}
    \mathcal{I}[V_\star,R_j] (x,y)= \int_{\R} (R_j^{-1}\circ V_\star)(x,y,z) \, {\rm d}z
\end{align*}
is the ideal noiseless 2-D projection obtained from the volume $V_\star$, and the 3-D rotation $R_j \in \mathcal{SO}(3)$  modeling the unknown orientation of the molecule in the $j$-th image.  Here,  $$(R_j^{-1}\circ V_\star)(x,y,z) := V_\star(R_j(x,y,z)^T)$$ denotes the rotation of $V_\star$ by  $R_j^{-1}$ evaluated at the given spatial coordinates.   We use $R_j^{-1}$ rather than $R_j$ to rotate the volume so that the model is consistent with the formalism used in \cite{Sharon_2020}, which we build upon in Section~\ref{sec:choicebasis}.  The operator $\mathcal{H}_j$ denotes the inverse Fourier transform of the contrast transfer function (CTF) applied to the projection image $\cI[V_\star, R_j]$ and the convolution $*$ is defined as  
\begin{align*}
    (f*g)(x)=\int_{\R^2} f(x-y)g(y)\mathrm{d}y.
\end{align*}
While CTF estimation is an important practical consideration, it is not the focus of this paper. Hence, we assume that the inverse Fourier transforms of the CTFs $\mathcal{H}_j$ are known.    We assume that $\epsilon_j(\cdot)$  is an isotropic, uncorrelated Gaussian random field (white noise)  with  variance $\sigma^2$. Likewise, we assume that $\sigma^2$ is known, although in practice it may be estimated from the pixel intensities in the corners of the images.  The goal of cryo-EM reconstruction is to recover the unknown volume  $V_\star$  from the observed images $\{I_j\}_{j=1}^N$. In typical experimental settings,  the sample size of images  $N$ is large, while the SNR can be quite low.

We denote the Fourier transform of a function $f:\R^d\to \R$ by
\begin{align*}
    \hat f(\xi) = \int_{\R^d} f(x) e^{-2\pi i\xi\cdot x}\,{\rm d}x 
\end{align*}
where $\xi\in \R^d$ denotes the  variable in Fourier domain. For any three-dimensional volume $V:\R^3 \to \R$, the 2-D projection along the $z$-axis is given by 
\begin{align*}
    \mathcal{I}(x,y)=\int_\R V(x,y,z)\, {\rm d}z 
\end{align*}
and the Fourier transform of $I(x,y)$ can be obtained from the Fourier transform of $V$ by slicing  through the origin of $V$ parallel to the $z$-axis, i.e., 
\begin{align}
\label{eqn:fourierslice}
     \hat{\mathcal{I}}(\xi_x,\xi_y) =  \hat V (\xi_x,\xi_y,0).
\end{align}
This is known as the \textit{Fourier slice theorem} \cite{Bracewell1990NumericalT}, which implies that 
\begin{align}
\label{eqn:fourierprojection}
     \hat \cI[ V_\star, R_j] (\xi_x,\xi_y) 
    = (R_j^{-1} \circ \hat V_\star)(\xi_x,\xi_y,0),  
\end{align}
and the image formation model~\eqref{eqn:cryoemimages} can  be written in Fourier domain as 
\begin{align}
    \hat I_j(\xi_x,\xi_y) =  \hat{\mathcal{H}}_j(\xi_x,\xi_y)  (R_j^{-1} \circ  \hat V_\star)(\xi_x,\xi_y,0) + \hat \epsilon_j(\xi_x,\xi_y) 
\end{align}
where we have used  the convolution theorem  $\widehat{f * g}=\hat f  \hat g.$
Modeling the image formation model in frequency space is therefore easier than in real space, so it is conventional to reconstruct $\hat V$ from $\hat I_j$ and then transform the resulting volume to real space. In the remainder of this paper, we will also refer to $\hat V$ and $\hat I_1,\ldots,\hat I_j$ as the volume and 2-D image data, respectively.

\subsection{Method of moments in cryo-EM}
\label{sec:MoMforCryoEM}
In this section, we introduce the main idea of the method of moments, which serves as the foundation for several recent reconstruction algorithms in cryo-EM~\cite{Sharon_2020,pnas.2216507120,doi:10.1137/22M1498218,khoo2023deep,doi:10.1137/22M1503828}. 
We assume that the   rotation $R \in \mathcal{SO}(3)$, modeling the unknown viewing orientation,  is drawn from a probability  density $\mu$ over $\mathcal{SO}(3)$, referred to as the \textit{rotational density}.   The first three moments of the Fourier 2-D projections~\eqref{eqn:fourierprojection}  are  given by 
\begin{align}
\label{eqn:M1}
    &  \cM^{(1)}[ V,\mu] (\vec{\xi})=  \Int \hat \cI[V, R](\vec{\xi} ) \mu(R)\,{\rm d}R, \\
\label{eqn:M2}
    &  \cM^{(2)}[ V,\mu] (\vec \xi_1, \vec \xi_2)=\Int \hat \cI[V, R](\vec \xi_1)\hat \cI[ V, R]^*(\vec \xi_2) \mu(R)\,{\rm d}R, \\
\label{eqn:M3}
    &  \cM^{(3)}[ V,\mu](\vec\xi_1,\vec\xi_2,\vec\xi_3) = \Int \hat \cI[V,R](\vec \xi_1)\hat \cI[ V,R](\vec \xi_2) \hat \cI[ V,R](\vec \xi_3) \mu(R)\,{\rm d}R,
\end{align}
where $\vec \xi, \vec \xi_1, \vec \xi_2, \vec \xi_3\in \R^2$  and ${\rm d}R$ denotes the \textit{Haar measure}  over $\mathcal{SO}(3)$~\cite{chirikjian2016harmonic}.
We note that the complex conjugate in~\eqref{eqn:M2}  ensures that $\mathcal{M}^{(2)}$ is Hermitian when viewed as a matrix. We further assume that the volume $V$ is band-limited by $1/2$, i.e.,  its Fourier transform  is compactly supported in $[-1/2,1/2]^3$. Consequently, the $\vec\xi$'s in~\eqref{eqn:M1}-\eqref{eqn:M3} are restricted to the square $[-1/2,1/2]^2$. Using the Euler angle description of rotation matrices,   the Haar measure on $\mathcal{SO}(3)$ can be expressed as
\begin{align}
\label{eqn:haar}
    \mathrm{d}R = \frac{1}{8\pi^2} \sin\beta \,{\rm d}\alpha \, {\rm d}\beta \, {\rm d}\gamma, \qquad \alpha \in [0,2\pi], \quad
    \beta \in [0,\pi], \quad \gamma \in [0,2\pi],
\end{align}
where $(\alpha,\beta,\gamma)$ are the Euler angles under the ``ZYZ convention''. A brief introduction to Euler angles is provided in  Appendix \ref{sec:constructSO3quadratures}.

Typically, we   observe the 2-D Fourier images only on a discrete set of  $m \times m$  grid points, denoted by  $\Xi\in \R^2$. Consequently, we  treat the 2-D projection in~\eqref{eqn:fourierslice}, when evaluated on $\Xi$, as a vector in $\C^d$, where $d=m^2$ represents the total number of pixels (image dimension). Similarly,  the discretized first three moments have sizes $m^2, m^2\times m^2$, and $m^2\times m^2\times m^2$, respectively.  Throughout the paper, we use the vector $\hat\cI[V,R](\Xi)$ to denote the evaluation of \eqref{eqn:fourierslice} on $\Xi$ and the tensors $\cM^{(k)}[V, R](\Xi^{\otimes k})$, for $k=1,2,3$ to denote the evaluation of   \eqref{eqn:M1}-\eqref{eqn:M3} on $\Xi^{\otimes k}$.

Let $\hat V_\star$ and $\mu_\star$ denote the unknown ground truth volume and rotational density. For $k=1,2,3$, let $\Mest^{(k)} \in (\C^{d})^{\otimes k}$ denote the sample estimates of the moments $\cM^{(k)}[V_\star, \mu_\star](\Xi^{\otimes k})$.  We provide the details of obtaining the estimates in Section~\ref{sec:reduced_MoMs}, Section~\ref{sec:CTFs} and Appendix~\ref{sec:debias}. Given these estimates, the method of moments attempts to recover the pair $(V_\star, \mu_\star)$ by solving the optimization problem:
\begin{align}
\label{eqn:originalMoM}
    \min_{(V,\mu)\in \mathcal{C}} \sum_{k=1}^3 \lambda_k \Vert   \cM^{(k)}[V, \mu](\Xi^{\otimes k})-\Mest^{(k)}\Vert^2_F,
\end{align}
where $\mathcal{C}$ denotes a constraint set requiring $V$ to be in $L^2(\R^3)$ and be band-limited, and $\mu$ to be a probability density function over $\mathcal{SO}(3)$, and $\lambda_k \ge 0$ are suitably chosen weights.  In practice, finding the global minima of \eqref{eqn:originalMoM} is  challenging due to its non-convexity and the significant computational cost of forming the moments. For instance, forming  the empirical third-order moment from $N$ images  incurs a cost of $\mathcal{O}(Nd^3)$  where $d$ is the number of pixels in an image. Addressing these computational challenges  is the central focus of this paper.

\section{Subspace method of moments}
\label{sec:subspaceMoM}
In this section, we present our main contribution: an efficient algorithm for solving the method of moments. Firstly, we consider the case in which approximations of the ranges of the unfolding matrices of the estimated moments are available. For $k=1,2,3$, let $U^{(k)}\in \C^{d\times r_k}, r_k \ll d$ denote the basis matrices corresponding to the approximated ranges with orthonormal columns. Intuitively, our approach relies on  compressing the estimated moments using the  Tucker decomposition, as described in~\eqref{eqn:tucker}. More precisely, we define the projection operators  $\mathcal{P}^{(k)}_{U^{(k)}}: (\C^{d})^{\otimes k} \to (\C^{r_k})^{\otimes k}$, for $k=1,2,3$, and apply them to the estimated moments as follows:
\begin{align}
\label{eqn:compressM1}
    & \mathcal{P}^{(1)}_{U^{(1)}}(\overline M^{(1)})  = (U^{(1)})^*\overline M^{(1)} , \\ 
    \label{eqn:compressM2}
    &\mathcal{P}^{(2)}_{U^{(2)}}(\overline M^{(2)})  = (U^{(2)})^* \overline M^{(2)}(U^{(2)}), \\
    \label{eqn:compressM3}
    &\mathcal{P}^{(3)}_{U^{(3)}}(\overline M^{(3)})  =  \overline M^{(3)} \times_1 U^{(3)} \times_2 U^{(3)} \times_3 U^{(3)}.
\end{align}
This compression yields lower-dimensional representations of the estimated moments,  significantly reduces their sizes, especially for the estimated third-order moment $\overline M^{(3)}$, whose size   is reduced from $\mathcal{O}(d^3)$ to $\mathcal{O}(r_3^3)$.  Similarly, we apply the same projection operators $\mathcal{P}^{(k)}_{U^{(k)}}$ to the moments $\cM^{(k)}[V, \mu](\Xi^{\otimes k})$, for $k=1,2,3$, which are formed from the volume $\hat V$ and the rotational density $\mu$. We refer to the resulting compressed representations, $\cM_c^{(k)}[V, \mu](\Xi^{\otimes k})=\mathcal{P}_{U^{(k)}}(\cM^{(k)}[V, \mu](\Xi^{\otimes k}))$, for $k=1,2,3$ as the \textit{subspace moments}. 
 
Using these compressed representations, we may replace the original MoM problem~\eqref{eqn:originalMoM} with a reduced optimization  problem, given by  
\begin{align}
\label{eqn:reducedcost}
  \min_{(V,\mu)\in  \mathcal{C}} \quad  \sum_{k=1}^3 \lambda_k \Vert \mathcal{P}^{(k)}_{U^{(k)}} (\cM^{(k)}[  V,\mu] (\Xi^{\otimes k}))- \mathcal{P}^{(k)}_{U^{(k)}}(\overline M^{(k)}) \Vert_F^2 .
\end{align}
We refer to this formulation as the \textit{subspace method of moments (SubspaceMoM)}. A similar idea has been applied in iterative refinement algorithms,  where it is known as \textit{SubspaceEM} \cite{DVORNEK2015200}.

We now briefly discuss the connection  and the difference between the two problems in~\eqref{eqn:originalMoM} and~\eqref{eqn:reducedcost}. When $r_k = d$, the matrices $U^{(k)}$ span the entire ranges of the moments, so the operators $\mathcal{P}^{(k)}_{U^{(k)}}$ are invertible.  In this case, the optimization problem in \eqref{eqn:reducedcost} is equivalent to the original MoM problem~\eqref{eqn:originalMoM}. When $r_k < d$, the operators $\mathcal{P}^{(k)}_{U^{(k)}}$ are not invertible, and thus the reduced problem~\eqref{eqn:reducedcost} is no longer equivalent to \eqref{eqn:originalMoM}. Nonetheless, by the orthogonality of the columns of $U^{(k)}$, we have the inequality
\begin{align}
\label{eqn:subspace_MoM_intro}
\Vert \mathcal{P}^{(k)}_{U^{(k)}} (\cM^{(k)}[V,  \mu] (\Xi^{\otimes k}))- \mathcal{P}^{(k)}_{U^{(k)}}(\overline M^{(k)}) \Vert_F^2 
 \leq \Vert   \cM^{(k)}[  V,  \mu](\Xi^{\otimes k})-\Mest^{(k)}\Vert^2_F. 
 \end{align}
We also have the following upper bound:
\begin{align}
    \Vert   \cM^{(k)}[  V,  \mu](\Xi^{\otimes k})-\Mest^{(k)}\Vert_F & \le  \Vert \mathcal{P}^{(k)}_{U^{(k)}} (\cM^{(k)}[V,  \mu] (\Xi^{\otimes k}))- \mathcal{P}^{(k)}_{U^{(k)}}(\overline M^{(k)}) \Vert_F \nonumber \\ 
    & +    \Vert \mathcal{P}^{(k)\,\perp}_{U^{(k)}} ( \cM^{(k)}[  V,  \mu](\Xi^{\otimes k})) \Vert_F +  \Vert \mathcal{P}^{(k)\,\perp}_{U^{(k)}}(\Mest^{(k)}) \Vert_F
\end{align}
where $\mathcal{P}^{(k)\,\perp}_{U^{(k)}}$ denotes the projection operator onto the orthogonal complement of the range of  $U^{(k)}$. In particular, the orthogonal projection components are given by:
\begin{align*}
    \mathcal{P}^{(1)\,\perp}_{U^{(1)}}(\Mest^{(1)})  &= (I_d-U^{(1)}(U^{(1)})^*) \overline M^{(1)}, \\ 
    \mathcal{P}^{(2)\,\perp}_{U^{(2)}}(\Mest^{(2)})  &= (I_d-U^{(2)}(U^{(2)})^*)  \overline M^{(2)} (I_d-U^{(2)}(U^{(2)})^*), \\
    \mathcal{P}^{(3)\,\perp}_{U^{(3)}}(\Mest^{(3)})  &=  \overline M^{(3)} \times_1 (I_d-U^{(3)}  (U^{(3)})^*) \times_2 (I_d-U^{(3)}  (U^{(3)})^*) \times_3 (I_d-U^{(3)}  (U^{(3)})^*).
\end{align*}
If $U^{(1)}=I_d$ and the equation \eqref{eqn:tucker} approximately holds for $\overline M^{(k)}$ with the range matrices given by the provided $U^{(k)}$ and the core tensor given by $\mathcal{P}^{(k)}_{U^{(k)}}(\overline M^{(k)})$ when $k=2,3$. In this case, the quantity $\sum_{k=1}^3 \lambda_k \Vert \mathcal{P}^{(k)\,\perp}_{U^{(k)}}(\Mest^{(k)}) \Vert_F^2$ is  expected to be small. Consequently, the discrepancy between the SubspaceMoM and the original MoM problem is approximately given by
\begin{align*}
\lambda_2 \Vert \mathcal{P}^{(2)\,\perp}_{U^{(2)}} ( \cM^{(2)}[  V,  \mu](\Xi^{\otimes 2})) \Vert_F^2+\lambda_3 \Vert \mathcal{P}^{(3)\,\perp}_{U^{(3)}} ( \cM^{(3)}[  V,  \mu](\Xi^{\otimes 3})) \Vert_F^2.
\end{align*}
This gap reflects the portion of the moments that lies outside the range of the estimated subspaces. It can be made small if $(V,\mu)$ approximates $(V_\star,\mu_\star)$, provided that accurate low-rank approximations to the moments are available.

In order to numerically solve the reduced optimization problem~\eqref{eqn:reducedcost}, we need to discretize it.  We assume that both the volume $\hat V$ and the rotational density $\mu$ are expanded in some suitable basis functions. The integrals over $\mathcal{SO}(3)$ appearing in~\eqref{eqn:M1}-\eqref{eqn:M3} are then approximated using numerical quadrature. In addition to yielding numerical approximations of the moments, quadrature has an additional benefit: it provides a natural approximate decomposition of $\cM^{(k)}[V, \mu], k=1,2,3$, as a sum of rank-$1$ tensors. This property plays a key role in our proposed approach.  The details of formulating the SubspaceMoM problem using basis expansions and numerical integration are provided in Section~\ref{sec:reduced_MoMs}. In Section~\ref{sec:estimation},  we describe how to determine the projection operators  $\mathcal{P}^{(k)}_{U^{(k)}}$ using some randomized numerical linear algebra techniques.  In Section~\ref{sec:steerablePCA}, we  introduce symmetry properties of the rotational distribution and offer an approach for improving subspace moment estimation. Finally in Section~\ref{sec:numerical_optim}, we presents further details on the numerical optimization of the reduced problem in~\eqref{eqn:reducedcost}.

\subsection{Details of the reduced moment problem} 
\label{sec:reduced_MoMs}
In this Section, we present an approach for efficiently discretizing the reduced moment problem and evaluating the associated objective cost function in Section~\ref{sec:compressed_rep}, and we discuss the computational complexity in Section~\ref{sec:complexity}.

\subsubsection{Compressed representation of moments}
\label{sec:compressed_rep}
We assume that the Fourier transform of the volume and the rotational density can be expressed in terms of suitable basis functions, as follows:
\begin{align}
\label{eqn:volrep}
    \hat V_a (\xi_x,\xi_y,\xi_z) &= \sum_{i \in \mathcal{S}_V} a_i \phi_i(\xi_x,\xi_y,\xi_z), \\ 
\label{eqn:viewrep}
    \mu_b(R) &= \sum_{j \in \cB_\mu} b_j \psi_j(R).
\end{align}
where $\phi_i:\R^3\to \C $ and $ \psi_j:\mathcal{SO}(3)\to \C$ are basis functions with  index sets $\cB_V$ and $\cB_\mu$. The coefficients $a \in \C^{|\cB_V|}$ and $b \in \C^{|\cB_\mu|}$ are the expansion coefficients, subject to additional constraints that ensure $\hat V_a$ is the Fourier transform of a real-valued 3-D volume and $\mu_b$ is a probability density over $\mathcal{SO}(3)$. Specific examples of suitable basis functions and constraints are discussed in Section \ref{sec:choicebasis} and Section \ref{sec:constraints}.

Using the Fourier slice theorem \eqref{eqn:fourierslice}, the noiseless 2-D projection  generated by the volume  $\hat V_a$ and  3-D rotation $R$ is  given by 
\begin{align}
\label{eqn:projrep}
    \hat\cI[V_a,R](\xi_x,\xi_y) = \sum_{i \in \cB_V} a_i  (R^{-1}\circ \phi_i)(\xi_x,\xi_y,0) 
\end{align}
for $(\xi_x,\xi_y)\in \Xi$. To express this projection in matrix-vector form,    we define the $d\times |\cB_V|$ matrix $\Phi[R,\Xi]$ by 
\begin{align*}
(\Phi[R,\Xi])_{i,j} = (R^{-1}\circ \phi_j)(\vec\xi_i,0),
\end{align*}
for  $\vec\xi_i \in \Xi, i=1,\ldots,d$ and $j \in \cB_V.$
Using this notation, the projection~\eqref{eqn:projrep} becomes
\begin{align}
\label{eqn:imageexpand}
    \hat\cI[V_a,R](\Xi) = \Phi[R,\Xi] \cdot a . 
\end{align}
To represent the rotational density compactly, we define the mapping $\Psi: \mathcal{SO}(3) \to \C^{1\times |\cB_\mu|}$ by
\begin{align*}
    (\Psi[R])_j = \psi_j(R),\quad j\in \cB_\mu
\end{align*}
which allows us to rewrite~\eqref{eqn:viewrep} as
\begin{align}
\label{eqn:densityexpand}
    \mu_b(R) = \Psi[R] \cdot  b. 
\end{align}

Inserting the equations \eqref{eqn:imageexpand} and \eqref{eqn:densityexpand} into the definitions of the moments given  in \eqref{eqn:M1}-\eqref{eqn:M3},  we can express  the  moments  $\cM^{(k)}[a, b](\Xi^{\otimes k})$ in terms of expansion coefficients $a$ and $b$. Specifically, we obtain
 \begin{align}
 \label{eqn:newM1}
     &\cM^{(1)}[a,b](\Xi)  =  \Int (\Phi[R,\Xi] \cdot a) (\Psi[R] \cdot b) \,{\rm d}R, \\
     \label{eqn:newM2}
     &\cM^{(2)}[a, b](\Xi^{\otimes 2})  =  \Int (\Phi[R,\Xi] \cdot a)(\Phi[R,\Xi] \cdot a)^* (\Psi[R] \cdot b) \,{\rm d}R, \\ 
     \label{eqn:newM3}
     &\cM^{(3)}[a, b](\Xi^{\otimes 3})  =  \Int (\Phi[R,\Xi] \cdot a)^{\otimes 3} (\Psi[R] \cdot b) \,{\rm d}R.
 \end{align}
The corresponding   subspace moments are obtained by applying the projection operators defined in~\eqref{eqn:compressM1}–\eqref{eqn:compressM3}:
\begin{align*}
\cM_c^{(k)}[a,b](\Xi^{\otimes k})=\mathcal{P}_{U^{(k)}}^{(k)}(\cM^{(k)}[a,b](\Xi^{\otimes k})), \quad k=1,2,3.
\end{align*}

To evaluate the integrals over $\mathcal{SO}(3)$, we use numerical quadrature with weights and nodes $\{(w_i,R_i),i=1,\ldots,Q\}$ for $w_i> 0$ and $R_i\in \mathcal{SO}(3)$. The moments can then be approximated as follows: 
\begin{align*}
    &\cM^{(1)}[a,b](\Xi) \approx  \sum_{i=1}^Q w_i (\Phi[R_i,\Xi] \cdot a) (\Psi[R_i] \cdot b), \\
    &\cM^{(2)}[a,b](\Xi^{\otimes 2})  \approx \sum_{i=1}^Q w_i (\Phi[R_i,\Xi] \cdot a) (\Phi[R_i,\Xi] \cdot a)^* (\Psi[R_i] \cdot b), \\ 
     &\cM^{(3)}[a,b](\Xi^{\otimes 3}) \approx \sum_{i=1}^Q w_i (\Phi[R_i,\Xi] \cdot a)^{\otimes 3} (\Psi[R_i]\cdot  b) .
\end{align*}
We can find accurate quadrature rules  when the basis functions are appropriately chosen. However, evaluation using these rules may be prohibitively expensive in practice. As a result, we often choose to use  inexact quadrature, which is computationally more efficient. Importantly, it has been shown in~\cite{khoo2023deep} that high-quality \textit{ab initio} models can be obtained from the first two moments, even when inexact quadrature is used during reconstruction, provided that the unknown 3-D volume is represented by a neural network.
  Additional evidence is presented in Section~\ref{sec:recoverEstimatedMoMs}, where we perform \textit{ab initio} modeling using synthetic image datasets. In that setting, the use of exact quadrature would be computationally infeasible on standard desktop machines. Nevertheless, our reconstruction algorithms are able to recover satisfactory \textit{ab initio} models efficiently using inexact quadrature throughout the computation.

Using the same quadrature rule, we  approximate the subspace moments as 
\begin{align}
\label{eqn:Mc1}
&\cM_c^{(1)}[a, b](\Xi) \approx  \sum_{i=1}^Q w_i \cdot ((U^{(1)})^* \Phi[R_i,\Xi] \cdot a) (\Psi[R_i] \cdot b), \\
\label{eqn:Mc2}
&\cM_c^{(2)}[a, b](\Xi^{\otimes 2})  \approx \sum_{i=1}^Q w_i \cdot ((U^{(2)})^* \Phi[R_i,\Xi]  \cdot a ) ((U^{(2)})^* \Phi[R_i,\Xi] \cdot  a )^* (\Psi[R_i] \cdot b), \\ 
\label{eqn:Mc3}
&\cM_c^{(3)}[a ,b](\Xi^{\otimes 3}) \approx \sum_{i=1}^Q w_i\cdot  ((U^{(3)})^* \Phi[R_i,\Xi] \cdot a )^{\otimes 3} (\Psi[R_i] \cdot b).
\end{align}
The estimated subspace moments  $\overline M_c^{(k)}=\mathcal{P}_{U^{(k)}}^{(k)}(\overline M^{(k)})$, for $k=1,2,3$, can be computed efficiently from the data using the method in Section \ref{sec:estimation}. Substituting these expressions into the SubspaceMoM problem in~\eqref{eqn:reducedcost}, we reformulate it as an optimization problem over the expansion coefficients $(a,b)$, given by
\begin{align}
\label{eqn:subspaceMoM}
\min_{(a,b) \in \mathcal{C}} \sum_{k=1}^3 \lambda_k \Vert \cM^{(k)}_c[a,b](\Xi^{\otimes k}) - \overline M_c^{(k)}\Vert_F^2 
\end{align}
where the subspace moments $\cM^{(k)}_c[a,b](\Xi^{\otimes k})$ are evaluated using~\eqref{eqn:Mc1}-\eqref{eqn:Mc3}. We use the same notation $\mathcal{C}$ to denote the constraints on $(a,b)$ as in \eqref{eqn:originalMoM} and \eqref{eqn:reducedcost}. In this work,  we choose the weights to be $\lambda_k=\frac{1}{\Vert \overline M_c^{(k)} \Vert_F^2}$, which is equivalent to minimizing the relative error with equal weighting across the three moment orders. Alternative choices for the weights have been proposed in~\cite{Sharon_2020,doi:10.1137/22M1503828}. While our algorithm can be adapted to accommodate those alternatives, we find that the current choice leads to significantly faster convergence in practice.

\subsubsection{Computational complexity of cost-gradient evaluation}
\label{sec:complexity}
We now discuss the computational cost of evaluating the cost function in \eqref{eqn:subspaceMoM} during a gradient-based optimization.  The  terms $(U^{(k)})^*\Phi[R_i,\Xi] \in \C^{r_k \times |\cB_V|}$, for $k=1,2,3$, and $i=1,\ldots,Q$ as well as $\Psi[R_i] \in \C^{1 \times |\cB_\mu|}$, for $i=1,\ldots,Q$  can be precomputed and reused across iterations.   The time complexity needed  for the precomputation is $\mathcal{O} \left( Q|\cB_{\mu}|(d+1)\sum_{k=1}^3r_k\right)$
and the space complexity is $\mathcal{O}(Q |\cB_V| \sum_{k=1}^3 r_k + Q |\cB_\mu|)$. In Section~\ref{sec:choicebasis}, we provide details for the specific choice of basis used in our implementation and, in particular,  provide bounds for $|\cB_\mu|$ and $|\cB_V|$ in terms of the corresponding truncation limits.
During each evaluation of the cost function, we compute the matrix-vector products $\left\{(U^{(k)})^*\Phi[R_i, \Xi] a\right\} \in \C^{r_k}$  and  $\Psi[R_i] b \in \C$, which together require $\mathcal{O}(Q|\cB_V|(r_1+r_2+r_3)+Q|\cB_\mu|)$ work.  For each $k=1,2,3$, we compute the order-$k$ outer product of the vector  $\left\{(U^{(k)})^*\Phi[R_i,\Xi] \cdot a\right\}$, scale it by the weight $w_i \cdot (\Psi[R_i] \cdot b)$,  and accumulate to form the subspace moment. This process incurs a total cost of $\mathcal{O}(Q(r_1+r_2^2+r_3^3))$. Lastly, we compute the Frobenius distances between the  subspace moments and the estimated moments, which has an additional cost of  $\mathcal{O}(r_1+r_2^2+r_3^3).$  Therefore, the overall time complexity for a single cost evaluation is
$\mathcal{O}\left(Q \left( |\cB_V| |\cB_\mu|\sum_{k=1}^3 r_k+ \sum_{k=1}^3 r_k^k + |\cB_\mu|\right)\right).$
The evaluation of the gradient of the cost function in \eqref{eqn:subspaceMoM}  is of the same order, and we omit the details here.
\begin{remark}
For clarity, we have assumed all the subspace moments are evaluated using the same quadrature rule. In practice, the number of quadrature nodes used to integrate the third subspace moments may be larger than those used in the quadrature rules for the first two subspace moments. In Appendix~\ref{sec:constructSO3quadratures}, we discuss constructing different quadrature rules for the first three moments when the basis functions are chosen according to Section~\ref{sec:choicebasis}.
\end{remark}

\subsection{Finding subspaces via randomized numerical linear algebra}
\label{sec:estimation}
In this section, we discuss a reduced complexity method for forming the subspace moments $\overline M_c^{(k)}$, for $k=1,2,3$. We assume the subspaces determined by the matrices $U^{(k)}\in \C^{d\times r_k}$ are given. The procedure for computing these subspaces will be discussed in Section~\ref{sec:findingsubspaces} and Section~\ref{sec:CTFs}.

For clarity, we assume that the data are given by $\hat I_j=\hat \cI[ V_\star,R_j] (\Xi),  j=1,\ldots,N$, which are not contaminated by white noise or contrast transfer functions (CTFs).  This assumption approximately holds, for example, when the images have been successfully denoised~\cite{BHAMRE201672,9506435,SHI2022107018,marshall_mickelin_shi_singer_2023}. We will discuss how to handle  CTFs  without performing image denoising in Section~\ref{sec:CTFs}.  When white noise is present, the bias terms in the estimated subspace moments can be corrected, according to the method in Appendix~\ref{sec:debias}.  Under these assumptions, the estimated third moment is given by
\begin{align}
\label{eqn:estimatedM3}
\overline M^{(3)} = \frac{1}{N}\sum_{j=1}^N \hat I_j^{\otimes 3}.
\end{align}
Applying the projection operator in \eqref{eqn:compressM3}, we obtain the estimated  third subspace moment 
\begin{align}
\label{eqn:estimatedMc3}
\overline M_c^{(3)} = \frac{1}{N}\sum_{j=1}^N \mathcal{P}_{U^{(3)}}^{(3)} (\hat I_j^{\otimes 3})= \frac{1}{N} \sum_{j=1}^N ((U^{(3)})^* \hat I_j)^{\otimes 3}.
\end{align}
From~\eqref{eqn:estimatedMc3}, we see  that  for each $I_j$ where $j=1,\ldots,N$,  one can first compute $(U^{(3)})^* \hat I_j$ and  form $((U^{(3)})^* \hat I_j)^{\otimes 3}$, and finally accumulate the result. This leads to a total computational cost of $\mathcal{O}(Ndr_3+Nr_3^3)$.  This computational strategy,  commonly referred to as  ``streaming''~\cite{doi:10.1137/19M1257718,doi:10.1137/18M1201068}, avoids the explicit formation of $\overline M^{(3)} \in \C^{d\times d \times d}$, which would necessitate  $\mathcal{O}(Nd^3)$ operations. 

The formation of the estimated first two moments is computationally and conceptually simpler, and is given by 
\begin{align}
\label{eqn:estimatedMc1}
\overline M_c^{(1)}& =\frac{1}{N}\sum_{j=1}^N (U^{(1)})^* \hat I_j, \\ 
\label{eqn:estimatedMc2}
\overline M_c^{(2)}& =\frac{1}{N}\sum_{j=1}^N ((U^{(2)})^* \hat I_j) ((U^{(2)})^* \hat I_j)^*. 
\end{align}
These can be computed simultaneously with  $\overline M_c^{(3)}$, resulting in a total  complexity of  $\mathcal{O}(Nd(r_1+r_2)+N(r_1+r_2^2))$.

\subsubsection{Finding subspaces via randomized range-finding}
\label{sec:findingsubspaces}
To determine the matrices $U^{(k)}$ $\in$ $\C^{d\times r_k}$, for $k=1,2,3$, used in the previous sections, we use a specific type of randomized range-finding algorithm, as describe in Section~\ref{sec:randomrangefinding}. This algorithm is applied to the   matrices $\overline M^{(2)} \in \C^{d\times d}$ and $\overline M^{(3)}_{[1]} \in \C^{d\times d^2}$.  For clarity, we assume that the images are noiseless.

To sketch the estimated second-order moment,  
\begin{align}
\label{eqn:empM2}
\overline M^{(2)}=\frac{1}{N} \sum_{j=1}^N \hat I_j \hat I_j^*,  
\end{align}
we apply a Gaussian random matrix $G \in \C^{d\times s}$  $s < d$ being a sufficiently large sampling size. The resulting smaller matrix of size $d \times s$ is then computed as
\begin{align}
\label{eqn:sketch_M2}
\overline M^{(2)} G = \frac{1}{N} \sum_{j=1}^N \hat I_j (G^* \hat I_j)^* . 
\end{align}
This computation can be carried out efficiently in a streaming fashion with a total cost of $\mathcal{O}(Nd s)$ operations. When the sampling size $s$ is sufficiently large, the range of $\overline{M}^{(2)}G$ closely approximates the range of $\overline{M}^{(2)}$. To select the bases in a data-dependent way, we compute the SVD of $\overline{M}^{(2)} G$ and denote the singular values by $\sigma^{(2)}_1\ge \sigma^{(2)}_2 \ge \ldots \ge \sigma^{(2)}_s\ge 0$. We determine $r_2>0$, i.e. the dimension of the subspace used to compress $\overline{M}^{(2)}$,  by selecting the smallest $r$ such that
\begin{align*}
\frac{\sum_{i=1}^r (\sigma^{(2)}_i)^2}{\sum_{i=1}^s (\sigma^{(2)}_i)^2} > 1-\tau^{(2)}
\end{align*}
where $\tau^{(2)}>0$ is a small value. The  range matrix $U^{(2)} \in \mathbb{C}^{d\times r_2}$    is  then formed from the first $r_2$ left singular vectors. This singular value thresholding procedure also achieves a denoising effect on the second moment~\cite{10.3150/24-BEJ1804}.

For the estimated third-order  moment,  sketching its unfolding matrix $\overline M^{(3)}_{[1]} \in \C^{d\times d^2}$ with a random Gaussian sketch matrix of size $d^2\times s$ through streaming incurs a complexity of  $\mathcal{O}(Nd^2 s)$, which can be prohibitively expensive for large image size $d$. 
To address this, we use a more efficient approach based on the \emph{face-splitting product} \cite{Ahle2019AlmostOT,slyusar1999face}. Let $G^{(1)},  G^{(2)} \in \C^{s\times d}$ be two independent Gaussian random matrices.  Their face-splitting product, denoted by $G^{(1)} \bullet  G^{(2)}: \C^d \times \C^d \to \C^s$, is defined such that for any $\vec x,\vec y \in \C^d$, 
\begin{align*}
    (G^{(1)} \bullet  G^{(2)})(\vec x \otimes \vec y)= (G^{(1)} \vec x)\odot (G^{(2)}\vec y)
\end{align*}
where $\odot$ denotes the element-wise (Hadamard) product. Using this operator, we can form a sketch of the third-order moment
\begin{align}
\label{eqn:tensor_sketch}
    Y^{(3)}  = \frac{1}{N} \sum_{j=1}^N  \hat I_j ((G^{(1)}\bullet G^{(2)})( \hat I_j \otimes  \hat I_j))^T = \frac{1}{N} \sum_{j=1}^N  \hat I_j ((G^{(1)} \hat I_j) \odot  (G^{(2)} \hat I_j))^T. 
\end{align}
Subsequently, we compute the SVD of $Y^{(3)}$ to obtain the singular values $\sigma^{(3)}_1\ge \sigma^{(3)}_2 \ge \ldots \ge \sigma^{(3)}_s \ge 0$ and the corresponding left singular vectors.  Similarly, we determine $r_3$, the dimension of the subspace for compressing $\overline{M}^{(3)}$, to be the first $r$ such that
\begin{align*}
    \frac{\sum_{i=1}^r (\sigma_i^{(3)})^2}{\sum_{i=1}^s (\sigma_i^{(3)})^2}>1-\tau^{(3)}
\end{align*}
for some small value $\tau^{(3)}>0$. The range matrix $U^{(3)}$ of the subspace  is formed by the corresponding first $r_3$ singular vectors. The parameters $s$, $\tau^{(2)}$, and $\tau^{(3)}$ need to be specified before running the estimating procedure. Larger $s$ and smaller $\tau^{(2)},\tau^{(3)}$ give a better approximation to the moment estimators but also incur greater computational cost. The effect of varying $\tau^{(2)}$ and  $\tau^{(3)}$ will be shown in Section~\ref{sec:thirdmomenteffect}. 

Alternative methods for compressing third-order moment tensors with provable guarantees have also been proposed. For instance, sketching based on Kronecker products of Johnson–Lindenstrauss transforms~\cite{10.1093/imaiai/iaaa028} provides faster computation, and recent work on symmetric Tucker decomposition via numerical optimization~\cite{doi:10.1137/23M1582928} avoids explicit formation of the full tensor. Nonetheless, the face-splitting approach described in~\eqref{eqn:tensor_sketch} is both simple to implement and empirically stable in our setting.

The first moment, $\overline{M}^{(1)}$, does not contribute useful information for compression. Therefore, in practice, we compress the first moment using $U^{(2)}$, which implies that $r_1 = r_2$. During reconstruction, we can apply the same quadrature formula to integrate both the first and second subspace moments, eliminating the need for additional precomputation of the first moment, which helps saving memory.

It is generally difficult to compute the approximation error of low-rank decompositions for the estimated second-order and third-order moments, as the original estimated moments are not fully available. To assess the quality of the approximation, we sample an index set $\mathcal{J} \subset [d]$ and evaluate the relative error over the sub-matrix or sub-tensor restricted to $\mathcal{\mathcal{J}}$. For the second-order moment, we define the relative error with respect to $\mathcal{J}$ as
\begin{align}
\label{eqn:relerr_m2}
    E^{(2)}_J = \frac{\|\bar M^{(2)}(\mathcal{J},\mathcal{J})-U^{(2)}(\mathcal{J},:) \bar M_c^{(2)} U^{(2)}(\mathcal{J},:)^*\|_F}{\|\bar M^{(2)}(\mathcal{J},\mathcal{J})\|_F}
\end{align}
and for the third-order moment, we define
\begin{align}
\label{eqn:relerr_m3}
    E^{(3)}_\mathcal{J} = \frac{\|\bar M^{(3)}(\mathcal{J},\mathcal{J},\mathcal{J})-M^{(3)}_c \times_1 U^{(3)}(\mathcal{J},:)\times_2 U^{(3)}(\mathcal{J},:) \times_3 U^{(3)}(\mathcal{J},:) \|_F}{\|\bar M^{(3)}(\mathcal{J},\mathcal{J},\mathcal{J})\|_F}.
\end{align}

\subsubsection{Finding subspaces via CUR decomposition}
\label{sec:CTFs}
We assume that the observed images are given by $\hat I_j^H= \hat H_j \odot  \hat I_j, j=1,\cdots,N,$ where $\hat H_j = \hat{\mathcal{H}}_j(\Xi) \in \C^{d\times 1}$ denotes the discretized CTF associated with the $j$-th image.  The  first three moments can be estimated  using the following least squares estimators \cite{marshall_mickelin_shi_singer_2023}
\begin{align}
    \overline{M}^{(1)} & = \left[ \sum_{j=1}^N \hat H_j\odot \hat I_j^H \right]  \oslash \left[\sum_{j=1}^N \hat H_j^{\odot^2} \right],  \\
    \overline{M}^{(2)} & =  \left [  \sum_{j=1}^N (\hat H_j\odot \hat I_j^H) (\hat H_j\odot \hat I_j^H)^* \right ] \oslash \left [ \sum_{j=1}^N \left(\hat H_j^{\odot^2}\right) \left(\hat H_j^{\odot^2}\right)^* \right], \\
    \overline{M}^{(3)} & = \left [  \sum_{j=1}^N \left(\hat H_j\odot \hat I_j^H \right)^{\otimes 3} \right] \oslash \left [\sum_{j=1}^N \left(\hat  H_j^{\odot^2} \right)^{\otimes 3}\right]
\end{align}
where $\oslash$ denotes the element-wise division. If the images contain additional white noise, we use the method in Appendix~\ref{sec:debias} to debias the moments.

To implement a streaming-based sketching method for compressing $\overline{M}^{(2)}$ and $\overline{M}^{(3)}$, the sketching operator needs to commute with element-wise division.   For this reason, we use the CUR decomposition introduced in Section~\ref{sec:CUR}, which relies only on row and column selection. We sample an index set $\mathcal{J}\subset [d]$ and construct a CUR approximation to the estimated second moment 
\begin{align*}
    \overline{M}^{(2)} \approx C^{(2)} (W^{(2)})^\dag (C^{(2)})^*
\end{align*}
where the $d \times |\mathcal{J}|$ matrix $C^{(2)}$ is formed as
\begin{align}
\label{eqn:form_C2_cur}
C^{(2)}  &= \overline{M}^{(2)}(:,\mathcal{J})  =    \left [  \sum_{j=1}^N (\hat H_j\odot \hat I_j^H) (\hat H_j\odot \hat I_j^H(\mathcal{J}))^* \right ] \oslash \left[ \sum_{j=1}^N \left(\hat H_j^{\odot^2}\right) \left(H_j^{\odot^2}(\mathcal{J})\right)^* \right], \\  
W^{(2)}   & = C^{(2)}(\mathcal{J},:). 
\end{align}
To further trim the rank of the decomposition, we compute the QR decomposition   $C^{(2)}=Q^{(2)} R^{(2)}$ where $Q^{(2)} \in \C^{d\times |\mathcal{J}|}$ and $R^{(2)}\in \C^{|\mathcal{J}|\times |\mathcal{J}|}$, followed by a truncated SVD of the $|\mathcal{J}| \times |\mathcal{J}|$ matrix 
$$
R^{(2)} (W^{(2)})^\dag (R^{(2)})^* \approx\tilde U^{(2)} \Sigma^{(2)} (\tilde U^{(2)})^*
$$ 
where $\tilde U^{(2)} \in \C^{|\mathcal{J}| \times r_2}, \Sigma^{(2)}  \in \C^{r_2\times r_2}$ with $r_2$ determined by singular value thresholding as described in Section~\ref{sec:findingsubspaces}. The resulting rank-$r_2$ approximation is $\overline M^{(2)} \approx  U^{(2)} \bar M^{(2)}_c (U^{(2)})^*$ with  $U^{(2)}=Q^{(2)} \tilde U^{(2)}, \bar M^{(2)}_c=\Sigma^{(2)}$. The computational cost is dominated by forming the submatrix in \eqref{eqn:form_C2_cur}, which is $\cO(Nd|\mathcal{J}|)$.

To construct a CUR-based Tucker approximation of the third-order moment tensor, we sample an index set $\tilde \cS\subset [d]$ and an index set $\tilde \cJ \subset [d]^2$. The set $\tilde \cJ$ is formed by independently sampling index sets $\cJ_1,\cJ_2\subset [d]$, then taking $\tilde \cJ=\cJ_1 \times \cJ_2$. 
The approximation has the form
\begin{align*}
    \overline{M}^{(3)} \approx Y  \times_1 (C^{(3)} (W^{(3)})^\dag )^T \times_2 (C^{(3)} (W^{(3)})^\dag)^T  \times_3 (C^{(3)} (W^{(3)})^\dag )^T, 
\end{align*}
where the  $d\times |\tilde \cJ|$ matrix $W^{(3)}$ is formed by
\begin{align}
\label{eqn:form_W3_cur}
 W^{(3)}  &= \overline{M}^{(3)}_{[1]}(:,\tilde \cJ)
 = \left[\left(\sum_{j=1}^N \left(\hat H_j\odot \hat I_j^H \right)^{\otimes 3} \right) \oslash 
\left(\sum_{j=1}^N \left(\hat H_j^{\odot 2} \right)^{\otimes 3}\right) \right]_{[1]}(:, \tilde \cJ) \nonumber \\
&=   \left[ \sum_{j=1}^N (\hat H_j \odot \hat I^H_j) \otimes  (\hat H_j(\cJ_1) \odot \hat I^H_j(\cJ_1))  \otimes  (\hat H_j(\cJ_2) \odot \hat I^H_j(\cJ_2)) \right]_{[1]} \nonumber \\
& \oslash  \left[  \sum_{j=1}^N \left(\hat H_j^{\odot 2} \otimes \hat H_j^{\odot 2}(\cJ_1) \otimes \hat H_j^{\odot 2}(\cJ_2) \right) \right]_{[1]},
\end{align}
 with the $|\tilde \cS|\times |\tilde \cJ|$ matrix $C^{(3)}=W^{(3)}(\tilde \cS,:)$  and the $|\tilde \cS|\times |\tilde \cS| \times |\tilde \cS|$ tensor 
\begin{align}
\label{eqn:form_Y_cur}
    Y  =\left[\sum_{j=1}^N \left(\hat H_j(\tilde \cS)\odot \hat I_j^H (\tilde \cS)\right)^{\otimes 3} \right] \oslash 
\left[\sum_{j=1}^N \left(\hat H_j^{\odot 2} (\tilde \cS)\right)^{\otimes 3}\right]. 
\end{align}
To trim this decomposition, we compute the QR decomposition  $C^{(3)}(W^{(3)})^\dag=Q^{(3)} R^{(3)}$ where $Q^{(3)} \in \C^{d\times |\tilde \cS|}$ and $R^{(3)} \in \C^{|\tilde \cS|\times |\tilde \cS|}$,  followed by forming the $|\tilde \cS|\times |\tilde \cS|\times |\tilde \cS|$ tensor 
$$X = Y \times_1 (R^{(3)})^T \times_2 (R^{(3)})^T  \times_3 (R^{(3)})^T. $$ 
We then compute the truncated SVD of the unfolding matrix
$
X_{[1]} =\tilde U^{(3)} \Sigma^{(3)} (V^{(3)})^*
$ with rank $r_3$ determined by singular value thresholding. The final approximation is given by $$\overline M^{(3)} \approx \bar M_c^{(3)} \times_1 (U^{(3)})^* \times_2 (U^{(3)})^* \times_3 (U^{(3)})^*$$ with $U^{(3)}=\overline{Q^{(3)} \tilde U^{(3)}}$ and $\bar M_c^{(3)} = X \times_1  \tilde U^{(3)} \times_2  \tilde U^{(3)} \times_3  \tilde U^{(3)}.$  The computational cost is dominated by forming the matrices and tensors in \eqref{eqn:form_W3_cur} and \eqref{eqn:form_Y_cur}, which is $\cO(Nd|\cJ_1||\cJ_2| +Nd\tilde |\tilde \cS|)$.

In this setting, the norm-weighted sampling  as described in~\eqref{eqn:cur_sampling} is not applicable because the full moment matrix or tensor is unavailable.  As an alternative, we propose a heuristic strategy motivated by the known decay of the Fourier volume. Specifically, Guinier law and Wilson statistics~\cite{Singer:ib5103} predict that the power spectrum of proteins decays rapidly at low frequencies and remains relatively flat at mid-to-high frequencies, up to the decay of the structure factors. Since we do not aim to recover mid-to-high frequency components, we may assume that these components have been removed via Fourier cropping. Then, by Fourier slice theorem,  the magnitude of the observed images $|\hat I_j(\xi_x, \xi_y)|$ should decay as $\xi_x^2 + \xi_y^2$ increases. This suggests that rows or columns associated with high-frequency indices are generally less informative and should be sampled with lower probability. 
We assign each index with a sampling probability given by the function 
\begin{align}
\label{eqn:sampling_rule}
    p(\xi_x,\xi_y)\propto (\xi_x^2+\xi_y^2)^{-1}
\end{align}
and always sample the index corresponding to $(\xi_x,\xi_y)=(0,0)$.  
While it is possible to construct better sampling distributions based on the estimated power spectrum via Wilson statistics, we find this heuristic to be simple and effective in practice.

When computing the pseudoinverse $W^\dag$ for the CUR decomposition, we cut off small  singular values in $W$ that are less than a prescribed tolerance $\epsilon>0$, which yields a  more stable low-rank approximation \cite{Nakatsukasa2020FastAS}, and we use $\epsilon=10^{-5}$ in this work.

\subsection{Moment formation under the presence of symmetry}
\label{sec:steerablePCA}

When forming empirical moments from images, the dataset can be augmented with in-plane rotated copies of the observed images. This augmentation does not alter the underlying distribution of the 3-D orientations, since there is typically no \emph{a priori} reason to assume a preferred in-plane direction in the 2-D projection space. Mathematically, we assume that the rotational distribution  is \emph{in-plane uniform}, meaning that $$\mu(R)=\mu(RR_z(\varphi)), \quad \forall \varphi \in [0,2\pi]
$$ 
where $R_z(\varphi)$ denotes the rotation about the $z$-axis by angle $\varphi.$

Exploiting this symmetry, one can perform principal component analysis (PCA)  using the estimated second moment more efficiently via a technique known as \emph{fast steerable PCA} \cite{Zhao2012FourierBesselRI}. Fast methods have been developed for computing second-order moments  using fast Fourier-Bessel expansions \cite{marshall2022fast} and exploiting the block-diagonal structure of the second moment matrix to accelerate PCA \cite{marshall_mickelin_shi_singer_2023}.

Randomized sketching provides an alternative to explicit basis expansions for achieving fast steerable PCA. Let the clean image  $\hat \cI[V,R]$ be expressed in polar coordinates  $(\kappa,\theta)$, where $0\le \kappa \le c=1/2$, $\theta \in [0,2\pi]$. Then the second moment can be written as
\begin{align*}
    \cM^{(2)}[V,\mu]((\kappa,\theta),(\kappa',\theta'))=\int_{\mathcal{SO}(3)} \hat  \cI[V,R](\kappa,\theta) \hat  \cI^*[V,R](\kappa',\theta') \mu(R)\, {\rm d}R. 
\end{align*}
If $\mu$ is inplane-uniform, the second moment satisfies
\begin{align*}
\cM^{(2)}[V,\mu]((\kappa,\theta),(\kappa',\theta'))=\cM^{(2)}[V,\mu]((\kappa,\theta+\varphi),(\kappa',\theta'+\varphi)), \quad \forall \varphi \in [0,2\pi].
\end{align*} 
Averaging over all in-plane rotations yields
\begin{align}
\label{eqn:augmentation}
\cM^{(2)}[V,\mu]((\kappa,\theta),(\kappa',\theta'))=\frac{1}{2\pi}\int_0^{2\pi}\cM^{(2)}[V,\mu]((\kappa,\theta+\varphi),(\kappa',\theta'+\varphi))\,{\rm d}\varphi, 
\end{align}
which acts as a symmetrizing operator equivalent to augmenting the dataset with all in-plane rotated versions of each image.

Sketching can be viewed as projecting the second moment against some  test functions.  Let $G^{(i)}(\kappa',\theta')$, for $i=1,\ldots,r$ be the test (sketch) functions. The $i$-th entry of the sketched second moment is computed as
\begin{align}
\label{eqn:continuous_sketch}
(G\circ \cM^{(2)})_i[V,\mu](\kappa,\theta) =\int_{0}^c \int_0^{2\pi}\cM^{(2)}[V,\mu]((\kappa,\theta),(\kappa',\theta'))  G^{(i)}(\kappa',\theta')\, \kappa'{\rm d }\kappa'{\rm d}\theta'. 
\end{align}  
Suppose that $G^{(i)}$ is a random linear combination of steerable Fourier--Bessel basis functions
\begin{align}
    F_{l,s}(\kappa,\theta) = f_{l,s}(\kappa) e^{il\theta}, \quad |l|\le L, \quad s\le S_l, 
\end{align}
 where $f_{l,s}$'s are some radial functions. Let $B=\sum_{l=-L}^L S_l$ be the total number of basis functions, then  $G^{(i)}=F W^{(i)}$ where $F\in \mathbb{C}^{d\times B}$ contains the basis functions and  $W^{(i)} \in \mathbb{C}^{B \times 1}$ is a random coefficient vector. Applying this sketch to the empirical moment  $\overline{\cM}^{(2)}$ under the in-plane symmetrization yields
\begin{align}
\label{eqn:moment_sym}
    (G\circ \overline{\cM^{(2)}})_i = \frac{1}{2\pi} \frac{1}{N}\int_0^{2\pi} \sum_{j=1}^N R(\varphi) \hat I_j \hat I_j^* R(-\varphi) G^{(i)}\, \mathrm{d}\varphi
\end{align}
where $R(\varphi)\in \mathbb{C}^{d\times d}$ is the rotation operator acting on images  with rotation angle $\varphi \in [0,2\pi]$. Because the Fourier--Bessel basis is steerable, i.e.
\begin{align*}
    F_{l,s}(\kappa,\theta+\varphi)=f_{l,s}(\kappa) e^{il(\theta+\varphi)} = F_{l,s}(\kappa,\theta)e^{il\varphi},
\end{align*}
the rotation operator satisfies  $R(-\varphi)F=FD(\varphi)$, where $D(\varphi) \in \mathbb{C}^{B\times B}$ is a diagonal matrix with entries $e^{il\varphi}, |l|\le L$. Then every entry of the integrand in~\eqref{eqn:moment_sym} is a function of $\varphi$ that can be represented using Fourier basis, such that 
\begin{align*}
    \frac{1}{N}\sum_{j=1}^N R(\varphi) \hat I_j \hat I_j^* R(-\varphi) G^{(i)} = \frac{1}{N}\sum_{j=1}^N  \hat I_j  \left(  \hat I_j^*   \left( D(\varphi) F W^{(i)} \right) \right) = A^{(i)} E(\varphi),
\end{align*}
where $E(\varphi) \in \mathbb{C}^{(2L+1)\times 1 }$ is the Fourier basis vector  $[e^{-iL\varphi},e^{-i(L-1)\varphi},\ldots,e^{iL\varphi}]^T$ and $A^{(i)}\in \C^{d\times (2L+1)}$ is a  coefficient matrix.  The computational cost of forming $A^{(i)}$ naively  is $\cO(NdL)$.  Once $A^{(i)}$ is formed, the integral
\begin{align}
    (G\circ \overline{\cM^{(2)}})_i = \frac{1}{2\pi} \int_0^{2\pi}  R(\varphi) A^{(i)} E(\varphi)\, \mathrm{d}\varphi
\end{align}
can be computed in $\cO(Q dL)$ using numerical quadrature with $Q$ nodes on $[0,2\pi]$. With $r$ sketch functions, the total computational cost is  $\cO(NdLr+QdLr)$, and the method is easily parallelizable. In contrast, data augmentation via explicit image rotation has cost  $\cO(NdQr)$, which is more expensive when $Q>L$. In practice, we can choose moderate $L$  for \textit{ab  initio} modeling to extract smooth features, while a large $Q\approx \sqrt{d}$ is required for fully sampling in-plane rotations and denoising.

This approach naturally accommodates radially symmetric CTFs and extends to third order moments. The efficient construction of the  $A^{(i)},i=1,\ldots,r$ and its analogy for the third order moment remains to be open, which we leave for future work.

\subsection{Important optimization details}
\label{sec:numerical_optim}
In this section, we discuss the optimization method used to solve the SubspaceMoM problem in~\eqref{eqn:subspaceMoM}. In Section \ref{sec:choicebasis},  we specify our choice of  basis functions used for representing the volume and rotational density. We require $\hat V$ to be the Fourier transform of a real-valued function and $\mu$ to be a probability density. Once the basis is chosen, these requirements correspond to constraints on the coefficients of the basis. We explain how to construct and enforce these constraints within the optimization in Section \ref{sec:constraints}. Subsequently, in Section \ref{sec:twostageopt}, we use a sequential moment matching   approach  \cite{Katsevich2020LikelihoodMA}, which in practice performs better than directly solving the full objective function in~\eqref{eqn:subspaceMoM}.

\subsubsection{Choice of basis}
\label{sec:choicebasis}
Since we have assumed that the volume is band-limited by $c=1/2$, we can represent its Fourier transform using the spherical Bessel basis $\phi_{l,m,s}$, defined by
\begin{align}
\label{eqn:sphFB}
\phi_{l,m,s}(\kappa,\theta,\varphi) = C_{l,s} \cdot f_{l,s}(\kappa)\cdot  Y_l^m(\theta,\varphi).
\end{align}
where the angular part is given by the spherical harmonics
\begin{align}
 Y_l^m(\theta,\varphi) = (-1)^m \sqrt{\frac{2l+1}{4\pi}\frac{(l-|m|)!}{(l+|m|)!}} P_l^{|m|}(\cos\theta) e^{im \varphi}
\end{align} 
with $l\ge 0$, $-l \le m \le l$, where $P_l^m$  are the associated Legendre functions \cite{brychkov2008handbook}. The spherical harmonics form an orthogonal basis with respect to the measure $\sin \theta\,{\rm d}\theta\,{\rm d}\varphi$ on the unit sphere.
The radial part is given by the following rescaled spherical Bessel function
\begin{align*}
    f_{l,s}(\kappa)=j_l\left (\frac{ \rho_{l,s}}{c} \kappa\right) , \quad c=1/2
\end{align*}
where $j_l$ denotes the $l$-th degree  spherical Bessel function \cite{brychkov2008handbook}, and $\rho_{l,s}$ is the $s$-th zero of $j_l$. The set $\{f_{l,s}\}$ forms an orthogonal basis on $[0,1/2]$ with respect to the measure $\kappa^2\,{\rm d}\kappa$. The constants $C_{l,s}$ are chosen such that $\phi_{l,m,s}$ are orthonormal with respect to the measure $\kappa^2 \sin\theta\, {\rm d}\kappa\, {\rm d}\theta\, {\rm d}\varphi$. 

The volume $\hat V$ is approximated using a truncated spherical Bessel expansion as follows:
\begin{align}
\label{eqn:volumesphBessel}
\hat V_a(\kappa,\theta,\varphi) = \sum_{l=0}^L \sum_{s=1}^{S(l)}\sum_{m=-l}^l
    a_{l,m,s} \phi_{l,m,s}(\kappa,\theta,\varphi).
\end{align}
where $a_{l,m,s}\in \C$ are  expansion coefficients, $L>0$ is the truncation limit, and $S(l)$ is the radial truncation limit for each $l$. In practice, since the volume is observed through discretized 2-D projections, only finitely many coefficients can be determined. We adopt the truncation criterion from \cite{Zhao2012FourierBesselRI}, where $S(l)$ is chosen as the largest integer $s$ satisfying
\begin{align}
\label{eqn:Nyquist_rule}
\frac{\rho_{l,s+1}}{\pi } \le A ,
\end{align}
where $A>0$ denotes the radius of a ball containing the images.   Therefore, the  number of coefficients used to represent the volume satisfies $|\mathcal{B}_V|=\cO\left ( \sum_{l=0}^L S(l) (2l+1)\right)$. Here typically $S(l)=\cO(L)$ because $S(l)$ is a decreasing function of $l$ and for sufficiently large $l$, we have $S(l)<2l+1$~\cite{Bhamre2017AnisotropicTF,Zhao2012FourierBesselRI}, so we can conclude that $|\mathcal{B}_V|=\cO(L^3).$

We approximate the rotational density using the Wigner D-matrix  $D_{u,v}^p(R)$ for $p\ge 0$ and $-p\le u,v \le p$. 
For each $p$, the $(2p+1)\times (2p+1)$ matrix $(D_{u,v}^p(R))_{u,v=-p}^p$ is orthonormal with the $(u,v)$-th  entry given by
\begin{align}
\label{eqn:wigner-D}
D_{u,v}^p(R) = e^{-iu\alpha} d_{u,v}^p(\beta) e^{-iv\gamma},
\end{align} 
where  $\alpha,\beta,\gamma$ are the Euler angles of $R $ and $d_{u,v}^p$ denotes the (small) Wigner d-matrices \cite[p 341]{chirikjian2016harmonic}. The entries themselves are orthogonal to each other and  can be used to form an orthonormal basis of $L^2(\mathcal{SO}(3))$. In this paper, we only consider densities that are invariant to in-plane rotations (in-plane uniform), as introduced in Section~\ref{sec:steerablePCA}. In this case, only the coefficients of the  functions in~\eqref{eqn:wigner-D} with $v = 0$ are nonzero~\cite{Sharon_2020}, so  the viewing direction density $\mu_b$ in~\eqref{eqn:densityexpand} can be expressed as
\begin{align}
\label{eqn:viewingdensityWignerD}
\mu_b(R) = \sum_{p=0}^P\sum_{u=-p}^p b_{p,u} D_{u,0}^p(R)
\end{align}
where $b_{p,u}\in \C$ are the expansion coefficients.
 Under this choice, the number of coefficients to represent the rotational density satisfies $|\mathcal{B}_\mu|=\cO\left( \sum_{p=0}^P (2p+1)\right)=\cO\left(P^2\right).$
In fact, the entries $D_{u,0}^p(R)$  are proportional to spherical harmonics, given by
\begin{align}
\label{eqn:wignerd_2_sph}
D_{u,0}^p(R)   =\sqrt{\frac{4\pi}{2p+1}}  Y_p^{u*}(\beta,\alpha) ,
\end{align}
where $(\alpha,\beta,\gamma)$ is the Euler angle representation of $R\in \mathcal{SO}(3)$. This implies that $\mu_b$ can be  determined by a probability distribution on the unit sphere, which we call the \textit{viewing direction distribution}. We define the normalized spherical harmonic basis:
\begin{align}
\label{eqn:spheicalbasis}
\zeta_{p,u}(\alpha,\beta) =  2\pi \sqrt{\frac{4\pi}{2p+1}}  Y_p^{u*}(\beta,\alpha), 
\end{align}
and write the viewing direction density $\nu_b$  as
\begin{align}
\label{eqn:viewing_density}
\nu_b(\alpha,\beta)=\sum_{p=0}^P\sum_{u=-p}^p b_{p,u} \zeta_{p,u}(\alpha,\beta). 
\end{align}
Consequently, we can factorize the rotational density function $\mu_b$ into the product of the above spherical density and a uniform in-plane density, given by
\begin{align}
\label{eqn:factorization_of_haar}
\mu_b(R(\alpha,\beta,\gamma)) = \frac{1}{2\pi}\nu_b(\alpha,\beta)
\end{align}
since the uniform in-plane density is  $\frac{1}{2\pi}$ for all $\gamma \in [0,2\pi]$. 

This representation  can be further simplified by enforcing invariance to in-plane reflection since we can augment the dataset by including copies of all projection images reflected through the origin. This can be realized by enforcing the spherical harmonic coefficients of odd degree in \eqref{eqn:viewing_density} to zero, i.e., $b_{p,u}=0$ if $p$ is odd~\cite{Zhang_Mickelin_Kileel_Verbeke_Marshall_Gilles_Singer_2024}. 

With this choice of basis functions, we can compute the mapping from the coefficients to the moments in~\eqref{eqn:newM1}-\eqref{eqn:newM3},
$
(a,b) \mapsto (\cM_c^{(1)}[ a, b](\Xi),\,\, \cM_c^{(2)}[a, b](\Xi^{\otimes 2}),\,\, \cM_c^{(3)}[a,b](\Xi^{\otimes 3})),
$
using numerical quadrature. Methods for constructing  numerical quadrature on $\mathcal{SO}(3)$ or the unit sphere can be found in \cite{Grf2009SamplingSA,Grf2013EfficientAF}, or in Appendix \ref{sec:constructSO3quadratures}. 

\begin{remark}
Our approach is not limited to the  basis functions  chosen in Section \ref{sec:choicebasis}. It can be easily extended to  other basis functions. For example, the spheroidal wave functions (PSWF) \cite{Greengard2018GeneralizedPS}  can also be used to model band-limited volumes. Other functions suitable for representing molecules include wavelets \cite{10.5555/130655,5872791,pnas.2216507120}, Gaussian mixtures \cite{doi:10.1137/22M1498218,Chen2021,Zickert_2022,CHEN2023168014,chen2023improving,KAWABATA20181,KAWABATA20084643},  and neural networks \cite{Zhong2021,Zhong2019ReconstructingCD,khoo2023deep,DONNAT2022107920,GIRI2023102536}.  The basis for rotational density functions  can for instance also be chosen as a mixture of von Mises–Fisher distributions \cite{ROSSI202241}, or a mixture of point masses on the unit sphere \cite{khoo2023deep}.  
\end{remark}

\subsubsection{Constraints for numerical optimizations}
\label{sec:constraints}
In this section, we describe several necessary constraints to ensure that solutions to the SubspaceMoM problem~\eqref{eqn:subspaceMoM} yield valid reconstructions.

Firstly, to guarantee that the reconstructed volume (i.e., the inverse Fourier transform of $\hat{V}_a$) is real-valued, the coefficients must satisfy
\begin{align*}
     a^*_{l,m,s} (-1)^m = a_{l,-m,s}(-1)^l.
\end{align*}
This condition can be enforced by the following re-parameterization:
\begin{align}
\label{eqn:vol_realness}
a_{l,m,s} = \begin{cases}
\tilde a_{l,m,s}-(-1)^{l+m}i \tilde a_{l,s,-m}, & m>0, \\
i^l \tilde a_{l,m,s}, & m=0,\\
i \tilde a_{l,m,s}+(-1)^{l+m} \tilde a_{l,s,-m},& m<0,
\end{cases}
\end{align}
for $\tilde a_{l,m,s}\in \R$.

Second, to ensure that the viewing direction density $\nu_b$ is real-valued, we require  
\begin{align*}
b_{p,u} = (-1)^u b_{p,-u}^*. 
\end{align*}
This condition can be enforced by the re-parameterization:
\begin{align}
\label{eqn:density_realness}
b_{p,u} =
\begin{cases}
    \tilde b_{p,u} +(-1)^u i \tilde b_{p,-u}, & p>0,\\
    \tilde b_{0,0}, & p=0,\\
    -i \tilde b_{p,u} + (-1)^u \tilde b_{p,-u}, & p<0,
\end{cases}
\end{align}
where $\tilde b_{p,u}\in \R$. For a discussion and derivation of the constraints in \eqref{eqn:vol_realness} and \eqref{eqn:density_realness}, we refer to \cite{Sharon_2020}.

Third, we require $\nu_b$ to be a probability density function, i.e.,
\begin{align*}
    \int_0^{\pi}\int_0^{2\pi} \nu_b(\alpha,\beta) \sin(\beta) \,{\rm d}\alpha\, {\rm d}\beta = 1.
\end{align*}
Since the spherical harmonics satisfy
\begin{align*}
    \int_0^\pi \int_0^{2\pi} Y_p^u(\beta,\alpha) \sin(\beta) \,{\rm d}\alpha\, {\rm d}\beta =  
    \begin{cases}
        \sqrt{4\pi}, & p=0,u=0,\\ 
        0, & \mbox{otherwise},
    \end{cases}
\end{align*}
we can simply fix  $b_{0,0} =  1$  and only solve for the rest of the parameters. 

Finally, we require the viewing direction density to be non-negative, i.e. 
\begin{align*}
    \nu_b(\alpha,\beta) = \sum_{p=0}^P\sum_{u=-p}^p b_{p,u} \zeta_{p,u}(\alpha,\beta) \ge 0, \quad \forall \alpha \in [0,2\pi], \beta \in [0,\pi]
\end{align*}
and we relax this requirement to a set of linear inequality constraints
\begin{align}
\label{eqn:inequalitiesconstraints}
    \sum_{p=0}^P \sum_{u=-p}^p b_{p,u} \zeta_{p,u}(\alpha_i,\beta_i) \ge 0, \quad i=1,\ldots,n
\end{align}
where $\{(\alpha_i,\beta_i)\}_{i=1}^n$  represents a set of  collocation points on the unit sphere.  To enforce reflection invariance, we can further  impose
\begin{align}
\label{eqn:reflection_invariance}
    b_{p,u} = 0, \quad -p\le u \le p,  \quad \text{ for all odd } p.
\end{align}
To realize this, we simply set these parameters to zero and  remove them from the set of optimization parameters. 

The final optimization problem involves a loss function that is a polynomial of the coefficients $(a,b)$ of degree up to six in $a$ and up to two in $b$, subject to linear inequality constraints on $b$ arising from~\eqref{eqn:inequalitiesconstraints}.

\subsubsection{Sequential moment matching} 
\label{sec:twostageopt}

Although the third-order moment contains important information and can substantially improve reconstruction quality, we observe in practice that directly solving the full moment-matching problem in \eqref{eqn:subspaceMoM} via numerical optimization may lead to unwanted local minima. In~\cite{lo2024methodmomentsestimationnoisy} where  a method of moments approach for reconstructing high-dimensional signals without projection is proposed, a careful initialization is critical when optimizing over the third-order moment, due to the nonconvex nature of the objective function.

Motivated by the \emph{frequency marching} approach in \cite{doi:10.1137/16M1097171}, which incrementally reconstructs the Fourier volume from low to high frequencies under an iterative refinement framework, we use a \emph{sequential moment matching} approach. Specifically, we first fit only the first   subspace moment, then incorporate the second subspace moment, and finally match all three  subspace moments. This progressive approach appears to avoid local trapping and stabilizes the optimization.

More concretely, we begin by solving the first stage optimization problem:
\begin{align}
\label{eqn:firststageoptm}
    \min_{(a,b)\in\mathcal{C}} \lambda_1 \Vert \cM_c^{(1)}[a,b](\Xi)-\overline M_c^{(1)}\Vert^2_F 
\end{align}
where the coefficients for the volume in~\eqref{eqn:vol_realness} are initialized using  i.i.d. Gaussian entries  and the coefficients for the rotational density in~\eqref{eqn:density_realness} are initialized  using the coefficients of a randomly generated non-uniform distribution. The resulting estimates are then used as the initial point to the second stage optimization problem:
\begin{align}
\label{eqn:secondstageoptm}
    \min_{(a,b)\in \mathcal{C}} \lambda_1 \Vert \cM_c^{(1)}[a,b](\Xi)-\overline M_c^{(1)}\Vert^2_F + \lambda_2 \Vert \cM_c^{(2)}[ a, b] (\Xi^{\otimes 2})-\overline M_c^{(2)}\Vert^2_F.
\end{align}
Finally, the output of \eqref{eqn:secondstageoptm} is used to initialize the third stage optimization problem involving all three subspace moments in~\eqref{eqn:subspaceMoM}. For all three optimization problems, we use the implementation of the sequential least squares programming (SLSQP) algorithm~\cite{NoceWrig06} provided by the SciPy optimization package.  We provide the objective function itself and its first-order gradient as inputs to this function. We set the objective function tolerance parameter $\texttt{ftol=1e-10}$, and the maximal number of iteration $\texttt{maxiter=1000}$, and we use the default values for the other parameters. For the numerical experiments in Section~\ref{sec:experiments}, we observe all optimizations exit successfully before reaching the maximal number of iterations.

The effectiveness of this sequential approach may be justified by analogy to the method of moments for estimating Gaussian mixture models. In~\cite{Katsevich2020LikelihoodMA}, it was shown that in the low SNR regime, estimators obtained via sequential moment matching will approximate the maximum likelihood estimator. This insight offers a potential explanation for the improved behavior observed in our setting.

\section{Numerical Results}
\label{sec:experiments}
In this section, we present numerical results that demonstrate the performance of SubspaceMoM in two different settings. In Section \ref{sec:recoverfromMoMs}, we provide examples of reconstructing a smooth volume using analytical subspace moments. In Section \ref{sec:recoverEstimatedMoMs},  we reconstruct volumes from subspace moments estimated from simulated noisy and CTF-affected 2-D projections, emulating the setting of performing \textit{ab initio} modeling in cryo-EM. The molecule data  we use is downloaded from the Electron Microscopy Data Bank (\href{https://www.ebi.ac.uk/emdb/}{https://www.ebi.ac.uk/emdb/}). All experiments were conducted  on a 16-core Linux machine with 12th Gen Intel(R) Core(TM) i9-12900 CPUs running at 2.40GHz and  125 GB of memory.  

The current implementation is written in Python 3 and does not utilize parallelization. The most computationally expensive parts are the tensor contractions involved in moment sketching and formation, as well as the evaluation of the gradient of the objective function. These operations are vectorized using the \texttt{einsum} function from the CPU version of the JAX library. In practice, we observe that \texttt{einsum} executes primarily in a single core. We expect that substantial speedups could be obtained by migrating the implementation to a GPU-based architecture.

\begin{figure}[!ht]
 \centering
 \begin{tikzpicture}[scale=1]
  \node[inner sep=0] at (-18,-3) {\includegraphics[width=0.4\textwidth]{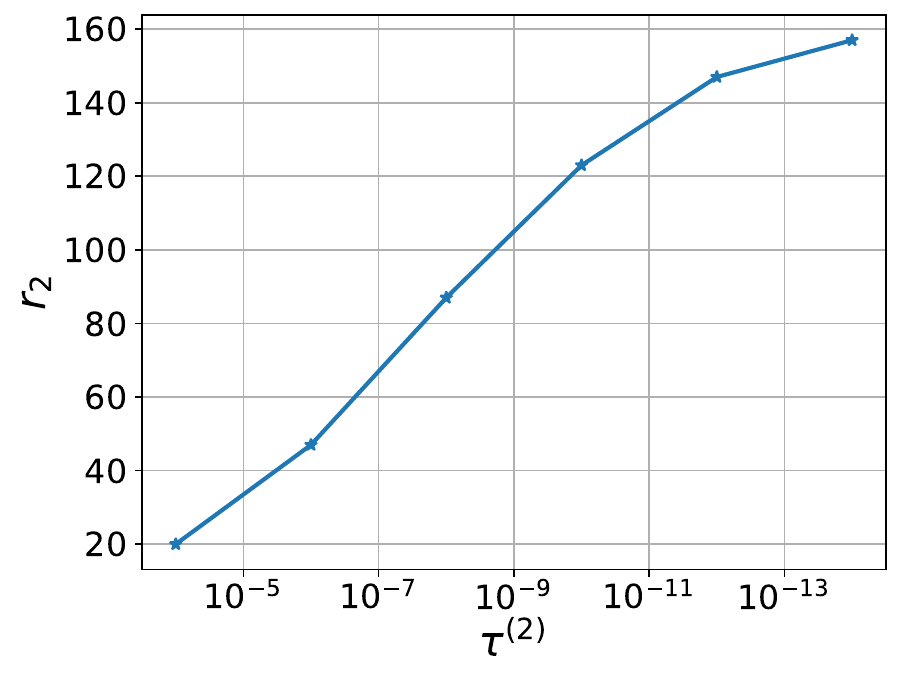}}; 
  \node[inner sep=0] at (-11,-3) {\includegraphics[width=0.4\textwidth]{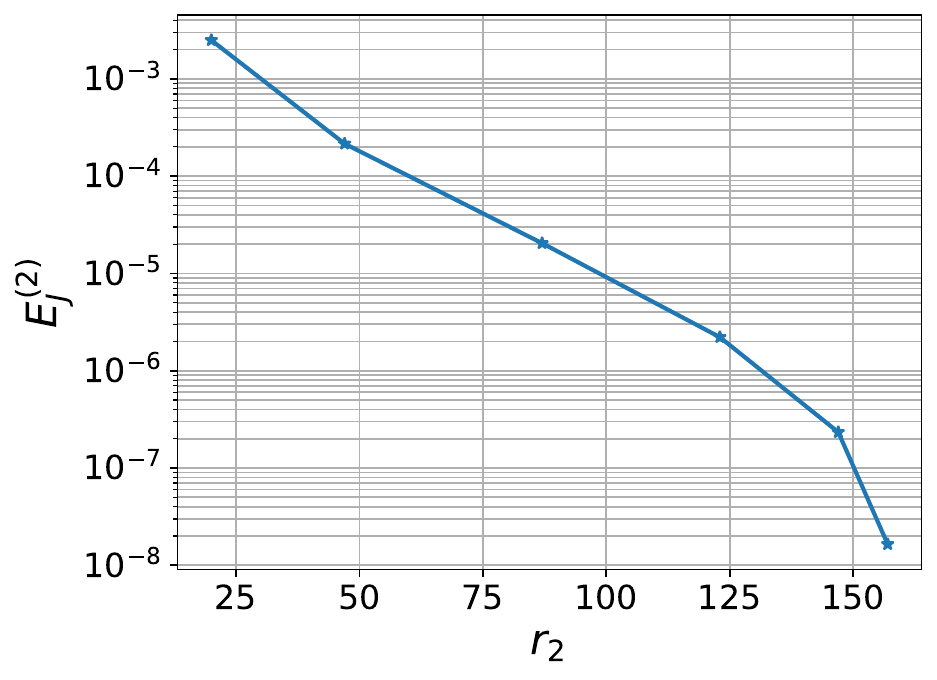}}; 
  \node[inner sep=0] at (-18,-8.5) {\includegraphics[width=0.4\textwidth]{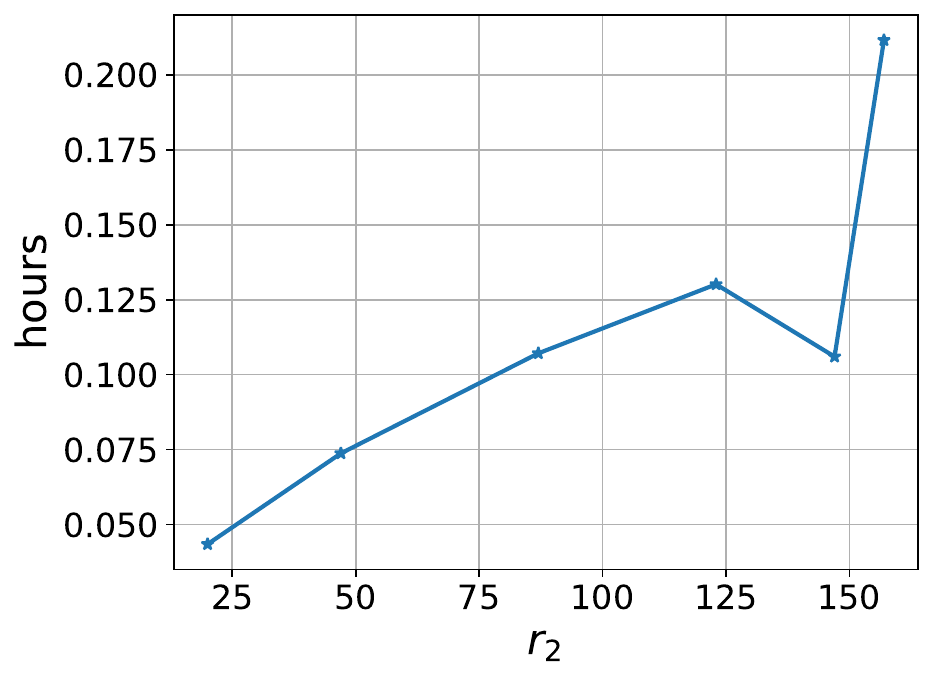}}; 
  \node[inner sep=0] at (-11,-8.5) {\includegraphics[width=0.4
\textwidth]{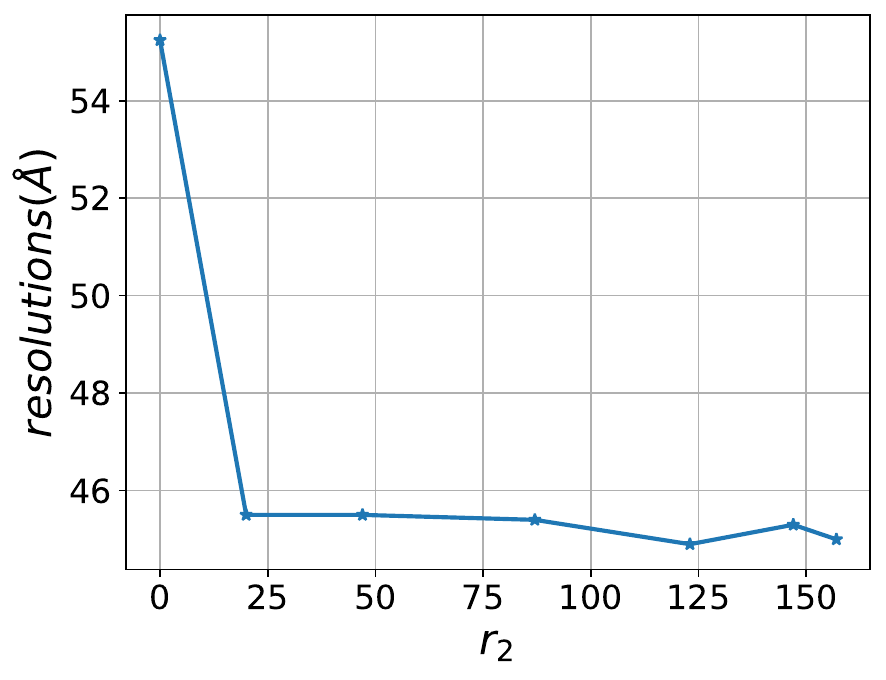}}; 
  \node[inner sep=0] at (-17.5,-5.5) {\footnotesize (A)};
  \node[inner sep=0] at (-10.5,-5.5) {\footnotesize (B)};
  \node[inner sep=0] at (-17.5,-11) {\footnotesize (C)};
  \node[inner sep=0] at (-10.5,-11) {\footnotesize (D)};
 \end{tikzpicture}
 \caption{(A) The dimension parameter $r_2$ for the second subspace moment as a function of the threshold value $\tau^{(2)}$, obtained in Section~\ref{sec:thirdmomenteffect}. Specifically, we obtain $r_2=20,  47,  87, 123, 147, 157$. (B) The relative errors evaluated on a $60\times60$ sub-matrix of the second moment, as a function of  the dimension parameter of the second subspace  moment $r_2$. (C) The running time (in hours) for completing the optimization~\eqref{eqn:secondstageoptm} that matches the first two subspace moments in the second stage as a function of the dimension parameter of the second subspace  moment $r_2$. (D) The resolutions (in angstrom$\A$) compared to the expanded ground truth, obtained from the first two (subspace) moments  as a function of the dimension parameter of the second subspace  moment $r_2$. }
 \label{fig:second_MoM_experiment}
\end{figure}

To evaluate the reconstruction resolution of our method, an alignment is necessary as the moments are invariant to orthogonal transformations of the volume \cite{Zhang_Mickelin_Kileel_Verbeke_Marshall_Gilles_Singer_2024}. We compute the alignment  using the   softwares~\cite{singer2023alignment,Harpaz_Shkolnisky_2023}. After alignment, we compute the \textit{Fourier-Shell correlation} (FSC) between the reconstruction and the ground truth, which is a standard metric for measuring reconstruction resolutions in cryo-EM. On a thin shell with radius $\kappa>0$ and a  small width $\delta \kappa>0$, i.e., $\{\xi \in \R^3 : \kappa \le \Vert \xi \Vert \le \kappa+\delta \kappa \}$, the FSC between two Fourier volumes is given by: 
\begin{align}
\label{eqn:fsc}
{\rm FSC}(\kappa) = \frac{\sum_{\kappa \le \Vert \xi\Vert \le \kappa+\delta \kappa} \hat V_1(\xi) \hat V_2^*(\xi)}{\sqrt{\sum_{\kappa \le \Vert \xi\Vert \le \kappa+\delta \kappa} |\hat V_1(\xi)|^2 \sum_{\kappa \le \Vert \xi\Vert \le \kappa+\delta \kappa} |\hat V_2(\xi)|^2}}.
\end{align}
The resolution is determined by the frequency $\kappa$ at which the FSC drops below a cutoff value, conventionally set to be $\frac{1}{7}$ or $\frac{1}{2}$ \cite{rosenthal2003optimal}. Here, we set the cut-off value to $\frac{1}{2}$, as  we are comparing the reconstructed volume to the underlying ground truth volume rather than two independent reconstructed volumes  from noisy images.

The visualizations of 3-D volumes in this Section are rendered using UCSF Chimera~\cite{https://doi.org/10.1002/jcc.20084}.

\subsection{Reconstruction from analytical moments}
\label{sec:recoverfromMoMs}
In this section, we apply our method to moments computed using the analytical formulas in \eqref{eqn:Mc1}-\eqref{eqn:Mc3}, which allows us to evaluate the performance of our approach under ideal conditions where there is no statistical error, and the target volume and viewing direction distribution can be exactly represented by our bases.  In Section~\ref{sec:thirdmomenteffect}, we show how incorporating the third subspace moment can significantly improve the quality of the reconstruction.  In Section~\ref{sec:truncationlevels}, we demonstrate the feasibility of using  higher truncation limits with inexact quadrature for reconstruction.

\begin{figure}[!ht]
 \centering
 \begin{tikzpicture}[scale=1]
  \node[inner sep=0] at (-18,-3) {\includegraphics[width=0.4\textwidth]{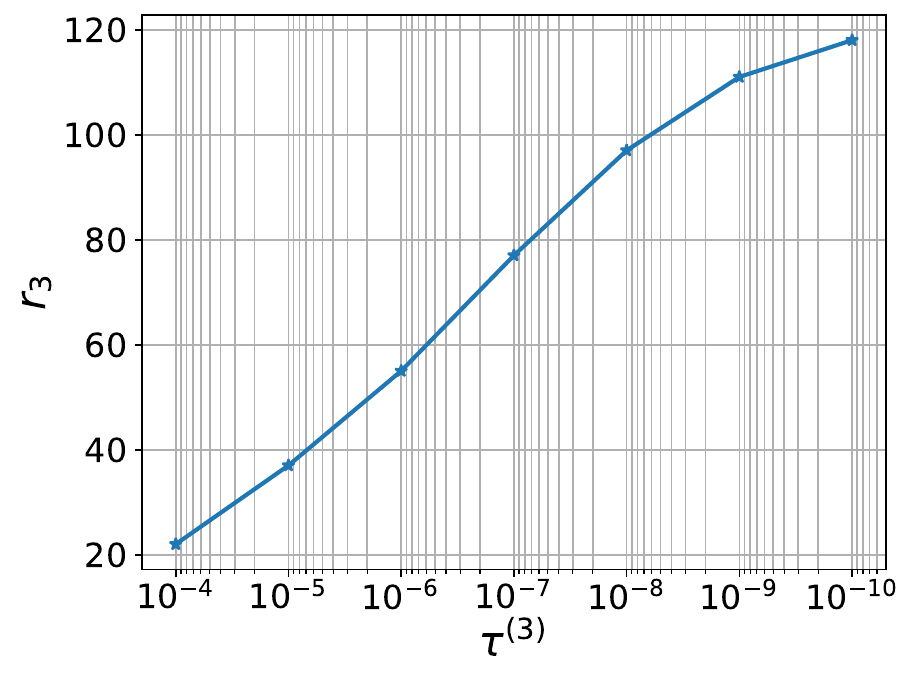}}; 
  \node[inner sep=0] at (-11,-3) {\includegraphics[width=0.4\textwidth]{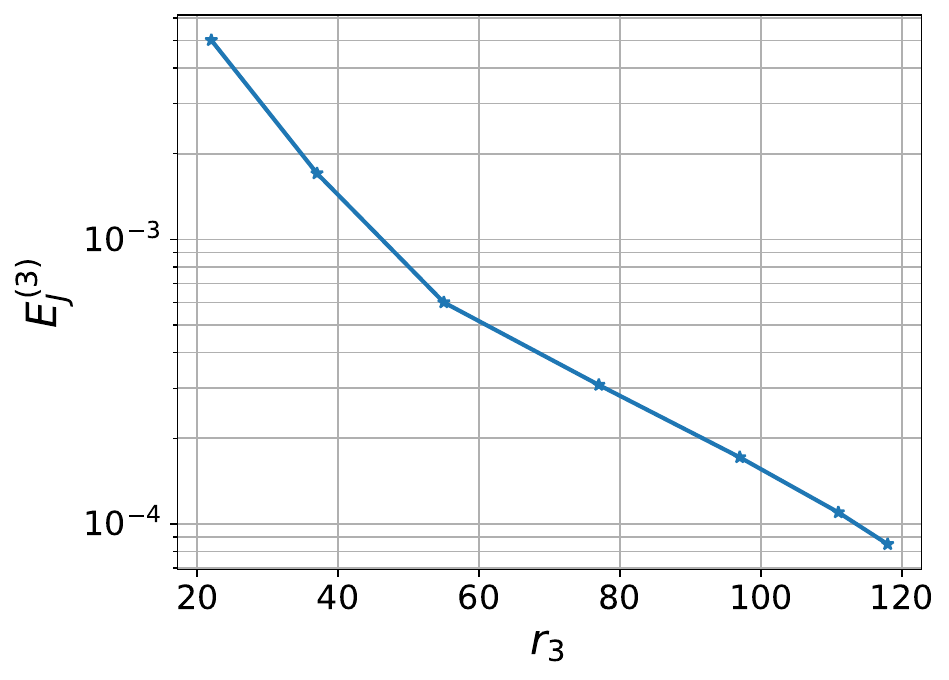}}; 
  \node[inner sep=0] at (-18,-8.5) {\includegraphics[width=0.4\textwidth]{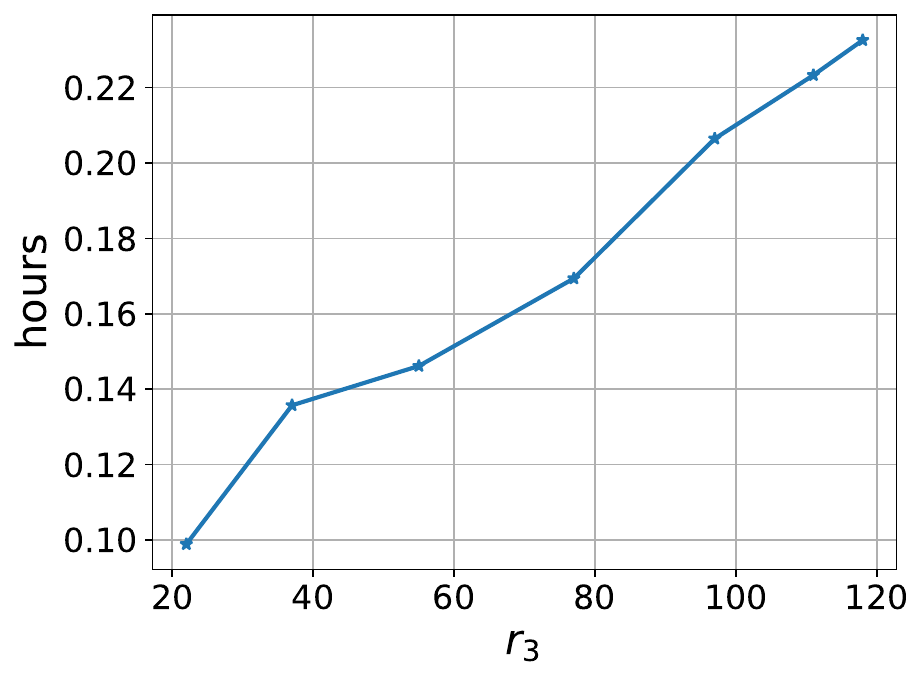}}; 
  \node[inner sep=0] at (-11,-8.5) {\includegraphics[width=0.4
\textwidth]{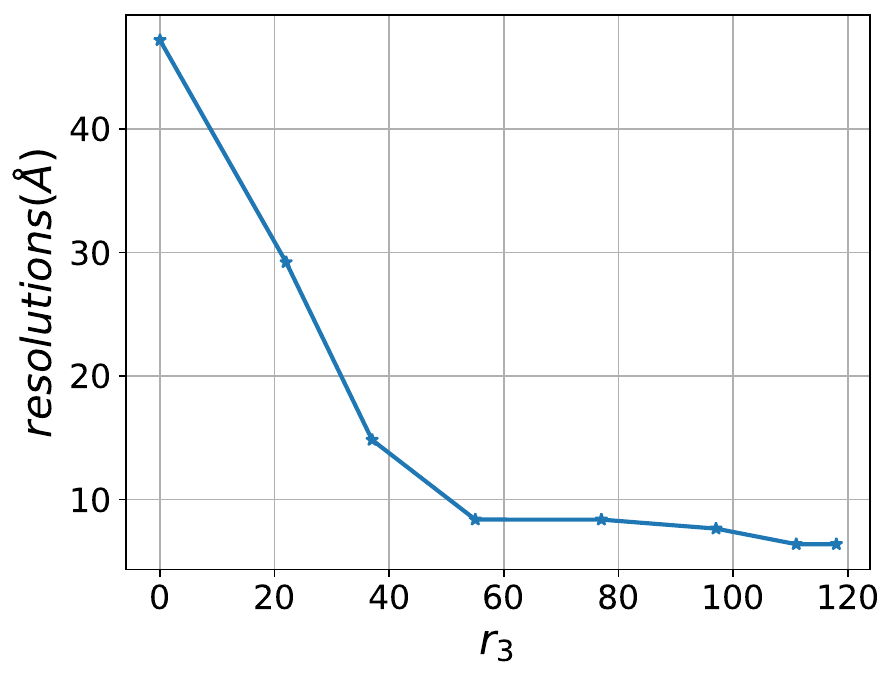}}; 
  \node[inner sep=0] at (-17.5,-5.5) {\footnotesize (A)};
  \node[inner sep=0] at (-10.5,-5.5) {\footnotesize (B)};
  \node[inner sep=0] at (-17.5,-11) {\footnotesize (C)};
  \node[inner sep=0] at (-10.5,-11) {\footnotesize (D)};
 \end{tikzpicture}
 \caption{(A) The dimension parameter $r_3$ for the third subspace moment as a function of the threshold value $\tau^{(2)}$ obtained in Section~\ref{sec:thirdmomenteffect}. Specifically, we obtain $r_3=22, 37, 55, 77, 97, 111, 118$. (B) The relative errors evaluated on a $60\times60\times 60$ sub-tensor of the third moment, as a function of  the dimension parameter of the second subspace  moment $r_3$. (C) The running time (in hours) for completing the optimization~\eqref{eqn:subspaceMoM} in the last stage as a function of the dimension parameter of the third subspace  moment $r_3$. (D) The reconstructed resolutions (in angstrom $\A$) compared to the expanded ground truth, obtained from the first three subspace moments  as a function of the dimension parameter of the third subspace  moment $r_3$. }
 \label{fig:third_MoM_experiment}
\end{figure}

\subsubsection{Importance of the third-order moment}
\label{sec:thirdmomenteffect}
\begin{figure}[!ht]
    \centering
    \begin{tikzpicture}
    \node[inner sep=0] at (0,0) {\includegraphics[width=0.2\textwidth]{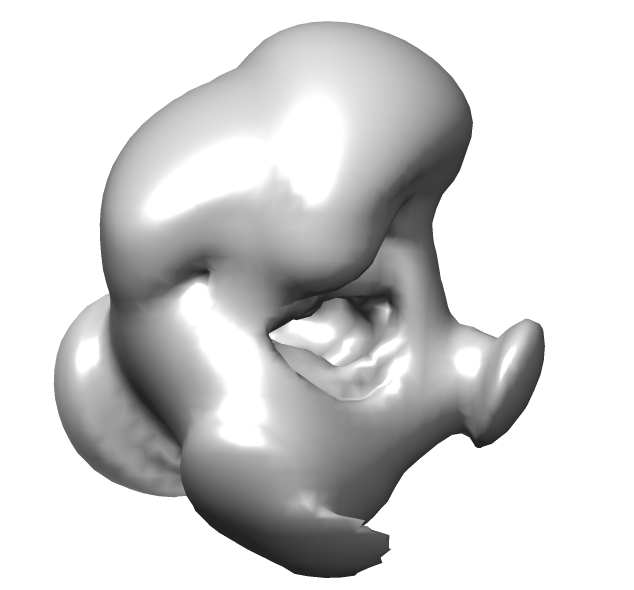}}; 
    \node[inner sep=0] at (4,0) {\includegraphics[width=0.2\textwidth]{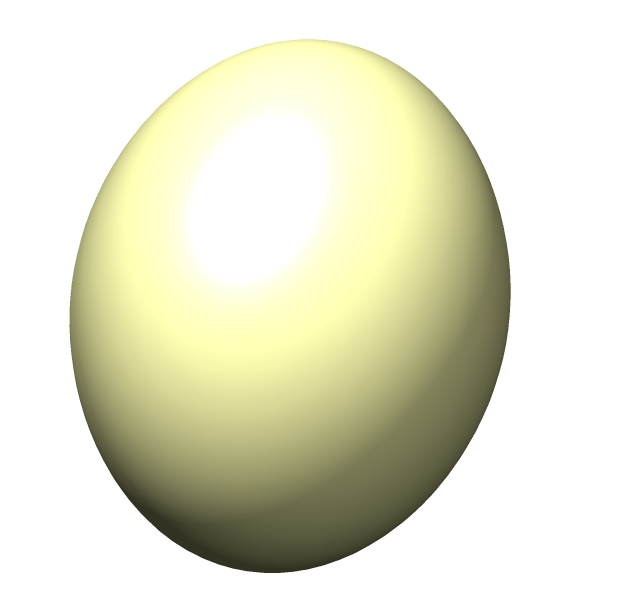}}; 
    \node[inner sep=0] at (8,0) {\includegraphics[width=0.2\textwidth]{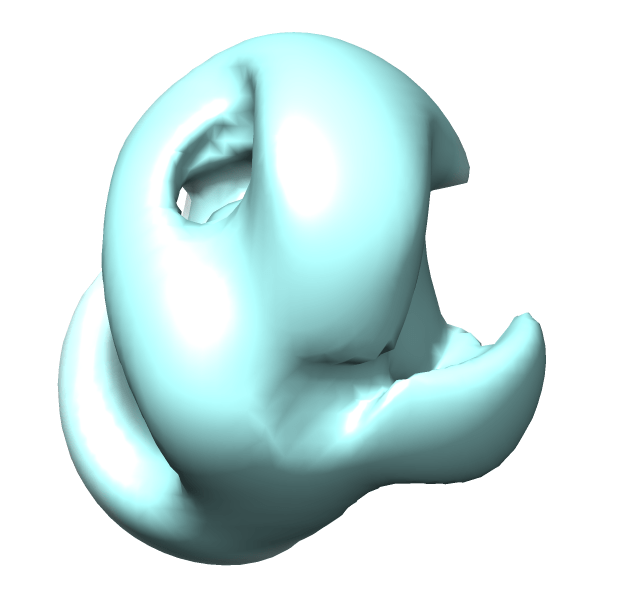}}; 
    \node[inner sep=0] at (12,0) {\includegraphics[width=0.2\textwidth]{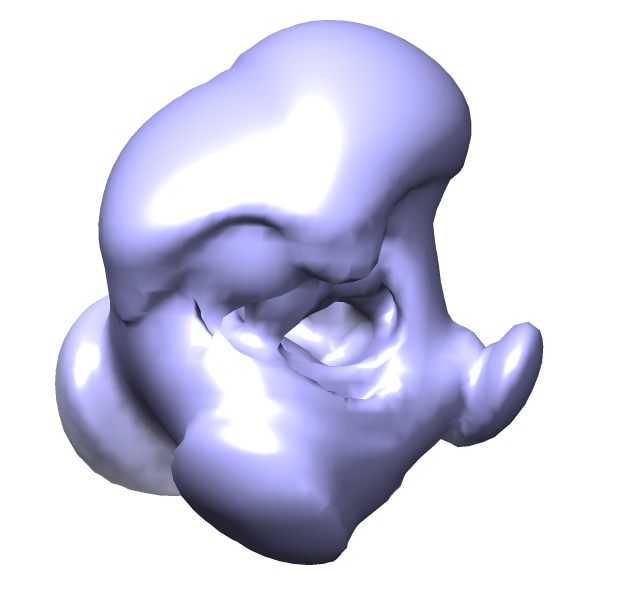}}; 
    \end{tikzpicture}
    \caption{The gray volume represents the ground truth obtained by expanding the downsampled EMD-34948 into the spherical Bessel basis with $L=5$ used in Section~\ref{sec:thirdmomenteffect}. The others are the reconstructions of the sequential moment matching at different stages. Specifically, the yellow volume represents the reconstructed volume obtained from the first subspace moment via solving~\eqref{eqn:firststageoptm}, achieving a resolution of  $54.94 \A$ compared to the expanded ground truth. The blue volume represents the reconstructed volume obtained from the first two subspace moments via solving~\eqref{eqn:secondstageoptm},  achieving a resolution of $45.03 \A$. The purple  volume represents the reconstructed volume obtained from the first three subspace moments via solving~\eqref{eqn:subspaceMoM}, achieving the Nyquist-limited resolution $6.37 \A$. The dimension parameters for the subspace moments are $r_1=400$, $r_2=157$ and $r_3=118$.}
    \label{fig:effect_of_m3_vol}
\end{figure}

\begin{figure}[!ht]
    \centering
    \begin{tikzpicture}[scale=0.99]
        \node[inner sep=0] at (0,0) {\includegraphics[width=0.24\textwidth]{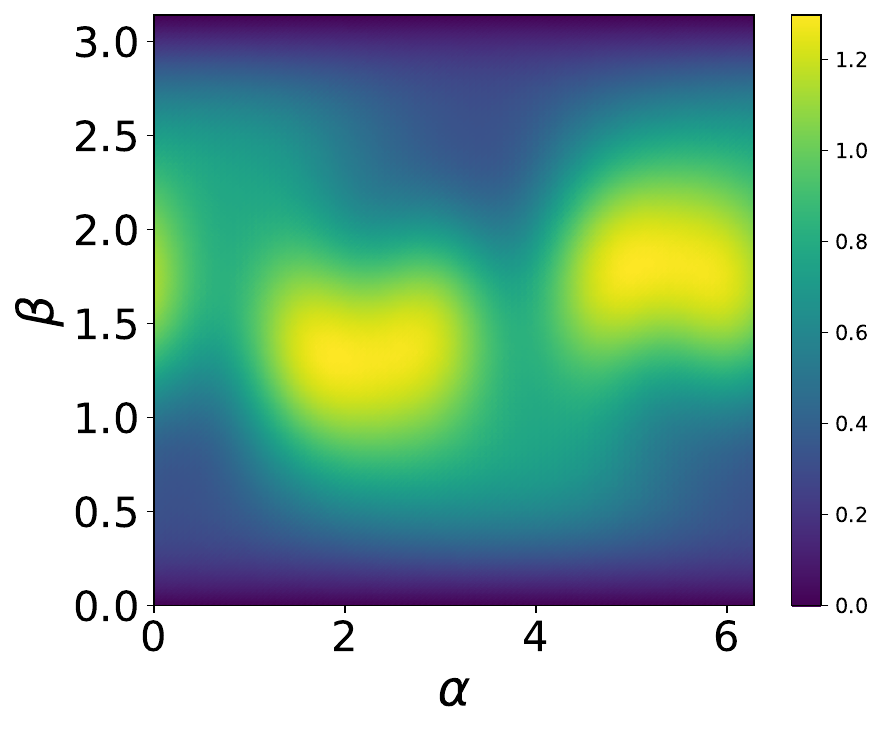}}; 
            
        \node[inner sep=0] at (4,0) {\includegraphics[width=0.24\textwidth]{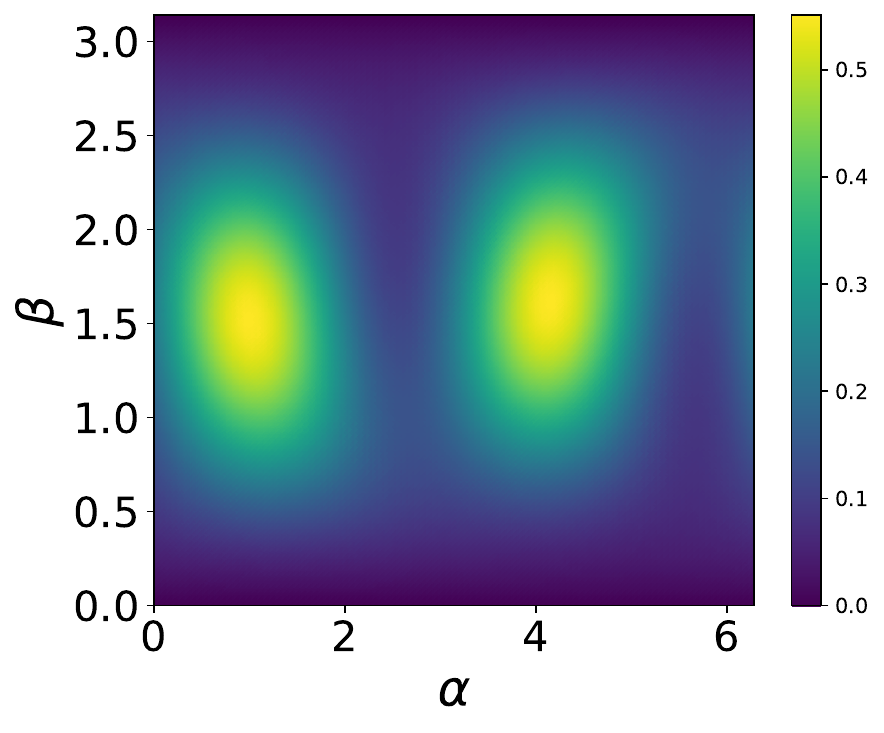}}; 
        
        \node[inner sep=0] at (8,0) {\includegraphics[width=0.24\textwidth]{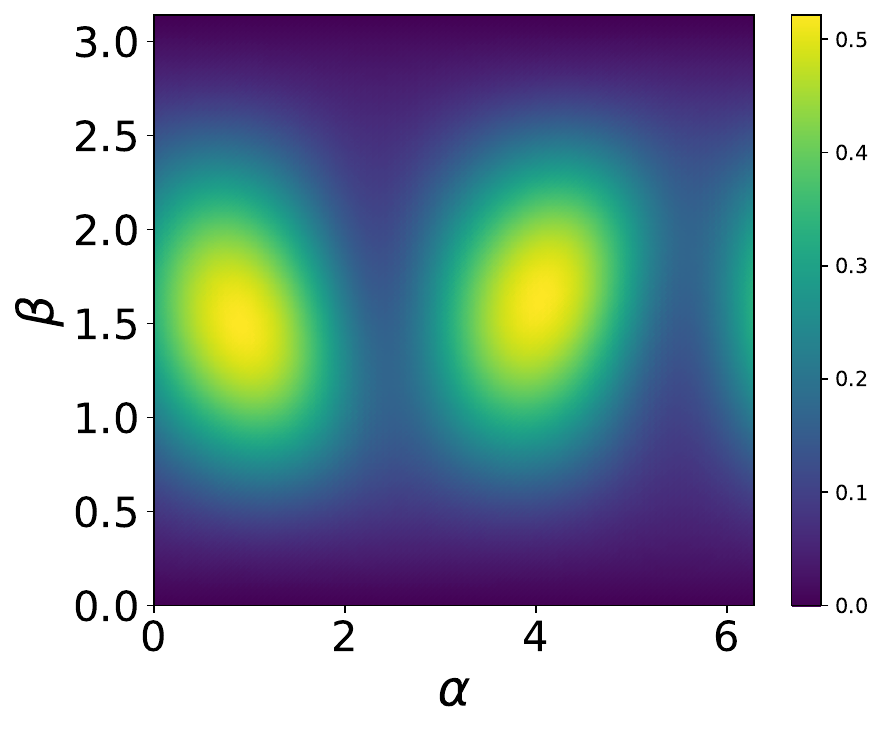}}; 
        \node[inner sep=0] at (12,0) {\includegraphics[width=0.24\textwidth]{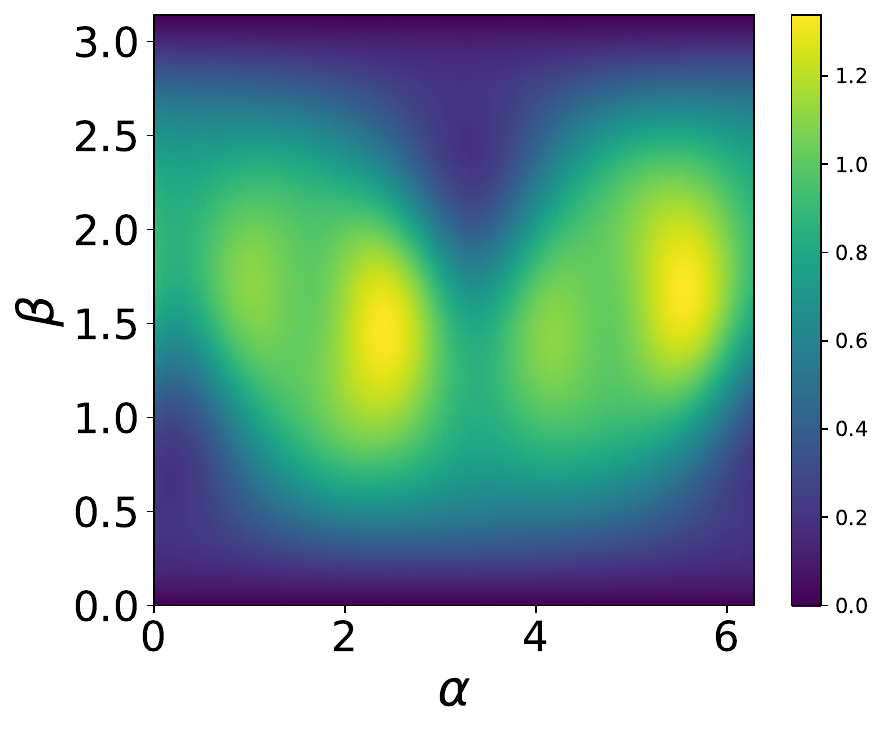}}; 
    \end{tikzpicture}
    \caption{The viewing direction densities plotted as  functions of first two Euler angles $\alpha\in [0,2\pi]$ and $\beta\in [0,\pi]$) that are invariant to in-plane reflection obtained  in Section~\ref{sec:thirdmomenteffect}. The leftmost figure shows the ground truth density.  The second figure shows the estimated density from the first subspace moment via solving~\eqref{eqn:firststageoptm}, which has a relative error $0.4086$ compared to the ground truth. The third figure shows the  estimated density obtained from the first two subspace moments via solving~\eqref{eqn:secondstageoptm}, which has a relative $L^2$ error $0.4080$. The last one shows the estimated density obtained from all three subspace moments via solving~\eqref{eqn:subspaceMoM}, which has a relative $L^2$ error $0.0811$. The dimension parameters for the subspace moments are $r_1=400$, $r_2=157$ and $r_3=118$.
    }
    \label{fig:effect_of_m3_mu}
\end{figure}

In this first example, we demonstrate the necessity of incorporating the third subspace moment. To create a smooth ground truth volume, we use EMD-34948~\cite{Wang2023_GPR132}  consisting of $196 \times 196 \times 196$ voxels, with a physical pixel size of  $1.04 \A$. For computational efficiency, we downsample the volume via Fourier cropping, to $64 \times 64 \times 64$ with physical pixel length  $3.185 \A$.  We expand it into the spherical Bessel basis with a truncation limit of $L=5$ to be our ground truth,  which contains $|\cB_V|=1079$ spherical Bessel coefficients. Compared to the original EMD-34948, this expanded ground truth has  a resolution of  $17.49 \A$. We generate a ground truth viewing direction density by expanding a mixture of von Mises–Fisher distributions (introduced in Appendix~\ref{sec:von-Mises–Fisher}) on the unit sphere using the spherical harmonic basis with a truncation limit of $P=4$, which contains  $|\cB_\mu|=14$ unknowns. We denote the ground truth parameters by the pair $(a_\star,b_\star)$.  The visualizations of the ground truth volume and viewing direction density can be found in \cref{fig:effect_of_m3_vol} and \cref{fig:effect_of_m3_mu}, respectively. With the truncation limits $L$ and $P$, forming the first three moments with $12$-digit accuracy requires $192,729$ and $2080$ quadrature nodes, respectively. These sets of nodes are also used during the optimization process. To enforce non-negativity of the reconstructed viewing direction density, we use a set of $322$ collocation points obtained from the quadrature rule on the unit sphere \cite{Grf2009SamplingSA,Grf2013EfficientAF} in all the following numerical experiments.

Because the first moment is not  informative,  solving the first-stage optimization problem \eqref{eqn:firststageoptm} with the uncompressed first moment only produces a  volume with resolution $55\A$ compared to the expanded ground truth.  To reduce computational cost and memory usage, we compress the first moment using a random projection matrix of size $4096 \times 400$, constructed via orthonormalizing a Gaussian random matrix. Matching this first subspace moment via~\eqref{eqn:firststageoptm} only takes 21 seconds. We use the resulting parameters to initialize the second-stage optimization problem \eqref{eqn:secondstageoptm} that matches the first two subspace moments. We apply the randomized range-finding algorithm described in Section~\ref{sec:findingsubspaces} to $\cM^{(2)}[a_{\star},b_{\star}]$ with a sampling size $s=250$. By varying the threshold value   $\tau^{(2)}= 10^{-4},10^{-6},\ldots, 10^{-14}$, we obtain the second subspace  moments $\cM_{\tau^{(2)}}^{(2)}[a_{\star},b_{\star}]\in \C^{r_2(\tau^{(2)})\times r_2(\tau^{(2)})}$ where $r_2(\tau^{(2)})$ denotes the dimension parameter obtained at $\tau^{(2)}$. As $\tau^{(2)}$ decreases,  the dimension  $r_2$ increases from $20$ to $157$ as plotted in~\cref{fig:second_MoM_experiment}(A). To estimate the errors of those low-rank approximations, we sample a sub-matrix of size $60\times 60$, and evaluate the relative errors \eqref{eqn:relerr_m2} on this sub-matrix. The  errors decay rapidly as $r_2$ increases, which are reported in~\cref{fig:second_MoM_experiment}(B).  The running time for solving the optimization of matching the first two subspace moments \eqref{eqn:secondstageoptm} is shown in~\cref{fig:second_MoM_experiment}(C). As illustrated in~\cref{fig:second_MoM_experiment}(D), the second subspace moment yields only limited improvement, with resolution plateauing at approximately $45 \A$ compared to the expanded ground truth. These results have shown the ill-posedness of the moment-matching problem when only the first two (subspace) moments are used.

We compress the third moment $\mathcal{M}^{(3)}[a_{\star},b_{\star}]$, using the threshold values $\tau^{(3)}$ $=$ $10^{-4}$, $10^{-5}$, $\ldots$, $10^{-10}$ resulting in the third subspace moments $\mathcal{M}^{(3)}_{\tau^{(3)}}[a_{\star},b_{\star}] \in \C^{r_3(\tau^{(3)}) \times r_3(\tau^{(3)}) \times r_3(\tau^{(3)})}$, each associated with a different dimension parameter $r_3(\tau^{(3)})$. As $\tau^{(3)}$ decreases, the dimension parameter $r_3$   increases from $22$ to $118$,  as shown in~\cref{fig:third_MoM_experiment}(A). In~\cref{fig:third_MoM_experiment}(B), we show the relative error of the low-rank approximation evaluated on a $60\times 60 \times 60$ sub-tensor of the third-order moment.  Using the fixed warm start obtained from the first two subspace moments, we solve the optimization problems in \eqref{eqn:subspaceMoM}  with the third subspace moments obtained by varying $\tau^{(3)}$'s. The running time of finishing the last stage optimization in~\eqref{eqn:subspaceMoM} increases with $r_3$, as shown  in~\cref{fig:third_MoM_experiment}(C). Compared to the expanded ground truth, the resolution improves to the Nyquist limit $6.37\A$  when $r_3=118$. The reconstructed resolutions for all values of $r_3$ are presented in~\cref{fig:third_MoM_experiment}(D).
 The visual comparison between the reconstructed volumes from the three stages of the sequential moment matching and the expanded ground truth is provided in~\cref{fig:effect_of_m3_vol}, from which we can observe the reconstruction using the third subspace moment highly resembles the ground truth volume. We evaluate the relative $L^2$ error between the estimated densities of viewing directions from the three stages of the sequential moment matching and the ground truth density. Here,   the relative $L^2$ error of an estimated density $\tilde \mu$  to the ground truth $\mu$, is defined as
\begin{align*}
   \mathrm{err}(\tilde \mu)= \frac{(\int_{0}^{2\pi}\int_0^\pi (\mu(\alpha,\beta)-\tilde \mu(\alpha,\beta))^2\, \sin\beta \, \mathrm{d}\beta\,\mathrm{d}\alpha)^{1/2}}{(\int_{0}^{2\pi}\int_0^\pi \mu(\alpha,\beta)^2\, \sin\beta \, \mathrm{d}\beta\,\mathrm{d}\alpha)^{1/2}}. 
\end{align*}
The computed relative errors of the estimated densities from the three stages are $0.4086$, $0.4080$ and $0.0811$. The visual comparison between the estimated viewing direction densities and the ground truth viewing direction density can be found in~\cref{fig:effect_of_m3_mu}. These numerical results suggest the importance of using the third subspace moment in the method of moments approach.

\begin{figure}[!ht]
    \centering
    \begin{tikzpicture}
        \node[inner sep=0] at (0,0) {\includegraphics[width=0.4\textwidth]{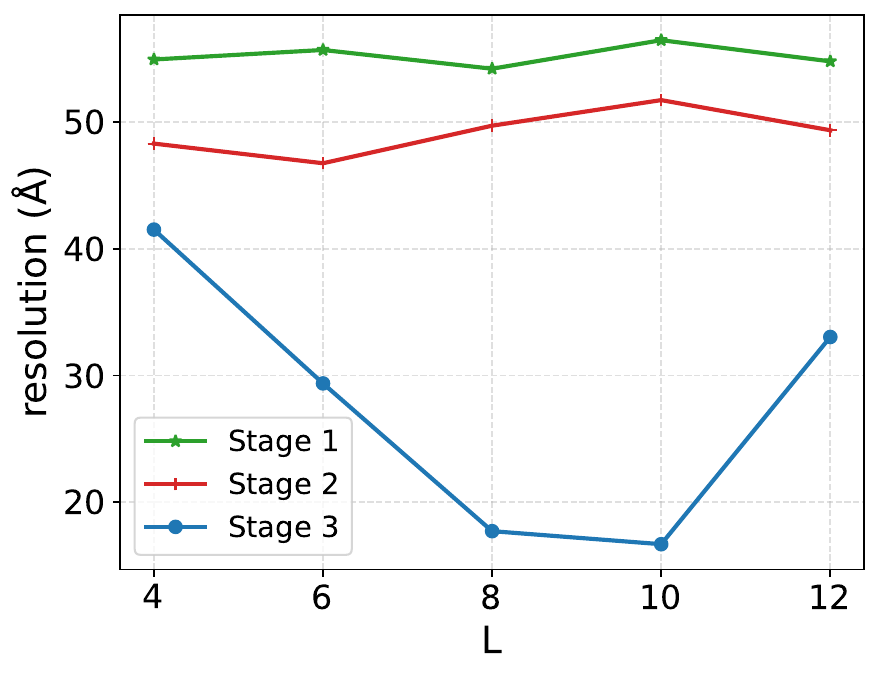}}; 
        \node[inner sep=0] at (8,0) {\includegraphics[width=0.4\textwidth]{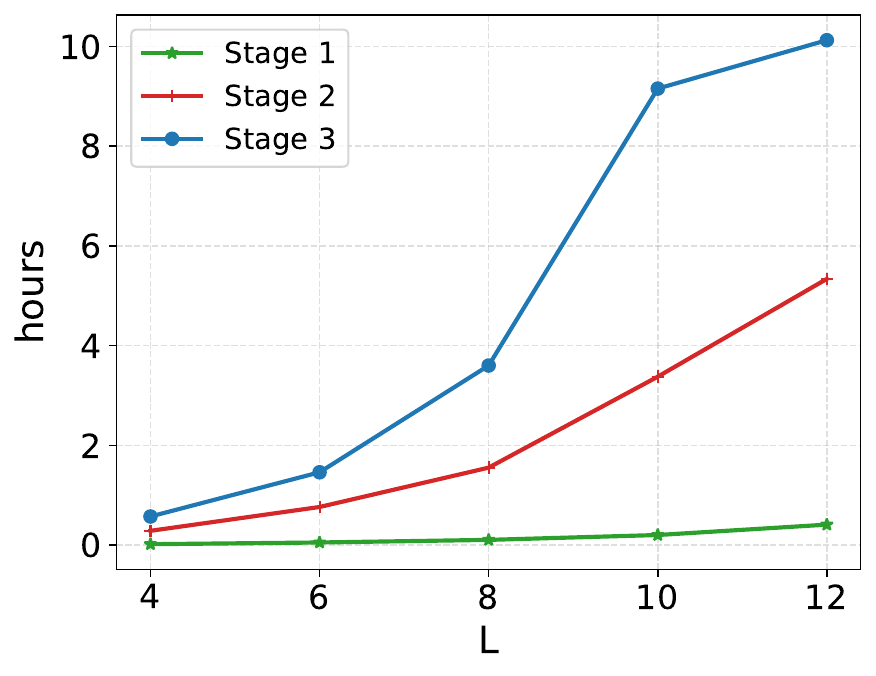}}; 
     \node[inner sep=0] at (0.4,-2.7) {\footnotesize (A)};
    \node[inner sep=0] at (8.4,-2.7) {\footnotesize (B)};
    \end{tikzpicture}
    \caption{(A) The resolutions (compared with the expanded ground truth at $L=12$) achieved by the reconstructed volumes  from different stages of the sequential moment matching, when the reconstructed volumes are expanded at different truncation limits $L$ in Section~\ref{sec:truncationlevels}.  (B) The running time (in hours) required to complete the optimization problems of different stages of the sequential moment matching.}
    \label{fig:effect_of_truncation}
\end{figure}

\subsubsection{Effect of truncation limits}
\label{sec:truncationlevels}
In the next experiment, we use the same ground truth viewing direction distribution as in Section~\ref{sec:thirdmomenteffect} and expand the downsampled EMD-34948 with a higher  truncation limit $L=12$ to be our ground truth volume. The expanded ground truth achieves a resolution of  $9.24 \A$ compared to the original EMD-34948. The analytical subspace moments are formed using the randomized range-finding algorithm in Section~\ref{sec:findingsubspaces} with a sampling size of  $s=250$. The threshold parameters are set to $\tau^{(2)}=10^{-8}, \tau^{(3)}=10^{-6}$, resulting in the subspace  moment dimensions $r_2=213$ and $r_3=103$. The subspace extracted from the second-order moment is also used to compress the first-order moment. For the viewing direction distribution, we aim to reconstruct all 14 spherical harmonic coefficients. We vary the truncation limit $L\in \{4,6,8,10,12\}$ used for the reconstructed volume,  corresponding to $760,1456,2352,3432,4680$ volume coefficients, respectively.   We use a quadrature rule with  $2184$ nodes to integrate all three analytical subspace moments during reconstruction.  This is the maximum number of nodes that can be accommodated by our memory when the reconstructed volume is expanded at $L=12.$  Since the quadrature rule is inexact, it introduces a relative integration error bounded by $10^{-3}$ for the second analytical subspace moment at  $L=8,10,12$ and for the third analytical subspace moment at $L=6,8,10,12$. This error is estimated by forming the analytical subspace moments using both the exact and the inexact quadrature rules, with fixed ground truth parameters as inputs, and computing the relative mean squared error between them, as listed in~\cref{tab:integral_errors}.

\begin{table}[!ht]
\begin{tabular}{|l|l|l|l|l|l|}
\hline
             & $L=4$                 & $L=6$                 & $L=8$                 & $L=10$                & $L=12$                \\ \hline
$1$-st   moment  & $5.58 \times 10^{-8}$ & $6.39 \times 10^{-8}$ & $5.68 \times 10^{-8}$ & $4.90 \times 10^{-8}$ & $4.85 \times 10^{-8}$ \\ \hline
$2$-nd   moment & $1.10 \times 10^{-7}$ & $1.25 \times 10^{-7}$ & $1.11 \times 10^{-7}$ & $1.16 \times 10^{-4}$ & $5.96 \times 10^{-4}$ \\ \hline
$3$-rd   moment  & $2.03 \times 10^{-7}$ & $1.92 \times 10^{-4}$ & $4.33 \times 10^{-4}$ & $5.86 \times 10^{-4}$ & $6.68 \times 10^{-4}$ \\ \hline
\end{tabular}
\caption{Relative squared errors when evaluating the analytical subspace moments using the numerical quadrature during the reconstruction process in Section~\ref{sec:truncationlevels}, for different truncation limits $L$. Accuracies are also affected by storing precomputed arrays using single-precision floating-point numbers.} 
\label{tab:integral_errors}
\end{table}

In~\cref{fig:effect_of_truncation}, we show the reconstructed resolutions compared to the expanded ground truth,  and the running times for various stages of the sequential moment matching. As the truncation limit used for reconstruction increases from $L=4$ to $L=10$, we observe a steady improvement in resolutions, with the best resolution of $16.68 \A$ achieved at $L=10$. However, the computational cost increases significantly due to the growing number of parameters, which scales as  $\mathcal{O}(L^3)$.  Notably, increasing the truncation limit from $L=10$ to $L=12$ results in a reduction of resolution to $33\A$,   likely due to the insufficient number of quadrature nodes used in evaluating the third-order subspace moment \eqref{eqn:Mc3}, constrained by available memory.   Thus, the optimal truncation limit for reconstruction is $L=8$ or $10.$

\begin{figure}
    \centering
    \begin{tikzpicture}
        \node[inner sep=0] at (0,0) {\includegraphics[width=\textwidth]{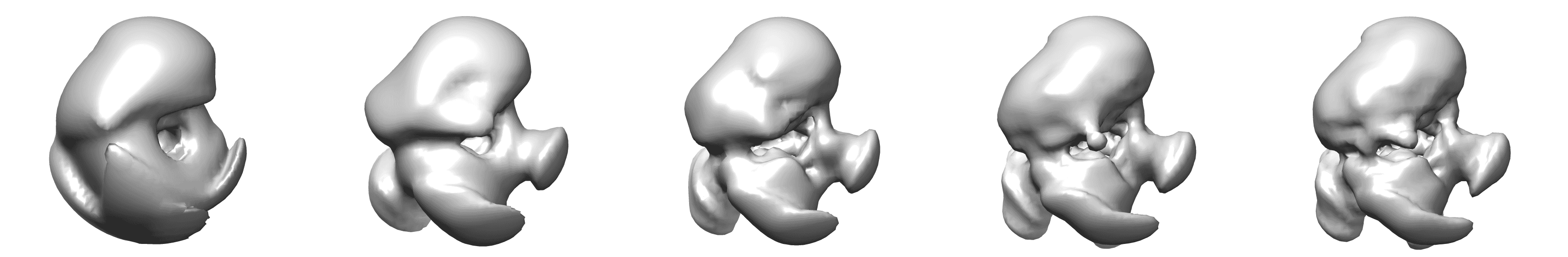}}; 
        \node[inner sep=0] at (0,-3) {\includegraphics[width=\textwidth]{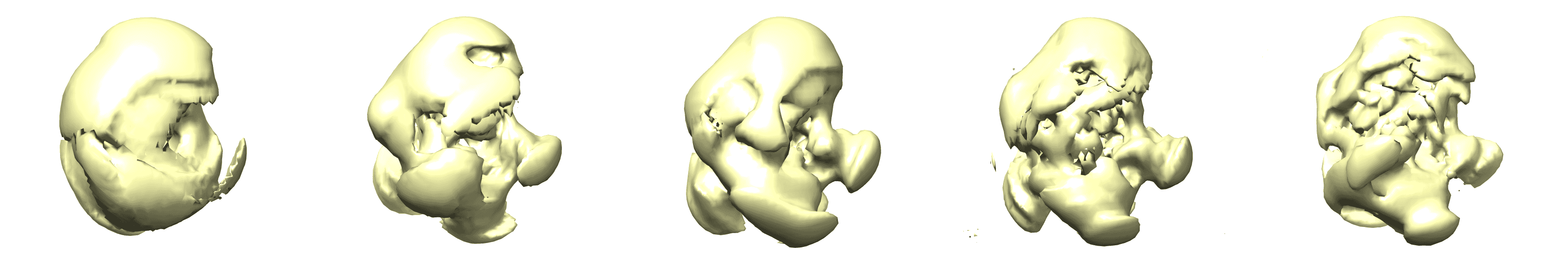}}; 
        \node[inner sep=0] at (-6.25,-1.5) {\footnotesize $24.74 \A$ };
        \node[inner sep=0] at (-3.25,-1.5) {\footnotesize $11.25 \A$ };
        \node[inner sep=0] at (0,-1.5) {\footnotesize $6.37 \A$ };
        \node[inner sep=0] at (3.25,-1.5) {\footnotesize $6.37 \A$ };
        \node[inner sep=0] at (6.25,-1.5) {\footnotesize $6.37 \A$ };

        \node[inner sep=0] at (-6.25,-4.5) {\footnotesize $41.52 \A$ };
        \node[inner sep=0] at (-3.25,-4.5) {\footnotesize $29.37 \A$ };
        \node[inner sep=0] at (0,-4.5) {\footnotesize $17.71 \A$ };
        \node[inner sep=0] at (3.25,-4.5) {\footnotesize $16.68 \A$ };
        \node[inner sep=0] at (6.25,-4.5) {\footnotesize $33.04 \A$ };
    \end{tikzpicture}
    \caption{The first row shows the expansion  (gray) of EMD-34948 data into the spherical Bessel basis with truncation limits $L=4,6,8,10,12$. The second row shows the reconstructed volumes (yellow) from the analytical subspace moments of the expanded EMD-34948 at $L=12$, which are expanded with truncation limits $L=4,6,8,10,12$, respectively. Annotated values indicate the resolution of each volume compared to the expanded ground truth at $L = 12$. Note that $6.37\A$ is the Nyquist-limited resolution.}
    \label{fig:effect_of_truncation_VOL}
\end{figure}

In~\cref{fig:effect_of_truncation_VOL}, we visualize the final reconstructed volumes  and the  spherical Bessel expansions of the downsampled EMD-34948 at different truncation limits $L$. Despite the resolution degeneration at $L=12$, the reconstructions remain visually consistent with the spherical Bessel approximations.

\subsection{Reconstruction from synthetic datasets}
\label{sec:recoverEstimatedMoMs}
In this section, we perform reconstructions using estimated subspace moments computed from synthetic noisy datasets.   In Section~\ref{sec:test_snrs}, we demonstrate the robustness of SubspaceMoM to noise by considering an idealized setting in which images are centered and corrupted only by additive white noise.  In Section~\ref{sec:ctfs}, we consider a more realistic setting in which images are additionally affected by CTFs.  In the same section, we compare the performance of SubspaceMoM with the \textit{ab initio} method implemented in RELION 5.0 via stochastic gradient descent (SGD). Our method currently assumes centered images so we only present in Section~\ref{sec:translation} a simple example illustrating its performance on uncentered images.

\begin{figure}[!ht]
    \centering
    \begin{tikzpicture}
    \node[inner sep=0] at (0,0)
{\includegraphics[width=0.2\textwidth]{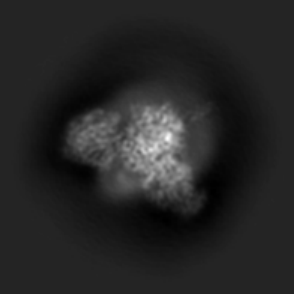}}; 
    \node[inner sep=0] at (4,0)
{\includegraphics[width=0.2\textwidth]{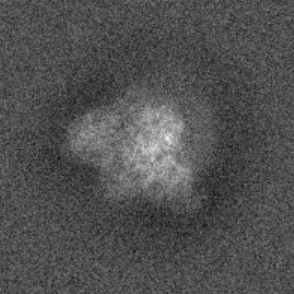}}; 
    \node[inner sep=0] at (8,0)
{\includegraphics[width=0.2\textwidth]{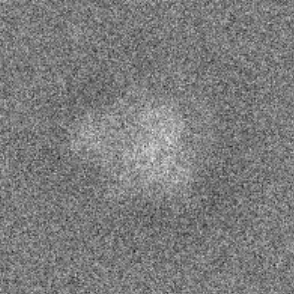}};   \node[inner sep=0] at (12,0)
{\includegraphics[width=0.2\textwidth]{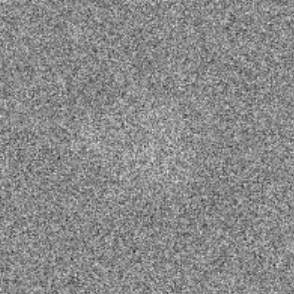}}; 
    
    \end{tikzpicture}
    \caption{Examples of synthetic image data used in  Section~\ref{sec:test_snrs}. The first image is a clean projection of EMD-34948. The second to last are images having the same orientation, but contaminated by white noise  with SNR approximately  $1,0.1$ and $0.01$, respectively.}
    \label{fig:synthetic_data}
\end{figure}

\subsubsection{Reconstruction from noisy images without CTF corruptions}
\label{sec:test_snrs}
We use the original EMD-34948  as the ground truth for generating synthetic projection images. The viewing directions in those  projections are sampled from a non-uniform distribution invariant under in-plane rotation and reflection and expanded in spherical harmonics up to $P=6$. We generate datasets with varying sample sizes $N \in \{ 2\times 10^4, 2\times 10^5, 2\times 10^6\}$ and varying SNR $\in \{1,0.1,0.01\}$, resulting in a total of $9$ datasets.  Here, the SNR of the noisy image  $I$ is  defined as $\Vert I-\epsilon \Vert_F^2/\Vert \epsilon \Vert_F^2$  where $\epsilon$ denotes the 2-D Gaussian noise.  For each target SNR level, we estimate the noise variance $\sigma^2$ by
\begin{align}
\label{eqn:variance}
   \sigma^2 = \frac{\sum_{j=1}^{1000} \|I_j^0\|_F^2}{\mathrm{SNR}\cdot\sum_{j=1}^{1000} \|\epsilon_j\|_F^2}
\end{align}
where $I_j^0$'s are randomly sampled clean projections and $\epsilon_j$'s are i.i.d. sampled white Gaussian noises with unit variance. The effect of varying SNR on images is illustrated in~\cref{fig:synthetic_data}.  

\begin{figure}[!ht]
    \centering
    \begin{tikzpicture}
    \node[inner sep=0] at (0,0)
    {\includegraphics[width=0.4\textwidth]{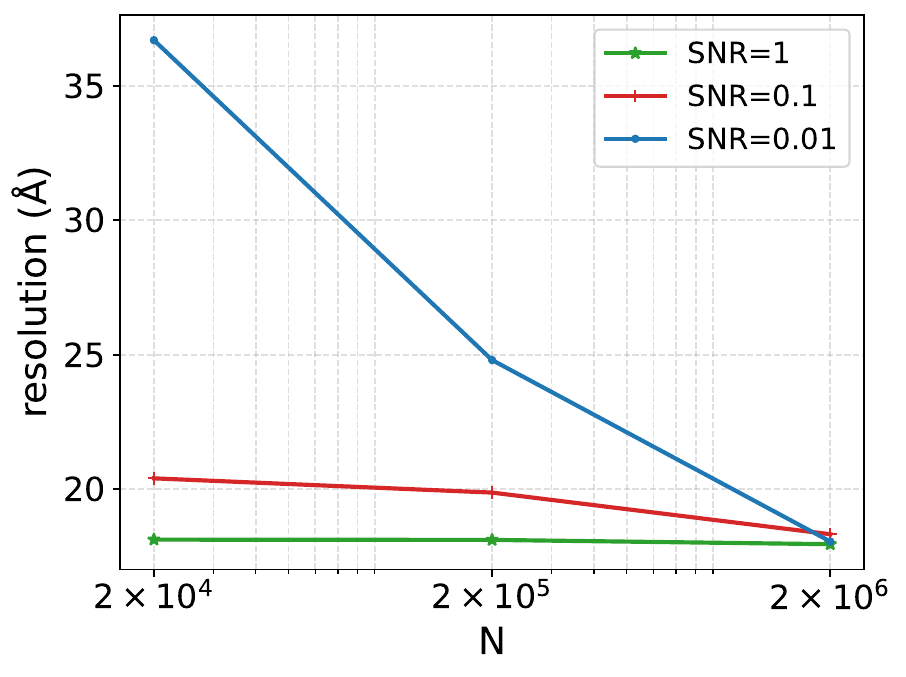}}; 
    \node[inner sep=0] at (8,0)
    {\includegraphics[width=0.4\textwidth]{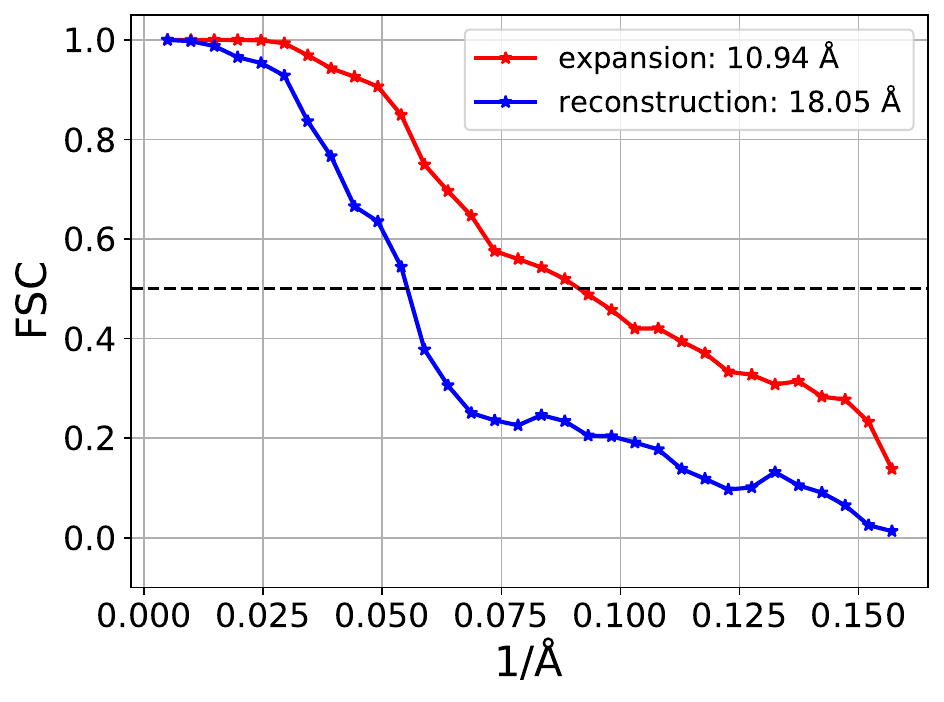}}; 
    \node[inner sep=0] at (0.5,-2.9) {\footnotesize (A)};
    \node[inner sep=0] at (8.5,-2.9) {\footnotesize (B)};
    \end{tikzpicture}
    \caption{(A) The final reconstructed resolutions compared to the ground truth EMD-34948 obtained from the synthetic datasets  generated with different sample sizes $N$ and $\mathrm{SNR}$ in Section~\ref{sec:test_snrs}.  (B) The blue line represents the FSC curve of the final reconstructed volume for $N=2\times 10^6$ and $\mathrm{SNR}\approx 0.01$, and the red curve corresponds to the  expansion of the ground truth volume at $L=10$. Both are compared to the ground truth EMD-34948.}
    \label{fig:m3_res}
\end{figure}

To make the computation feasible on our computer, all $196\times 196$ images are downsampled to  $64 \times 64$ pixels via Fourier cropping. We apply the randomized range-finding algorithm described in Section~\ref{sec:findingsubspaces} to compute the subspace moments, using a sampling size of  $s=250$ for the Gaussian sketch operator.  The debiasing procedure described in Appendix~\ref{sec:debias} is applied during moment formation.  Threshold parameters are fixed as  $\tau^{(2)}=10^{-8}, \tau^{(3)}=10^{-6}$ as in Section~\ref{sec:truncationlevels}. However, due to slower singular value decay in the presence of noise, we need to impose upper bounds on the subspace dimensions to prevent excessive memory usage: $r_2\le 220$ for the second subspace moment and  $r_3\le 120$ for the third subspace moment. In all 9 experiments, the maximum allowable ranks $r_2$ and $r_3$ are reached. The computational time required to form the subspace moments scales linearly with the number of images $N$, excluding the time needed to generate raw images. On our computer, the time for computing the subspace moments is 2 hours for $N = 2\times 10^6$.

\begin{figure}[!ht]
    \centering
\includegraphics[width=0.8\textwidth]{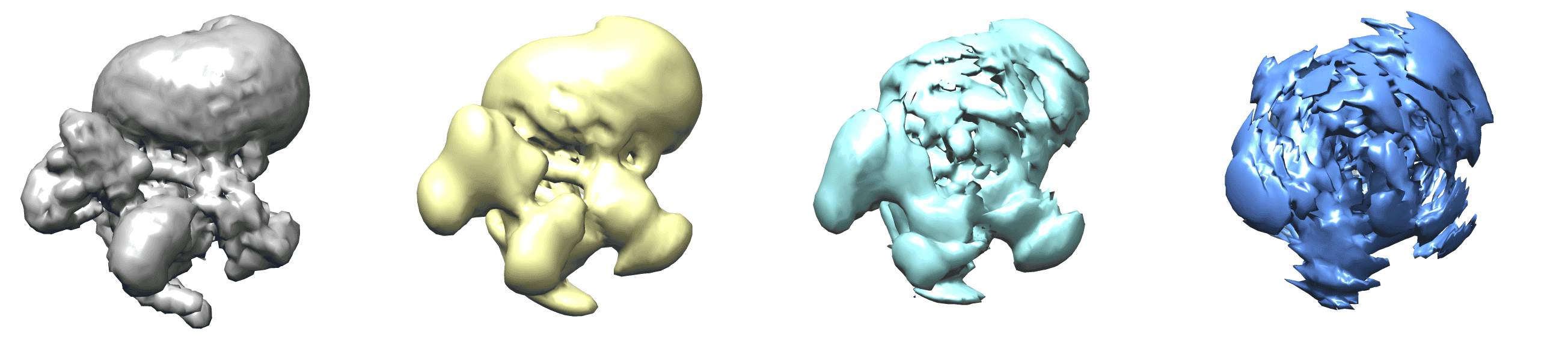} 
    \caption{From left to right: The ground truth volume (gray) used in Section~\ref{sec:test_snrs}, the expansion of the ground truth volume with truncation limit $L=10$ (yellow) with resolution $10.94 \A$, the reconstructed volume with  $L=10$ obtained by sequential moment matching  (sky blue) with resolution $18.05 \A$, the reconstructed volume  with  $L=10$  obtained by directly optimizing \eqref{eqn:subspaceMoM} (blue) with resolution $42.11 \A$. The two reconstructed volumes are obtained at  $N=2\times 10^6$ and $\mathrm{SNR} \approx 0.01$. }
    \label{fig:image_test_result}
\end{figure}

For the reconstruction, we expand the volume using $L=10$ and the viewing direction density using $P=6$, resulting in $|\cB_V| = 3432$ and $|\cB_\mu| = 28$ parameters, respectively.  We use a quadrature rule with $2496$ nodes, which is the maximum number of nodes that the memory can accommodate for the existing parameters.  The precomputation step takes about 1 hour and the sequential moment matching takes  $9$ to $12$ hours to complete.  On average, it takes about $13$ hours  to finish the full reconstruction process when $N=2\times 10^6$.

\begin{figure}[!ht]
    \centering
    \begin{tikzpicture}
        \node[inner sep=0] at (0,0)
            {\includegraphics[width=0.3\textwidth]{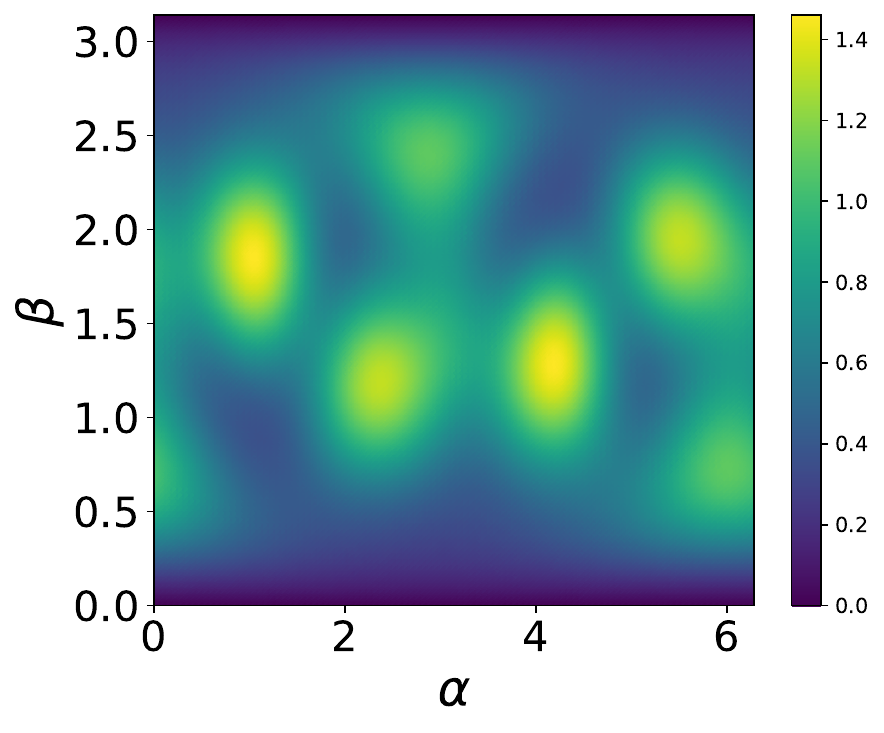}};
        \node[inner sep=0] at (5,0)
            {\includegraphics[width=0.3\textwidth]{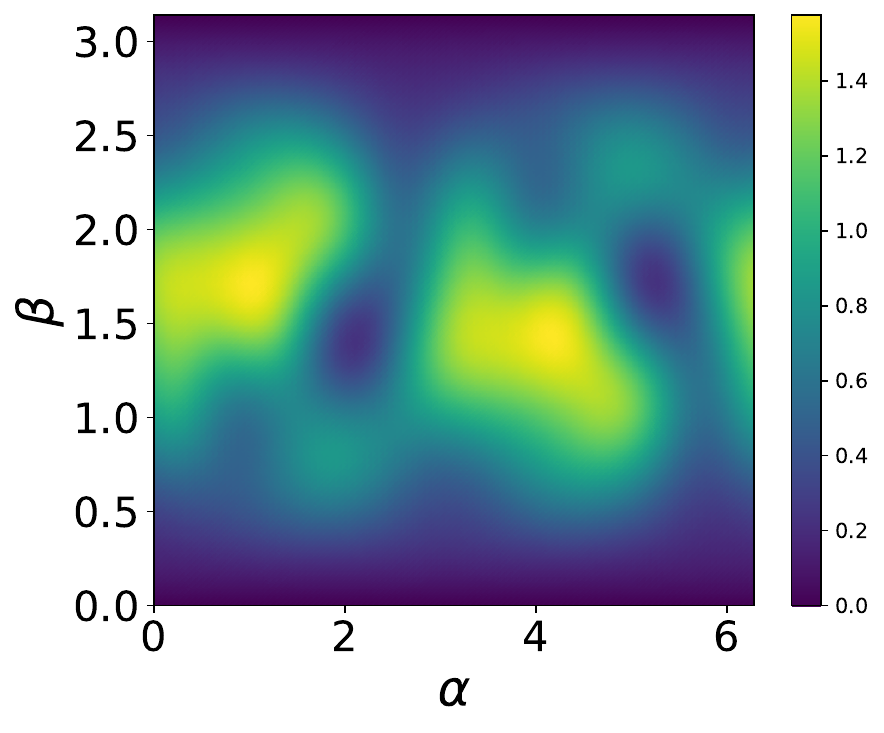}};
        \node[inner sep=0] at (10,0)
            {\includegraphics[width=0.3\textwidth]{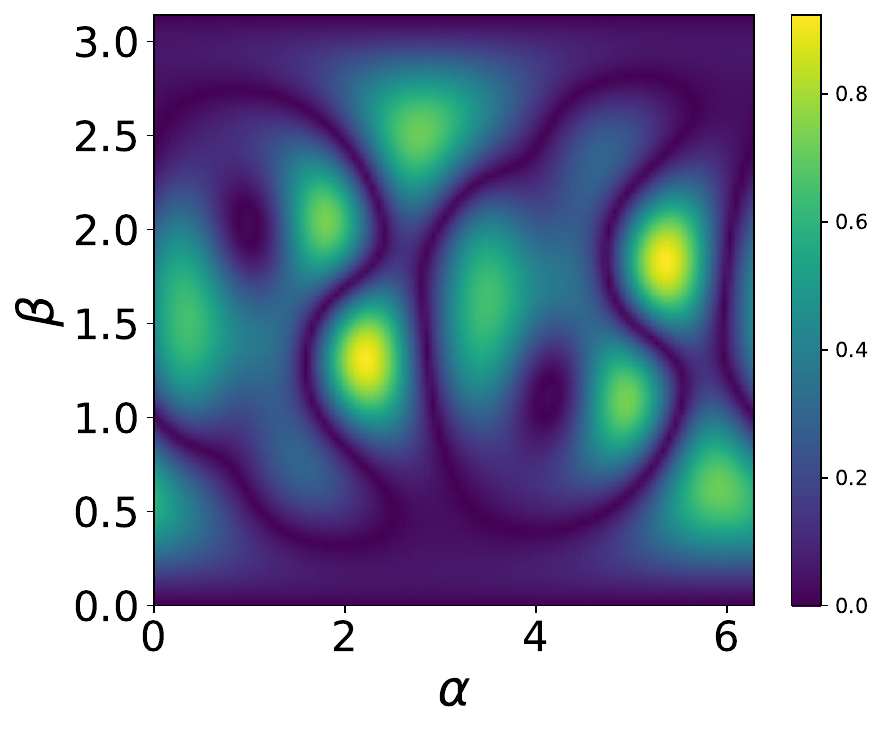}};
    \end{tikzpicture}
    \caption{The leftmost is the ground truth viewing direction density used in Section~\ref{sec:test_snrs}. The middle one is the estimated viewing direction density obtained at $N=2\times 10^6$ and $\mathrm{SNR} \approx 0.01$.  The rightmost shows their absolute point-wise differences.  They are shown as functions of the first two Euler angles $\alpha$ and $\beta$.}
    \label{fig:viewing_recons}
\end{figure}

We show the resolutions of all  final reconstructed volumes compared to the EMD-34948 in~\cref{fig:m3_res}(A), from which we can see that when the sample size is sufficiently large  $(N= 2\times 10^6)$, the final reconstructed resolution is around  $18 \A$ for all $\mathrm{SNR}$, demonstrating the robustness of the MoM approach to noise.  Notably, at $N=2\times 10^6$, the reconstructed volume at SNR $=0.1$ achieves a resolution of $18.32 \A$,which is slightly worse than that in the noisier case of SNR $=0.01$, with a resolution of $18.05 \A$. This counterintuitive outcome is possibly due to errors introduced during the alignment procedure. Nevertheless, the difference is sufficiently small that both reconstructions can be regarded as having the same quality. In~\cref{fig:m3_res}(B), we show the FSC curve of the reconstruction obtained at  $(N=2\times 10^6,\text{SNR} \approx 0.01)$ along with the FSC curve of the spherical Bessel expansion of the downsampled EMD-34948 at $L=10$. A noticeable gap remains between the reconstruction and the expansion, which is due to various approximation errors in our current method. We visualize the ground truth volume, the expansion of the ground truth and the final reconstructed volume obtained from the  sequential moment matching in~\cref{fig:image_test_result}. To demonstrate the effectiveness of the sequential moment matching, we also include, in the same figure, a visualization of the reconstructed volume obtained by directly optimizing the objective function in~\eqref{eqn:subspaceMoM}. This direct approach yields a significantly lower resolution of $42.11 \A$.  From this visualization, we can see that  although our reconstruction has a visual difference from the expansion and has noisy artifacts, it still approximates the ground truth reasonably well despite observing that there is a large discrepancy between the estimated viewing  direction density and the ground truth viewing  direction  density as shown in~\cref{fig:viewing_recons}.

\subsubsection{Reconstruction from images with CTF corruptions}
\label{sec:ctfs}

\begin{figure}[!ht]
    \centering
    \begin{tikzpicture}
    \node[inner sep=0] at (0,0)
    {\includegraphics[width=0.2\textwidth]{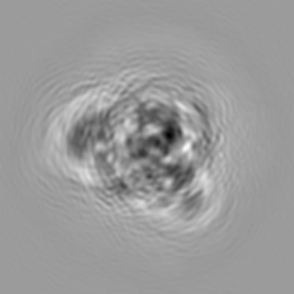}};  
     \node[inner sep=0] at (4,0)
    {\includegraphics[width=0.2\textwidth]{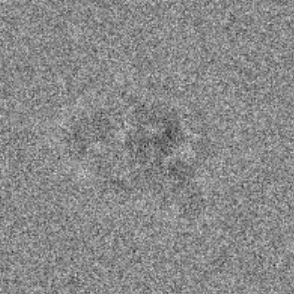}};  
     \node[inner sep=0] at (8,0)
    {\includegraphics[width=0.2\textwidth]{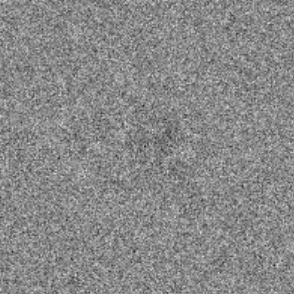}};    
    \node[inner sep=0] at (12,0)
    {\includegraphics[width=0.2\textwidth]{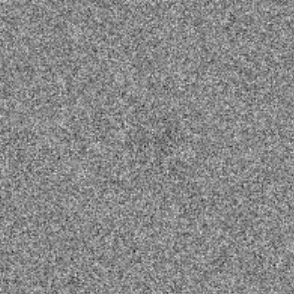}};   
    \end{tikzpicture}
    \caption{Examples of synthetic image data used in  Section~\ref{sec:ctfs}. The first image is a clean projection of EMD-34948 convolved with a radial CTF corresponding to a $300$ keV voltage and  $1$ $\mu$m defocus. The remaining images, having the same orientation and CTF, are contaminated by white noise with SNR approximately $1/16$, $1/64$, and $1/256$, respectively.}
    \label{fig:ctf_images}
\end{figure}

We generate $N=10^5$  projections of size $196\times 196$ from EMD-34948 using the ground truth viewing direction distribution in Section~\ref{sec:test_snrs}. The projections are corrupted by radial CTFs (introduced in Appendix~\ref{sec:ctf_intro}). The simulated CTF effects are implemented using the ASPIRE software~\cite{aspire}.   We split the data into  $100$ de-focus groups with  de-focus values ranging  from  1 $\mu$m to 3 $\mu$m. For all CTFs, we set the voltage as 300 keV and the spherical aberration as 2 mm.  After convolving the images with CTFs, we add Gaussian white noise to create three synthetic datasets with approximate SNR $1/16, 1/64$ and $1/256$, respectively.
The variance of the white noise is computed by~\eqref{eqn:variance} with $I_j^0$ being a clean sampled projection convolved by a randomly sampled CTF.  The effects of the CTF and white noise on the images are illustrated in~\cref{fig:ctf_images}.

\begin{figure}
    \centering
    \begin{tikzpicture}
        \node[inner sep=0] at (0,0) {\includegraphics[width=0.35\textwidth]{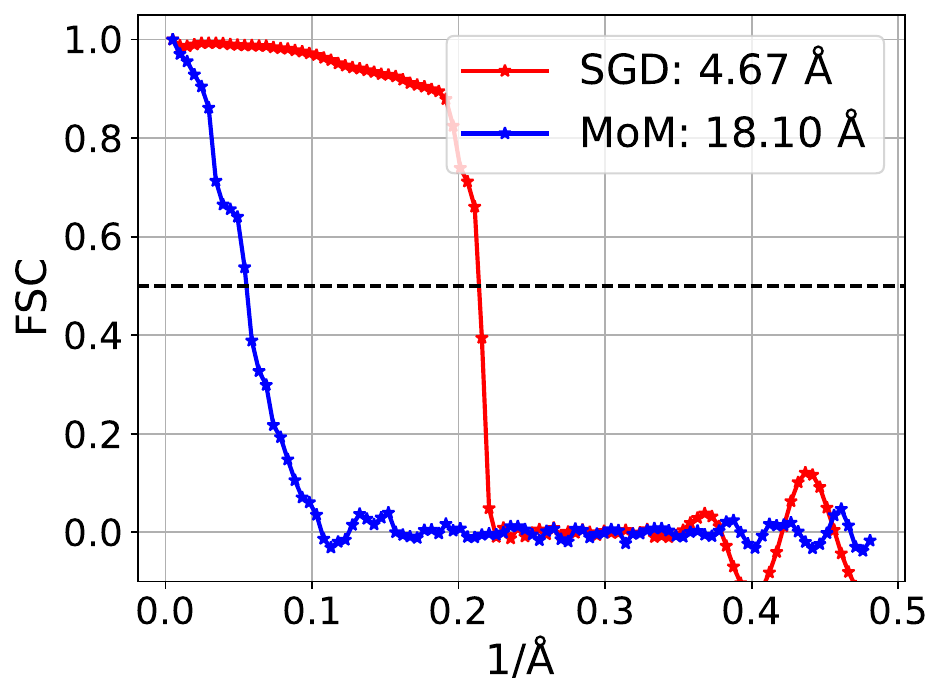}}; 
        \node[inner sep=0] at (4.5,0) {\includegraphics[width=0.22\textwidth]{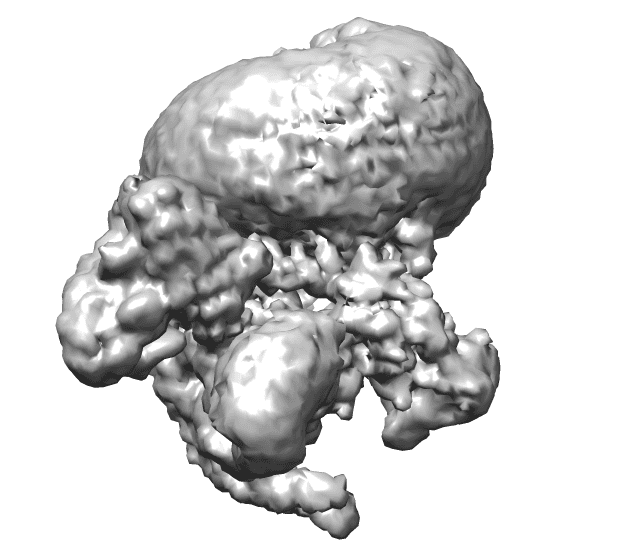}}; 
        \node[inner sep=0] at (7.5,0) {\includegraphics[width=0.22\textwidth]{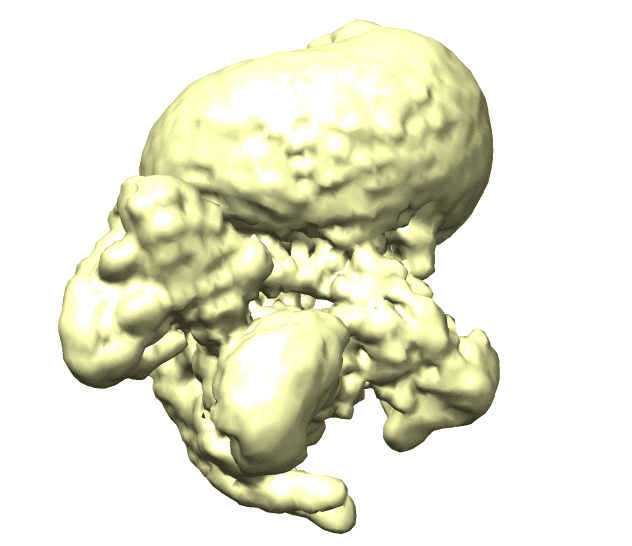}}; 
        \node[inner sep=0] at (10.5,0) {\includegraphics[width=0.22\textwidth]{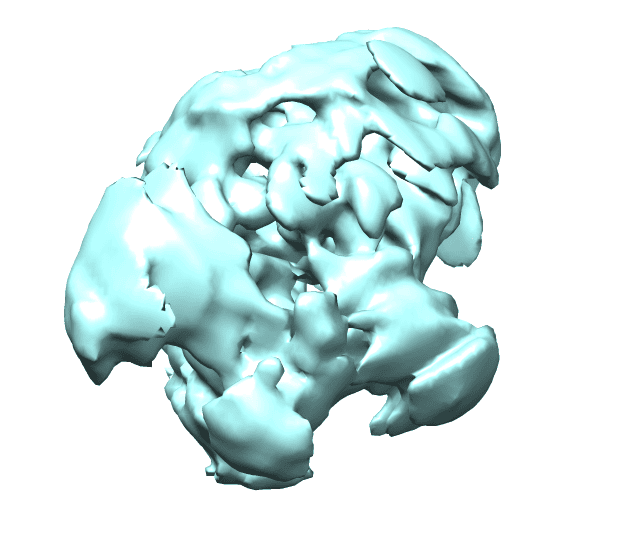}}; 
        \node[inner sep=0] at (0,-4) {\includegraphics[width=0.35\textwidth]{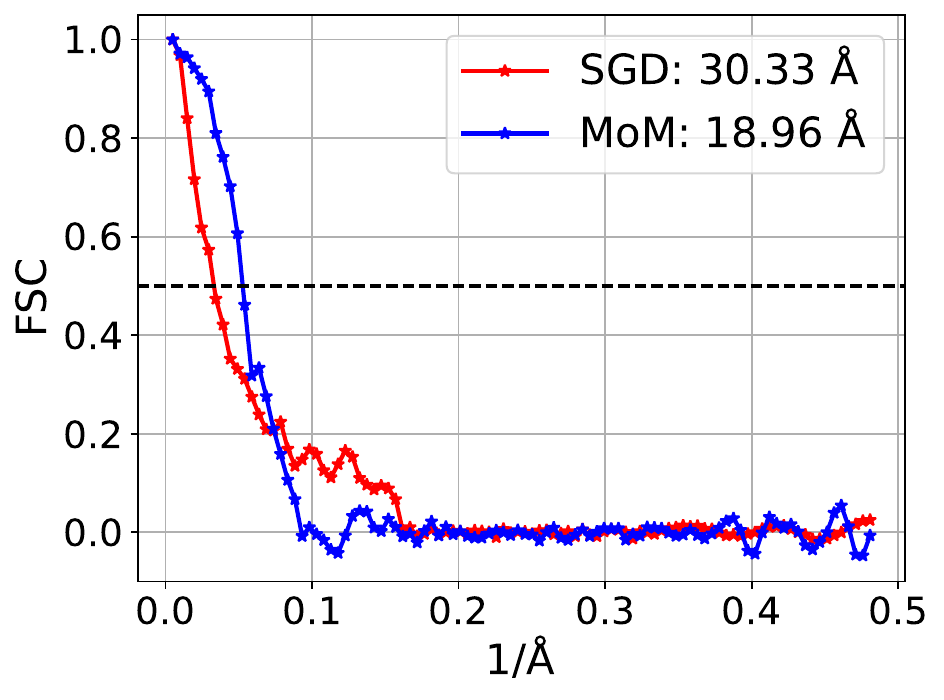}}; 
        \node[inner sep=0] at (4.5,-4) {\includegraphics[width=0.22\textwidth]{figures/vol_gt_relion.png}}; 
        \node[inner sep=0] at (7.5,-4) {\includegraphics[width=0.22\textwidth]{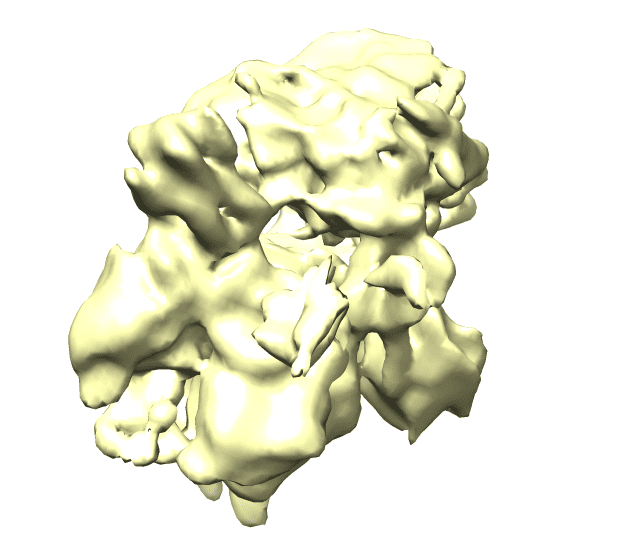}}; 
        \node[inner sep=0] at (10.5,-4) {\includegraphics[width=0.22\textwidth]{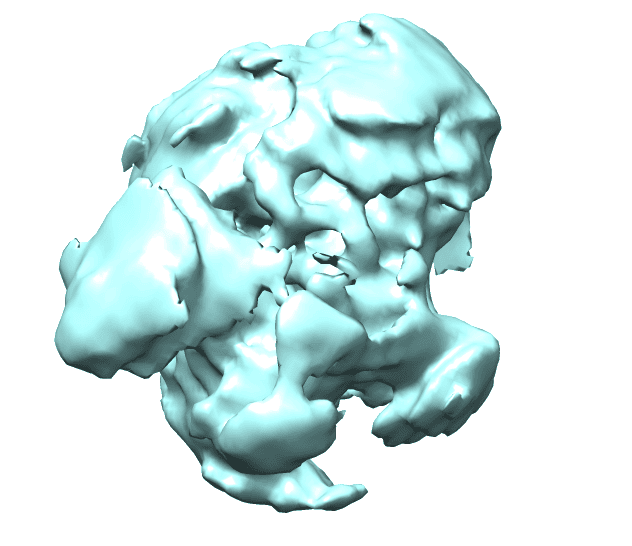}}; 
        \node[inner sep=0] at (0,-8) {\includegraphics[width=0.35\textwidth]{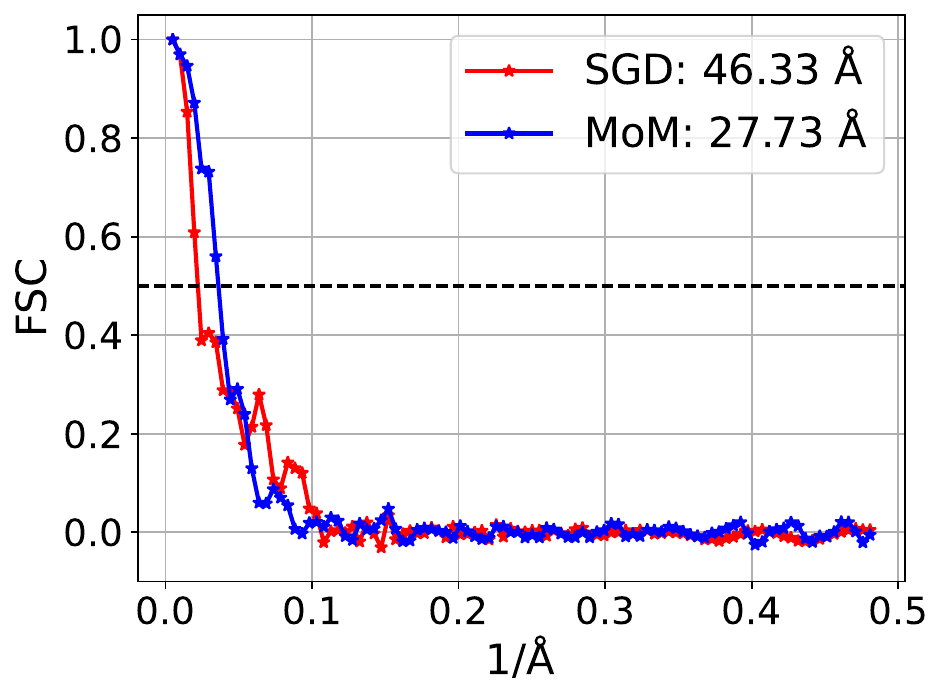}}; 
        \node[inner sep=0] at (4.5,-8) {\includegraphics[width=0.22\textwidth]{figures/vol_gt_relion.png}}; 
        \node[inner sep=0] at (7.5,-8) {\includegraphics[width=0.22\textwidth]{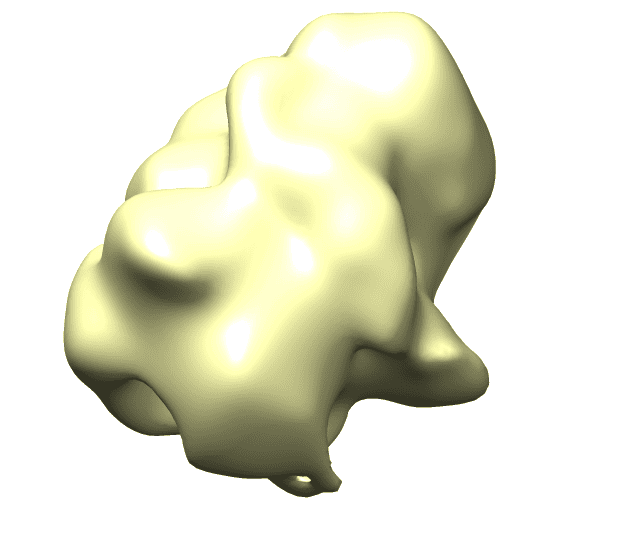}}; 
        \node[inner sep=0] at (10.5,-8) {\includegraphics[width=0.22\textwidth]{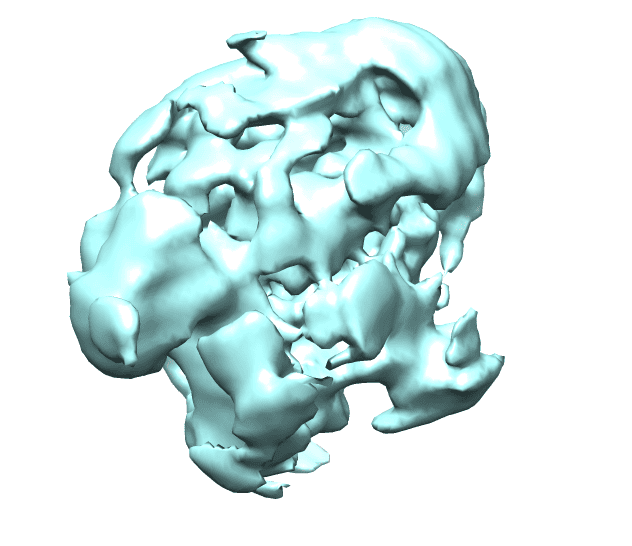}}; 
    \end{tikzpicture}
    \caption{Visualization of the reconstruction results from $10^5$ noisy images affected by CTFs in Section~\ref{sec:ctfs}. From left to right: the FSC curves obtained by the SGD \textit{ab initio} method in RELION and our SubspaceMoM compared to the ground truth volume, the ground truth volume (gray), the reconstructed volume by the SGD-based \emph{ab initio} method (yellow),  the reconstructed volume by SubspaceMoM (sky blue). We show the results obtained at SNR $\approx 1/16,1/64$ and $1/256$ from the first to the third rows. }
    \label{fig:comparisons}
\end{figure}

To obtain reference models, we run the SGD-based \emph{ab initio} reconstruction implemented in RELION 5.0~\cite{https://doi.org/10.1002/2211-5463.13873} on the synthetic datasets. We provide RELION with the ground-truth CTF information and set the offset range to zero, as the images are perfectly centered. The number of iterations is fixed at 100, as no significant changes are observed beyond this point.  For the dataset with SNR $1/16$, where the particle’s shape remains visible, the SGD method successfully produces a high-resolution reconstruction of $4.66 \A$. For the  datasets with SNR $1/64$ and $1/256$, where orientation estimation is considerably more challenging, the reconstructed resolutions are $30.33 \A$ and $46.32 \A$, respectively. 

\begin{figure}
    \centering
    \includegraphics[width=0.75\textwidth]{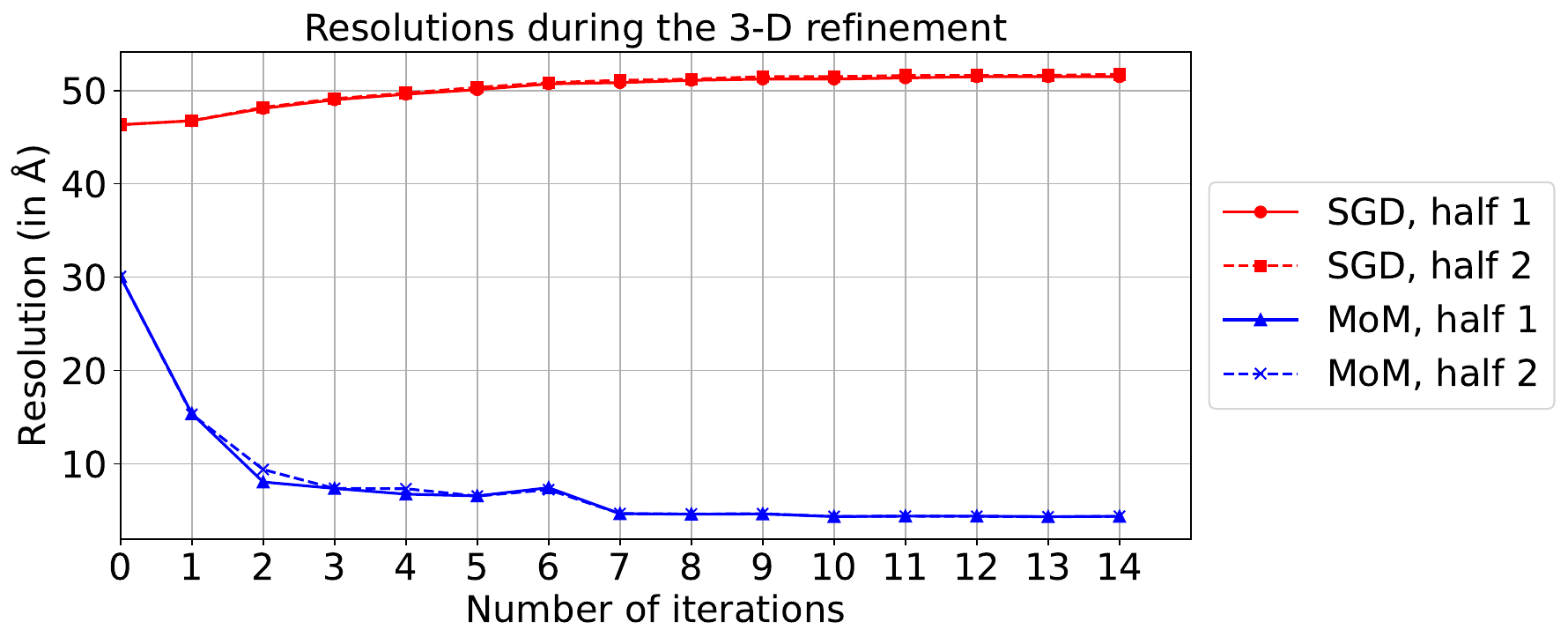}
    \caption{Resolutions progression during 3-D refinement on the two halves of the dataset with SNR = $1/256$.   The red curve corresponds to refinement initialized with the \emph{ab initio} model from the SGD-based method, while the blue curve corresponds to refinement initialized with the model from the SubspaceMoM.
}
    \label{fig:refinement_resolutions}
\end{figure}


\begin{figure}
    \centering
    \begin{tikzpicture}
        \node[inner sep=0] at (0,0) {\includegraphics[width=0.2\textwidth]{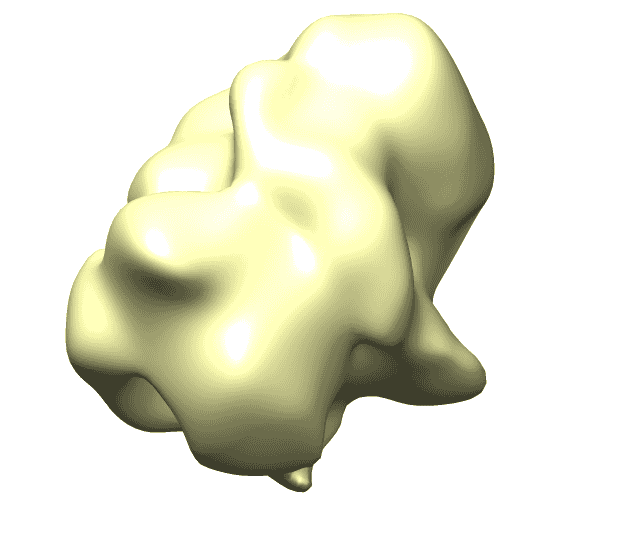}}; 
        \node[inner sep=0] at (3,0) {\includegraphics[width=0.2\textwidth]{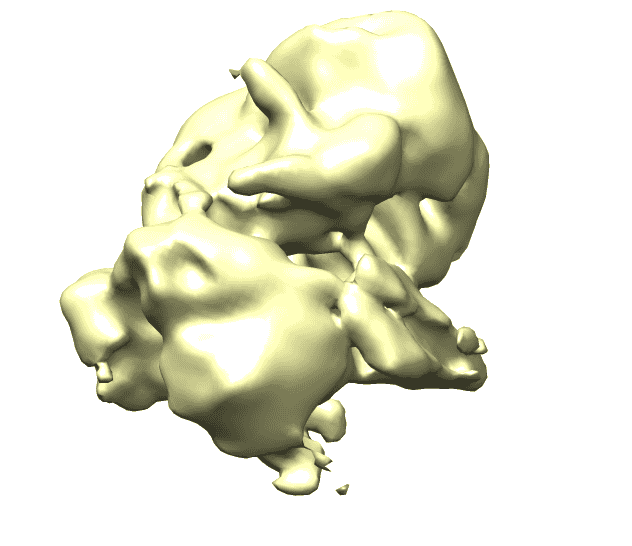}}; 
        \node[inner sep=0] at (6,0) {\includegraphics[width=0.2\textwidth]{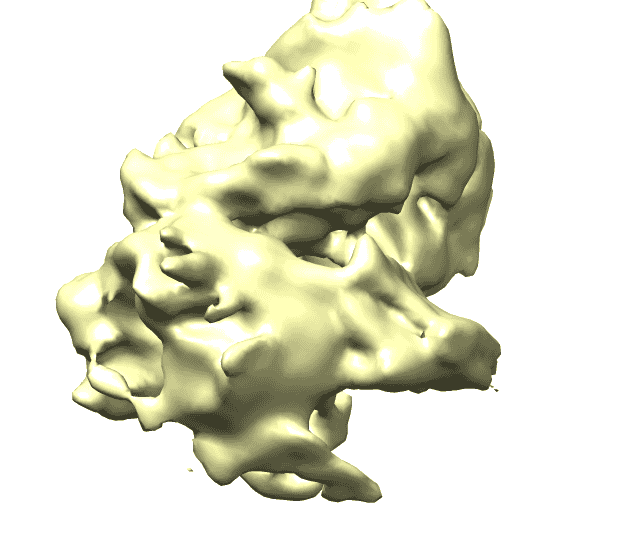}}; 
        \node[inner sep=0] at (9,0) {\includegraphics[width=0.2\textwidth]{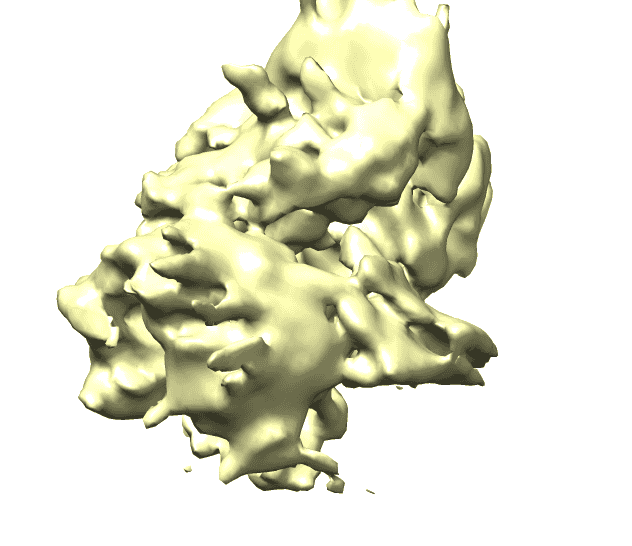}}; 
        \node[inner sep=0] at (12,0) {\includegraphics[width=0.2\textwidth]{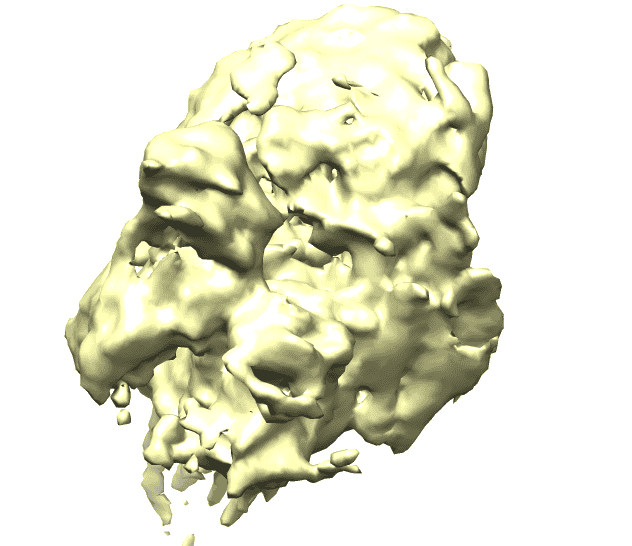}}; 

         \node[inner sep=0] at (0,-3) {\includegraphics[width=0.2\textwidth]{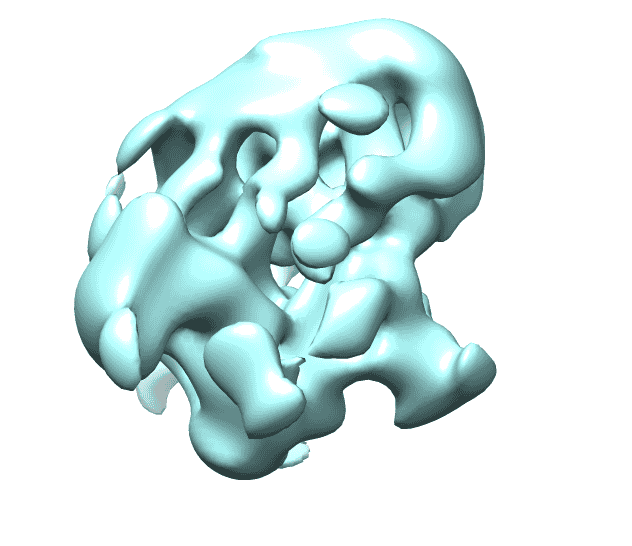}}; 
        \node[inner sep=0] at (3,-3) {\includegraphics[width=0.2\textwidth]{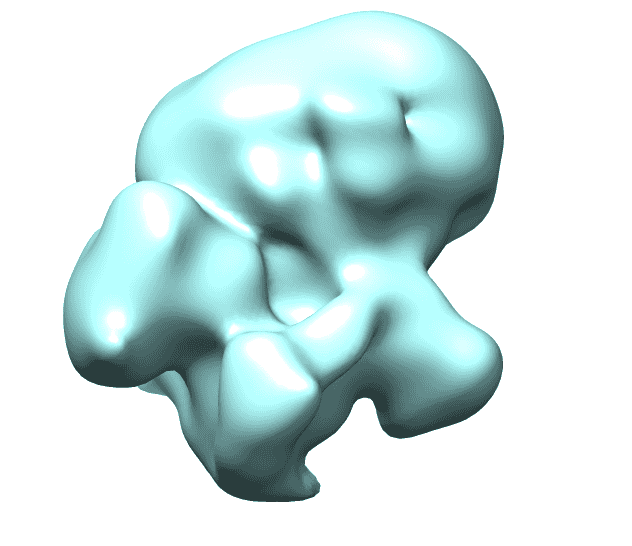}}; 
        \node[inner sep=0] at (6,-3) {\includegraphics[width=0.2\textwidth]{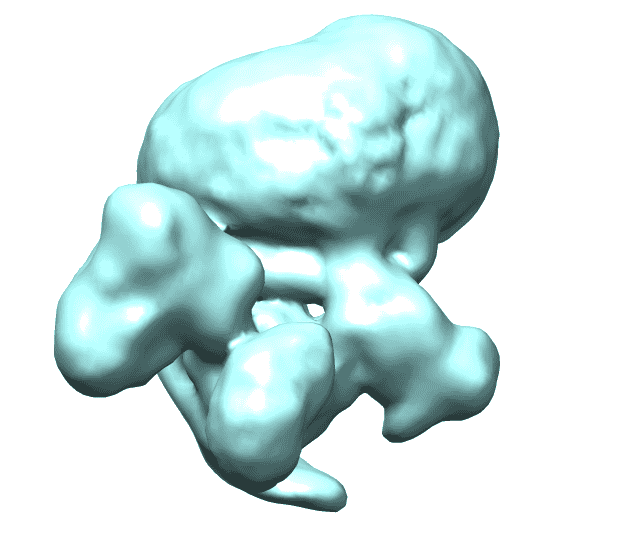}}; 
        \node[inner sep=0] at (9,-3) {\includegraphics[width=0.2\textwidth]{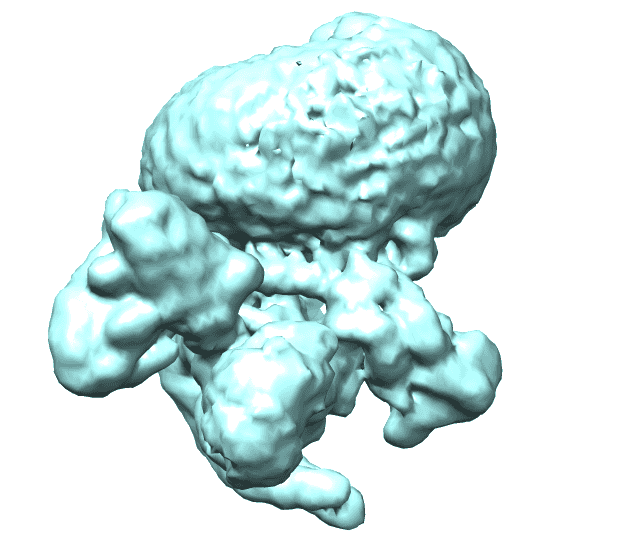}}; 
        \node[inner sep=0] at (12,-3) {\includegraphics[width=0.2\textwidth]{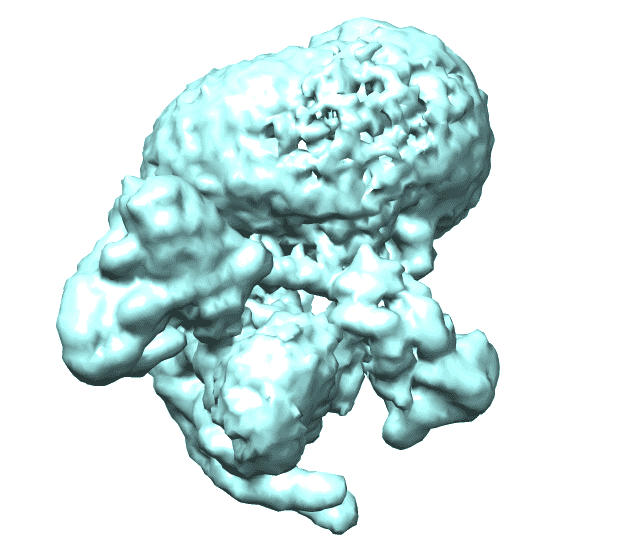}}; 

        \node[inner sep=0] at (0,-1.5) {\footnotesize Iter 0, $46.32 \A$};
        \node[inner sep=0] at (3,-1.5) {\footnotesize Iter 3, $49.00 \A$};
        \node[inner sep=0] at (6,-1.5) {\footnotesize Iter 6, $ 50.71 \A$};
        \node[inner sep=0] at (9,-1.5) {\footnotesize Iter 12, $51.61 \A$};
        \node[inner sep=0] at (12,-1.5) {\footnotesize Final, $41.77 \A$};

        \node[inner sep=0] at (0,-4.5) {\footnotesize Iter 0, $30.06 \A$};
        \node[inner sep=0] at (3,-4.5) {\footnotesize Iter 1, $15.29 \A$};
        \node[inner sep=0] at (6,-4.5) {\footnotesize Iter 2, $8.04 \A$};
        \node[inner sep=0] at (9,-4.5) {\footnotesize Iter 12, $4.38 \A$};
        \node[inner sep=0] at (12,-4.5) {\footnotesize Final, $4.21 \A$};

    \end{tikzpicture}
    \caption{Reconstructed volumes from the first half-set during RELION’s high-resolution 3-D refinement, shown at selected iterations, along with the final averaged reconstructions.   The first row (yellow volumes) corresponds to refinement initialized with the SGD-based \textit{ab initio} model, while the second row (sky blue volumes) corresponds to refinement initialized with the SubspaceMoM \textit{ab initio} model. The volumes at  iteration 0 are the \textit{ab initio} models, masked by a low-pass filter with a cutoff of $15 \A$.} 
    \label{fig:final_models}
\end{figure}

Next, we apply SubspaceMoM to the same datasets. The subspace moments are estimated using the CUR-based approach described in Section~\ref{sec:CTFs} with the sampling rule~\eqref{eqn:sampling_rule}.  When forming the subspace moments, we downsample the images to size $36\times 36$. For the second moment, we sample the index set $\mathcal{J}\subset [36^2]$ with size $|\mathcal{J}|=400$. For the third moment, we sample the index set $\tilde \cS\subset [36^2]$ with size $|\tilde \cS|=400$ and $\tilde \cJ \subset [36^4]$ with size $|\tilde \cJ|=14400$. These sampling sets are used for all datasets. Finally, we obtain the subspace moments with dimension parameters $r_1=r_2=220,r_3=120$ for all datasets. We reconstruct the volume and the viewing direction density with truncation parameters $L=10$ and $P=4$. We use the same quadrature rule as in Section~\ref{sec:test_snrs}. Forming the subspace moments only takes about 4 minutes and  the precomputation takes about  $14$ minutes. The  optimization involving only the first subspace moment~\eqref{eqn:firststageoptm} takes less than $3$ minutes, the optimization involving the first two subspace moments~\eqref{eqn:secondstageoptm} takes less than $2$ hours and the optimization involving all three subspace  moments~\eqref{eqn:subspaceMoM} takes less than $4$ hours. Therefore, our method completes each reconstruction within $7$ hours. The achieved resolutions are  $18.10 \A, 18.96 \A$ and $27.73 \A$ for datasets with SNRs of $1/16$, $1/64$, and $1/256$, respectively. Although our method is limited to low-resolution reconstructions, it is more robust than the SGD-based method in high-noise settings, as illustrated in~\cref{fig:comparisons}.

Finally, for the dataset with SNR $1/256$, we run RELION’s high-resolution 3-D refinement using the \textit{ab initio} models obtained by the SGD-based method and SubspaceMoM. We use RELION’s default parameters with fixed half-sets, except that we set the cutoff of the initial low-pass filter  to $15 \A$. Both refinements converge in 14 iterations. The refinement initialized with the SGD-based method becomes trapped in a  local minima with  final resolution $41.77 \A$, reflecting the difficulty of orientation estimation at this noise level. In contrast, the refinement initialized with SubspaceMoM produces a model closely resembling the ground truth with final resolution $4.21 \A$. These results indicate that the low-resolution model provided by SubspaceMoM serves as a reasonable initialization that improves the robustness of the whole reconstruction process. In~\cref{fig:refinement_resolutions}, we plot the per-iteration half-map resolutions for both refinements, and in~\cref{fig:final_models} we show selected reconstructed volumes from each refinement process.

\subsubsection{Effect of imperfect centering}
\label{sec:translation}
Our current method does not include a mechanism for handling uncentered images, which is a necessary extension for applications to experimentally obtained datasets. As a preliminary study, we evaluate the performance of  SubspaceMoM  on synthetic data subjected to random translations. To isolate the effects of translations, we simplify all other aspects of the experimental setup. Specifically, we downsample the EMD-34948 volume to a $32 \times 32 \times 32$ grid with a physical pixel size of $6.37 \text{\AA}$. The downsampled volume is then expanded in the spherical Bessel basis with a truncation limit $L=5$ to serve as the ground truth. We also construct a ground truth viewing direction distribution expanded in spherical harmonics with a truncation limit $P=4$. Using these ground truths, we generate $2 \times 10^5$ clean images of size $32 \times 32$, without either additive white noise or CTF effects. Each image is then translated by a random two-dimensional Gaussian shift vector with zero mean and variance $\sigma_t^2$, where the unit corresponds to one pixel length after downsampling.

\begin{figure}[!ht]
\centering
\includegraphics[width=0.4\textwidth]{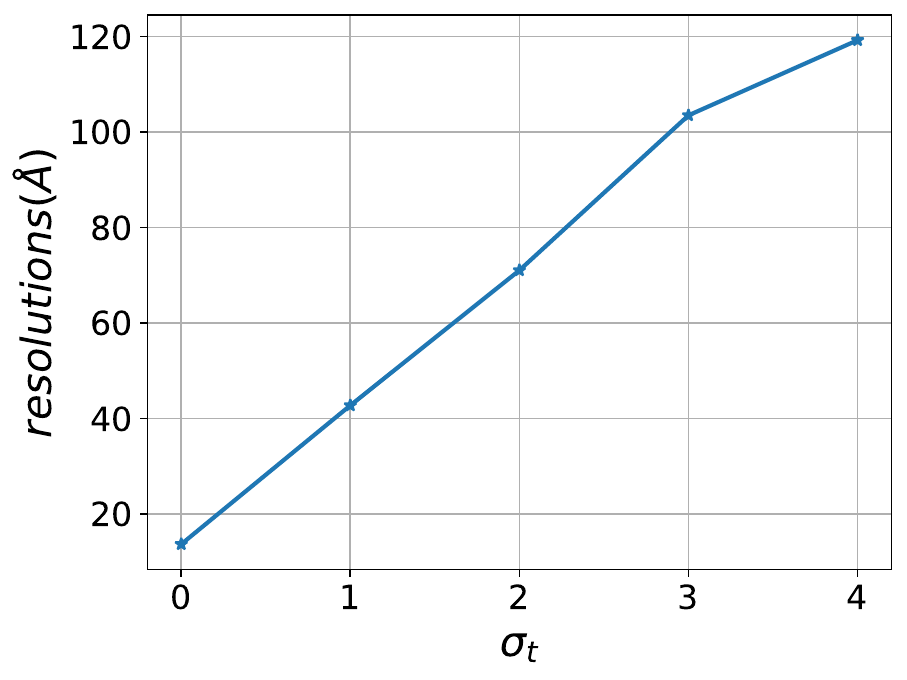}
\caption{The reconstructed resolutions compared to the expanded ground truth when images are randomly translated by two-dimensional Gaussian vectors with variance $\sigma_t^2$, where the unit corresponds to one pixel length $6.37 \A$ of a $32\times 32$ image.}
\label{fig:translation_res}
\end{figure}

To obtain the estimated subspace moments, we apply the randomized range-finding algorithm described in Section~\ref{sec:findingsubspaces} using a sample size $s = 250$ for the Gaussian sketch operators. The threshold parameters are set to $\tau^{(2)} = 10^{-10}$ and $\tau^{(3)} = 10^{-8}$. In the reconstruction process, we use high-order numerical quadratures to ensure integration accuracy up to machine precision. The resulting reconstruction problem involves 502 parameters for the volume and 14 parameters for the viewing direction distribution. The reconstructed resolutions, compared to the ground truth, are shown in~\cref{fig:translation_res}. As the figure indicates, the resolutions deteriorate significantly in the presence of translations. This suggests that further development is needed to incorporate a mechanism for handling imperfect centering within the method of moments framework prior to its application to experimentally obtained data.   We further discuss this issue in Section~\ref{sec:discussion}.

\section{Discussion}
\label{sec:discussion}

In this paper, we presented a modular and scalable approach for implementing the method of moments  for 3-D reconstruction in cryo-EM using the first three moments. To achieve computational efficiency, we applied subspace projections to compress the moments and employed a streaming strategy for projection and randomized sketching to reduce computational complexity. Additionally, we discretized the compressed representation using (inexact) numerical quadrature to further reduce complexity. Our numerical results showed that, on synthetic datasets under white noise and CTF contamination, our method  produced decent \textit{ab initio} models, using the first three subspace moments under a non-uniform viewing direction distribution. We also showed that initializing 3-D refinement with our \textit{ab initio} model leads to rapid convergence to a high-resolution reconstruction on a challenging synthetic dataset.

Our computational pipeline is flexible and can be extended to more realistic datasets. For instance, while the current approach assumes centered particle images, it can be adapted to account for random translational offsets. Since translations correspond to multiplicative phase factors in the Fourier domain, they can be incorporated into the model by modifying the forward model. Specifically, to translate an image $\cI(x,y)$ by an offset vector $\vec{t} = (t_1, t_2) \in \mathbb{R}^2$, defined by $\vec{t} \circ \cI(x,y) = \cI(x - t_1, y - t_2)$, we modify the Fourier transform of the forward model in~\eqref{eqn:fourierprojection} as
\[
\vec{t} \circ \hat \cI[ V_\star, R_j] (\xi_x,\xi_y)
= e^{i2\pi (\xi_x t_1 + \xi_y t_2)} (R_j^{-1} \circ \hat V_\star)(\xi_x,\xi_y,0).
\]
One may assume that the translation vectors are drawn from an unknown distribution satisfying mild conditions, such as radial symmetry or Gaussianity. The moments would then be formed by further integrating over the translational distribution.  Preliminary numerical results implementing this idea under a method of moments framework have been presented in~\cite{sounak2024}. It may also be possible to estimate this distribution jointly with the volume and the viewing direction distribution. Furthermore, it would be of interest to apply the techniques of this work to methods that compute higher-order moments directly from micrographs, thereby bypassing particle picking~\cite{doi:10.1137/22M1503828}.

As noted in Section~\ref{sec:experiments}, the images needed to be downsampled primarily due to memory constraints. While the computation of moments is tractable, the memory required to store intermediate quantities for efficient gradient evaluation can be substantial. Additionally, higher-resolution reconstructions would require significantly more spherical Bessel basis functions. We believe that with access to hardware acceleration and modern automatic differentiation frameworks, it will be feasible to compute gradients on the fly, which would allow larger image sizes to be used without downsampling.

Sparsity-promoting regularization, shown to improve reconstructions from second-order moments~\cite{pnas.2216507120}, can be readily incorporated into our framework. Molecular priors based on Wilson statistics~\cite{Singer:ib5103, GILLES2022106830} may further improve conditioning by reweighting the moments. There are also interesting numerical questions to be addressed, such as choosing appropriate weights in the optimization problem~\eqref{eqn:subspaceMoM}, and designing quadrature rules that balance integration error and computational cost. We leave these directions for future investigation.

Our results suggest that the method of moments can produce meaningful \textit{ab initio} models even for small molecules and low SNRs. These are regimes in which existing pipelines~\cite{SCHERES2012519, punjani2017cryosparc} often face challenges. However, particle picking remains feasible in such settings~\cite[Fig.~10f–h]{vinothkumar2016single}. Finally, the subspace method of moments presented here could also be potentially generalized to  3-D classification \emph{ab initio} modeling \cite{scheres2007disentangling}, that is, to produce a few representative 3-D classes in case the molecule exhibits structural variability. This extension may require extending the methods presented here to higher order moments (e.g., 4th order moment).

\section*{Acknowledgement}
Y. Khoo was partially supported by DMS-2111563, DMS-2339439, DE-SC0022232, and the Sloan Research Foundation. A. Singer and O. Mickelin were supported in part by AFOSR FA9550-20-1-0266 and FA9550-23-1-0249, the Simons Foundation Math+X Investigator Award, NSF DMS-2009753, and NIH/NIGMS R01GM136780-01.
Y. Wang was partially supported by funding from the University of Chicago Data Science Institute (DSI).  The authors would like to  thank Joshua Carmichael, Yunpeng Shi, Eric Verbeke, Fengyu Yang, Ruiyi Yang and Chugang Yi
for many helpful discussions. The authors would like to thank the editor and the anonymous reviewers for carefully reading the manuscript and providing many constructive suggestions. 


\appendix

\section{Debiasing the empirical (subspace) moments}
\label{sec:debias}

The noisy  Fourier images are generated from a ground truth structure $(V_\star, \mu_\star)$ according to the model
\begin{align}
\hat I_j^H = \hat H_j \odot   F \cI[V_\star,R_j] +  F \epsilon_j, \quad \epsilon_j \iid \cN(0,\sigma^2I_{d}), \quad R_j \iid \mu_\star\quad j=1,\ldots,N
\end{align}
where $N$ is the number of noisy images and $d$ is the dimension of the image in its vectorized form. The vector  $\hat H_j \in  \mathbb{C}^{d\times 1}$ encodes the CTF function and $\hat H_j=\vec 1_d$ for the case without CTF. The matrix $F \in \mathbb{C}^{d\times d}$ is the 2-D discrete Fourier transform matrix. For simplicity, we assume that $F$ is orthonormal, i.e. $FF^*=I_d$.  In this case, we also have   $FF^T=P$ where $P\in \R^{d\times d}$ is a permutation matrix.  The debiased  estimators for the second and the third order moments are given by 
\begin{align}
    \overline{M}^{(2)} & =  \left [  \sum_{j=1}^N (\hat H_j\odot \hat I_j^H) (\hat H_j\odot \hat I_j^H)^* -B^{(2)}\right ] \oslash \left [ \sum_{j=1}^N \left(\hat H_j^{\odot^2}\right) \left(\hat H_j^{\odot^2}\right)^* \right], \\
    \overline{M}^{(3)} & = \left [  \sum_{j=1}^N \left(\hat H_j\odot \hat I_j^H \right)^{\otimes 3} -B^{(3)}\right] \oslash \left [\sum_{j=1}^N \left(\hat  H_j^{\odot^2} \right)^{\otimes 3}\right]
\end{align}
where $B^{(2)} \in \C^{d\times d},B^{(3)}\C^{d\times d \times d}$ are the biases.  Let $F_k=F(:,k)$ for $k \in [d]$.  The biases can be written as 
\begin{align}
    B^{(2)}  = \sigma^2 \sum_{j=1}^N \mathrm{diag}(\hat H_j) FF^* \mathrm{diag}(\hat H_j) =\sigma^2 \sum_{j=1}^N \mathrm{diag}(\hat H_j^{\odot 2}).
\end{align}
\begin{align}
    B^{(3)} &= \sigma^2 \sum_{j=1}^N \sum_{k=1}^d (\hat H_j^2 \odot M^{(1)}) \otimes (\hat H_j\odot F_k) \otimes (\hat H_j \odot F_k)\nonumber \\
   & + \sigma^2 \sum_{j=1}^N \sum_{k=1}^d   (\hat H_j \odot F_k) \otimes (\hat H_j^2\odot M^{(1)}) \otimes (\hat H_j\odot F_k)\nonumber \\
    &+ \sigma^2 \sum_{j=1}^N \sum_{k=1}^d  (\hat H_j \odot F_k) \otimes (\hat H_j\odot  F_k) \otimes (\hat H_j^2\odot M^{(1)})  \nonumber \\
    &= 3\sigma^2 \sum_{j=1}^N \mathrm{Sym}\left[ (\hat H_j^2\odot M^{(1)}) \otimes \mathrm{diag}(\hat H_j) P \mathrm{diag}(\hat H_j) \right] 
\end{align}
where $\mathrm{Sym}:\C^{d\times d\times d}\to \C^{d\times d\times d}$ is a symmetrizing operator such that  for any $M\in \C^{d\times d\times d}$, 
\begin{align}
    \mathrm{Sym}(M)_{ijk} = \frac{1}{3}(M_{ijk}+M_{jik}+M_{jki}). 
\end{align}
In practice, we replace $M^{(1)}$ by the unbiased estimator $\overline{M}^{(1)}$.  In our streaming-based sketching implementation, these bias terms can be efficiently sketched and subtracted on the fly due to their simple  and sparse nature. We omit further implementation details.

\section{Quadrature rules for subspace moments}
\label{sec:constructSO3quadratures}
In this section, we introduce the   Euler angle parameterization of  3-D  rotations used in this paper. We can represent any rotation matrix  $R \in \mathcal{SO}(3)$ by the product of three elementary rotation matrices given by 
\begin{align*}
    R = R_z(\alpha)R_y(\beta)R_z(\gamma).
\end{align*}
Here, $\alpha \in [0,2\pi],\beta \in [0,\pi]$ and $\gamma \in [0,2\pi]$ are the Euler angles under the ZYZ convention, and $R_z,R_y$ are the elementary rotation matrices defined as 
\begin{align*}
    R_y(\theta) = \begin{pmatrix}
        \cos\theta & 0 & \sin \theta\\
        0 & 1 & 0 \\
        -\sin\theta & 0 & \cos \theta 
    \end{pmatrix},\qquad
    R_z(\phi)  = \begin{pmatrix}
        \cos\phi &  -\sin \phi & 0\\
        \sin \phi & \cos\phi & 0 \\
        0 & 0 & 1 
    \end{pmatrix}.
\end{align*}
The Haar integral of any function $f$ over $\mathcal{SO}(3)$ can be defined as  
\begin{align}
\label{eqn:SO3integraEuler}
    \Int f(R)\,{\rm d}R = \frac{1}{8\pi^2}\int_0^{2\pi} \int_0^\pi \int_0^{2\pi} f(R(\alpha,\beta,\gamma)) \, {\rm d} \alpha  \sin\beta {\rm d} \beta  {\rm d}\gamma. 
\end{align}
We assume that the function   $f$ can be expanded by Wigner D-matrices (\ref{eqn:wigner-D}) with a truncation limit $L$. To evaluate the integral in (\ref{eqn:SO3integraEuler}), we want to find a  quadrature rule that can integrate all the entries of the Wigner D-matrices
\begin{align*}
    \left \{D^{l}_{m,m'}(R) \bigg | 0\le l \le L, -l\le m,m' \le l \right \} 
\end{align*}
with respect to the Haar measure.  A particular choice is given by the Lemma 3.1 in  \cite{Grf2009SamplingSA} for constructing such a rule, which is formed by the product of two  quadrature rules. The first rule can integrate the spherical harmonics with truncation limit $L$:
\begin{align*}
    \left \{ Y^{m}_{l}(\beta,\alpha) \bigg | 0\le l \le L, -l \le m \le l \right \}. 
\end{align*}
with respect to the spherical measure ${\rm d} \alpha \, \sin\beta {\rm d} \beta $ over $[0,2\pi] \times [0,\pi]$. For this one, we use the Gaussian quadrature rules over $\mathbb{S}^2$ constructed in \cite{Grf2009SamplingSA,Grf2013EfficientAF}, which guarantee $12$ digits of accuracy. They are generously provided to download from the website via the \href{https://www-user.tu-chemnitz.de/~potts/workgroup/graef/quadrature/index.php.en}{link}. The second quadrature rule can  integrate the 1-D Fourier basis
\begin{align*}
    \left\{e^{im\gamma} \bigg | -L\le m \le L\right\}
\end{align*}
with respect to $\rm{d} \gamma$ over $[0,2\pi]$. We  use  a trapezoidal rule on $[0,2\pi]$ with $L+1$ equispaced quadrature  points. 

\begin{table}[!ht]
\small
        \centering
\begin{tabular}{|l|l|l|l|l|}
\hline
   & $P=2$   & $P=3$   & $P=4$   & $P=5$     \\ \hline
$L=4$  & $(90,378,926)$  & $(110,432,1066)$ & $(140,522,1274)$ & $(160,576,1352) $\\ \hline
$L=5$  &$ (132,638,1664) $& $(168,704,1952)$ & $(192,792,2080)$ & $(252,902,2368)$  \\ \hline
$L=6$  & $(196,936,2812) $& $(224,1066,2964)$ & $(294,1274,3382)$ &$ (336,1352,3534)$ \\ \hline
$L=7$  & $(256,1470,4092) $& $(336,1560,4620)$ & $(384,1830,4840)$ &$ (464,1950,5368)$ \\ \hline
$L=8$  & $(378,2074,6100) $& $(432,2210,6350)$ & $(522,2516,7050)$ & $(576,2652,7300)$\\ \hline
$L=9$  & $(480,2812,8176)$ & $(580,2964,9016)$ & $(640,3382,10192) $&$ (720, 3534 ,10192)$\\ \hline 
$L=10$  & $(638,3738,11284)$ & $(704,3906,12710)$ & $(792,4410,12710) $&$ (902, 4620,13082)$\\
\hline
$L=11$  & $(768,4830,14348)$ & $(864,5060,15572)$ & $(984,5612,16048) $&$ (1176, 5842,17272)$\\
\hline
$L=12$  & $(936,6100,18796)$ & $(1066,6350,19314)$ & $(1274,7050,24864) $&$ (1352, 7300,24864)$\\
\hline
\end{tabular}
\caption{The sizes of the constructed quadrature rules, which can form the first three (subspace) moments with $12$ digits of accuracy when the volume is represented by the spherical Bessel basis with a truncation limit of $L$ and the rotational density is represented by the Wigner-D matrices with a truncation limit of $P$. The Gaussian quadrature that can integrate the spherical harmonics with order $31$ is not provided on the website. Hence in the code, we use the Gaussian quadrature for $32$ as a replacement.} 
\label{tab:numquadnodes}
\end{table}

When  the volume and and the rotational density are represented by the spherical Bessel basis and the Wigner D-matrices  respectively, the moments can be fully represented by the Wigner D-matrices. Therefore, the constructed quadrature rules  in this Appendix can be applied to integrate the moments with machine precision. To see this, we use the fact that the representation of rotating the spherical harmonic $Y_l^m$ by $R^{-1}\in \mathcal{SO}(3)$ is given by the Wigner-D matrix $(D_{m,m'}^l(R))_{m,m'=-l}^l$,  through 
\begin{align}
    (R^{-1} \circ Y_l^m) (\theta,\varphi) = \sum_{m'=-l}^l \overline{D_{m,m'}^l}(R) Y_l^{m'}(\theta,\varphi)=\sum_{m'=-l}^l (-1)^{m-m'} D_{m',m}^l(R) Y_l^{m'}(\theta,\varphi)
\end{align}
Using the Fourier-slice theorem (\ref{eqn:fourierslice}), the 2-D projection of the Fourier volume is
\begin{align}
\label{eqn:imagesWignerD}
   \hat \cI[a,R] (\rho,\varphi) = \sum_{l=0}^L \sum_{s=1}^{S(l)} \sum_{m=-l}^l \sum_{m'=-l}^l a_{l,m,s} c_{l,s} j_l(\rho \frac{z_{l,s}}{1/2})   (-1)^{m-m'} Y_l^{m'}(\pi/2,\varphi) D_{m',m}^l(R) . 
\end{align}
We use another fact that the product of two Wigner D-matrices with orders $l$ and $l'$ can be linearly represented by the Wigner D-matrices with order less than or equal to $l+l'$,  through
\begin{align}
\label{eqn:prod_wigner_D}
    D_{m_1,m_2}^{l}(R)D_{m_1',m_2'}^{l'}(R) =& \sum_{p=|l-l'|}^{l+l'} \langle l,m_1,l,m_1' | p, m_1+m_1' \rangle  \langle l,m_2,l,m_2' | p, m_2+m_2' \rangle \nonumber \\ 
    & \cdot D^{p}_{m_1+m_1',m_2+m_2'} (R),
\end{align}
where $\langle \cdots , \cdots | \cdots \rangle$ denotes the Clebsch-Gordan coefficients \cite[p 351]{chirikjian2016harmonic}. The first three subspace moments are formed by outer products of terms  in the form of \eqref{eqn:imagesWignerD}, multiplied by a rotational density represented by Wigner-D matrices \eqref{eqn:viewingdensityWignerD}. Combing that with \eqref{eqn:wigner-D}, we know that they can be represented by the Wigner D-matrices with orders less than or equal to $L+P,2L+P$ and $3L+P$ respectively.   In Table \ref{tab:numquadnodes}, we list the sizes of the constructed quadrature rules used to  form the first three moments as functions of some truncation limits $L$ and $P$, with $12$ digits of accuracy.

\section{von Mises–Fisher distribution}
\label{sec:von-Mises–Fisher}
The von Mises–Fisher distribution defined on $\bS^{d-1}$  has the probability density function 
\begin{align}
\label{eqn:vmf_density}
    f(x|\mu,\kappa) = C(\kappa) \exp(\kappa \mu^T x), 
\end{align}
where $\kappa>0$ is the concentration parameter and the normalizing constant $C(\kappa)$ is given by
\begin{align*}
    C(\kappa) = \frac{\kappa^{d/2-1}}{(2\pi)^{d/2} I_{d/2-1}(\kappa)}.
\end{align*}
Here, $I_v$ denotes the  modified Bessel function of the first kind \cite{brychkov2008handbook}. In Section \ref{sec:experiments}, we use a mixture of von Mises–Fisher distributions to model the non-uniform  viewing direction distribution, whose density is then given by  
\begin{align*}
    \rho(x) = \sum_{j=1}^n w_j  f(x|\mu_j,\kappa)
\end{align*}
with  $w_j>0$ and  $\sum_{j=1}^n w_j=1$.

\section{Contrast transfer function (CTF)}
\label{sec:ctf_intro}
Let $\kappa>0$ be the radial frequency. The radial contrast transfer function is defined as
\begin{align}
    CTF(\kappa)=\sin\left(-\pi\lambda \delta \kappa^2+\frac{\pi}{2}C_s\lambda^3\kappa^4-\alpha\right)e^{-B\kappa^2/4}
\end{align}
where $\lambda$ is the electron wavelength computed from the voltage used, $\delta$ is the defocus value, $C_s$ is the spherical abberation, $\alpha$ is  amplitude contrast and $B$ is the B-factor of Gaussian blurring. We use $\alpha=0.1$ and $B=0$.

\bibliographystyle{siam} 
\bibliography{refs} 

\end{document}